\documentclass[12pt,english,refpage,intoc,bibliography=totoc,index=totoc,BCOR7.5mm,captions=tableheading]{amsart}
\usepackage{lmodern}

\usepackage[T1]{fontenc}
\usepackage[latin9]{inputenc}
\usepackage{geometry}
\usepackage{color}
\usepackage{babel}
\usepackage{textcomp}
\usepackage{enumitem}
\usepackage{amstext}
\usepackage{amsthm}
\usepackage{amssymb}
\usepackage[unicode=true,
 bookmarks=true,bookmarksnumbered=true,bookmarksopen=false,
 breaklinks=false,pdfborder={0 0 1},backref=false,colorlinks=true]
 {hyperref}
\hypersetup{pdftitle={The LyX User's Guide},
 pdfauthor={LyX Team},
 pdfsubject={LyX},
 pdfkeywords={LyX},
 linkcolor=black, citecolor=black, urlcolor=blue, filecolor=blue, pdfpagelayout=OneColumn, pdfnewwindow=true, pdfstartview=XYZ, plainpages=false}

\makeatletter
\numberwithin{equation}{section}
\numberwithin{figure}{section}
\theoremstyle{definition}
    \ifx\thechapter\undefined
      \newtheorem{example}{\protect\examplename}
    \else
      \newtheorem{example}{\protect\examplename}[chapter]
    \fi
\theoremstyle{definition}
    \ifx\thechapter\undefined
      \newtheorem{defn}{\protect\definitionname}
    \else
      \newtheorem{defn}{\protect\definitionname}[chapter]
    \fi
\theoremstyle{remark}
    \ifx\thechapter\undefined
      \newtheorem{rem}{\protect\remarkname}
    \else
      \newtheorem{rem}{\protect\remarkname}[chapter]
    \fi
\theoremstyle{plain}
    \ifx\thechapter\undefined
      \newtheorem{prop}{\protect\propositionname}
    \else
      \newtheorem{prop}{\protect\propositionname}[chapter]
    \fi
\theoremstyle{plain}
    \ifx\thechapter\undefined
	    \newtheorem{thm}{\protect\theoremname}
	  \else
      \newtheorem{thm}{\protect\theoremname}[chapter]
    \fi
\theoremstyle{plain}
    \ifx\thechapter\undefined
      \newtheorem{lem}{\protect\lemmaname}
    \else
      \newtheorem{lem}{\protect\lemmaname}[chapter]
    \fi
\theoremstyle{plain}
    \ifx\thechapter\undefined
      \newtheorem{question}{\protect\questionname}
    \else
      \newtheorem{question}{\protect\questionname}[chapter]
    \fi
\theoremstyle{plain}
    \ifx\thechapter\undefined
      \newtheorem{conjecture}{\protect\conjecturename}
    \else
      \newtheorem{conjecture}{\protect\conjecturename}[chapter]
    \fi

\usepackage[figure]{hypcap}

\let\myTOC\tableofcontents
\renewcommand\tableofcontents{%
  \frontmatter
  \pdfbookmark[1]{\contentsname}{}
  \myTOC
  \mainmatter }

\usepackage{fancyhdr}

\@ifundefined{extrarowheight}
 {\usepackage{array}}{}
\makeatother

\providecommand{\conjecturename}{Conjecture}
\providecommand{\definitionname}{Definition}
\providecommand{\examplename}{Example}
\providecommand{\lemmaname}{Lemma}
\providecommand{\propositionname}{Proposition}
\providecommand{\questionname}{Question}
\providecommand{\remarkname}{Remark}
\providecommand{\theoremname}{Theorem}

\begin{document}
\title{Infinite Series Whose Topology of Convergence Varies From Point to
Point}
\author{by Maxwell C. Siegel (University of Southern California)}
\address{3620 S. Vermont Ave., KAP 104, Los Angeles, CA 90089-2532}
\curraddr{1626 Thayer Ave., Los Angeles, CA, 90024}
\email{siegelmaxwellc@ucla.edu}
\date{4 June, 2023}
\begin{abstract}
Let $p$ and $q$ be distinct primes, and consider the expression
$S_{p,q}\left(\mathfrak{z}\right)$ defined by the formal series $\sum_{n=0}^{\infty}q^{\#_{1}\left(\left[\mathfrak{z}\right]_{2^{n}}\right)}/p^{n}$,
where $\mathfrak{z}$ is a $2$-adic integer variable, $\left[\mathfrak{z}\right]_{2^{n}}$
is the integer in $\left\{ 0,\ldots,2^{n}-1\right\} $ congruent to
$\mathfrak{z}$ mod $2^{n}$, and where, for any integer $m\geq0$,
$\#_{1}\left(m\right)$ is the number of $1$s in the binary expansion
of $m$. When $\mathfrak{z}\in\left\{ 0,1,2,\ldots\right\} $, $S_{p,q}\left(\mathfrak{z}\right)$
reduces to a geometric series with $1/p$ as its common ratio. This
series converges in the topology of $\mathbb{R}$, and its sum is
a rational number which is also a $q$-adic integer. For $\mathfrak{z}\in\mathbb{Z}_{2}\backslash\left\{ 0,1,2,\ldots\right\} $,
$\#_{1}\left(\left[\mathfrak{z}\right]_{2^{n}}\right)\rightarrow\infty$
as $n\rightarrow\infty$, and so $S_{p,q}\left(\mathfrak{z}\right)$
converges in the $q$-adic topology to a $q$-adic integer. In this
way, we can define a ``$\left(2,q\right)$-adic function'' $X_{p,q}:\mathbb{Z}_{2}\rightarrow\mathbb{Z}_{q}$
by the sum of $S_{p,q}\left(\mathfrak{z}\right)$ in $\mathbb{R}$
for $\mathfrak{z}\in\mathbb{Z}_{2}\backslash\left\{ 0,1,2,\ldots\right\} $
and as the sum of $S_{p,q}\left(\mathfrak{z}\right)$ in $\mathbb{Z}_{q}$
for all other $\mathfrak{z}$. Thus, while $X_{p,q}\left(\mathfrak{z}\right)$
is well-defined as a $q$-adic integer for all $\mathfrak{z}\in\mathbb{Z}_{2}$,
the topology of convergence used to sum the series representation
of $X_{p,q}\left(\mathfrak{z}\right)$ to compute its value at any
given $\mathfrak{z}$ depends on the value of $\mathfrak{z}$. This
represents an entirely new type of point-wise convergence, one where
the topology in which the limit of a sequence of functions $\left\{ f_{n}\right\} _{n\geq1}$
is evaluated depends on the point at which the sequence is evaluated.
In a manner comparable to the adèle ring of a number field, functions
defined by $\mathcal{F}$-series require considering different metric
completions of an underlying field in order to be properly understood.
This paper catalogues a variety of examples of this new, unstudied
convergence phenomenon, and presents the concept of a ``frame'',
a rigorous formalism for defining and studying $\mathcal{F}$-series
and their peculiar modes of convergence.
\end{abstract}

\maketitle

\section*{\label{sec:Introduction}Introduction}

Let $p$ and $q$ be indeterminates, and, for any $x\in\mathbb{R}$,
let $\mathbb{N}_{x}$ denote the set of integers $\geq x$. We adopt
the convention of defining any sum:
\[
\sum_{n=n_{\textrm{lower}}}^{n_{\textrm{upper}}}
\]
to be $0$ whenever $n_{\textrm{upper}}<n_{\textrm{lower}}$. Similarly,
we define any product:
\[
\prod_{n=n_{\textrm{lower}}}^{n_{\textrm{upper}}}
\]
to be $1$ whenever $n_{\textrm{upper}}<n_{\textrm{lower}}$. For
brevity, we write $\mathbb{Z}_{2}^{\prime}$ to denote $\mathbb{Z}_{2}\backslash\mathbb{N}_{0}$
(the set of all $2$-adic integers which are not in $\mathbb{N}_{0}$).
For any integer $m\in\mathbb{N}_{0}$, let $\#_{1}\left(m\right)$
denote the number of $1$s in the binary expansion of $m$. It is
not difficult to show that $\#_{1}$ is the \emph{unique} function
$\mathbb{N}_{0}\rightarrow\mathbb{N}_{0}$ satisfying the system of
functional equations:
\begin{align}
\#_{1}\left(2m\right) & =\#_{1}\left(m\right)\label{eq:number of 1s functional equation}\\
\#_{1}\left(2m+1\right) & =\#_{1}\left(m\right)+1
\end{align}
for all $m\in\mathbb{N}_{0}$. Next, for any $m\in\mathbb{N}_{0}$,
let $\lambda_{2}\left(m\right)$ the number of binary digits of $m$,
also known as the number of \textbf{bits }of $m$; an explicit formula
for this function is $\lambda_{2}\left(m\right)=\left\lceil \log_{2}\left(m+1\right)\right\rceil $.
Like $\#_{1}$, $\lambda_{2}$ is the unique function $\mathbb{N}_{0}\rightarrow\mathbb{N}_{0}$
satisfying the system of functional equations:
\begin{align}
\lambda_{2}\left(2m\right) & =\lambda_{2}\left(m\right)+1\label{eq:lambda_2 functional equation}\\
\lambda_{2}\left(2m+1\right) & =\lambda_{2}\left(m\right)+1
\end{align}
Finally, for any $2$-adic integer $\mathfrak{z}\in\mathbb{Z}_{2}$,
we define the formal expression $S_{p,q}\left(\mathfrak{z}\right)$
by:
\begin{equation}
S_{p,q}\left(\mathfrak{z}\right)\overset{\textrm{def}}{=}\sum_{n=0}^{\infty}\frac{q^{\#_{1}\left(\left[\mathfrak{z}\right]_{2^{n}}\right)}}{p^{n}}\label{eq:The definition of X_p,q}
\end{equation}
Here $\left[\mathfrak{z}\right]_{2^{n}}$ denotes the  projection
of $\mathfrak{z}$ mod $2^{n}$; that is, the unique integer in $\left\{ 0,\ldots,2^{n}-1\right\} $
which is congruent to $\mathfrak{z}$ mod $2^{n}$. Following \cite{my dissertation}\textemdash the
author's PhD dissertation\textemdash we will refer to $S_{p,q}\left(\mathfrak{z}\right)$
as an example of an \textbf{$\mathcal{F}$-series} (F for ``Frame'').

The paper presents an apparently novel form of series convergence
that occurs when attempting to realize formal expressions like $S_{p,q}\left(\mathfrak{z}\right)$
as bonafide functions of the $2$-adic integer variable $\mathfrak{z}$.
For psychological reasons, rather than start with an abstract definition,
we will begin with a simple example that captures the essential features
of our findings.
\begin{example}
\label{exa:1}Let $p=2$ and $q=3$, so that:
\begin{equation}
S_{2,3}\left(\mathfrak{z}\right)=\sum_{n=0}^{\infty}\frac{3^{\#_{1}\left(\left[\mathfrak{z}\right]_{2^{n}}\right)}}{2^{n}}\label{eq:S_2,3}
\end{equation}
Now, let us try to sum this series.

We recall a basic fact of $2$-adic analysis: every $2$-adic integer
$\mathfrak{z}$ can be represented as a formal sum of the shape:
\begin{equation}
\mathfrak{z}=\sum_{n=0}^{\infty}d_{n}\left(\mathfrak{z}\right)2^{n}\label{eq:Hensel series}
\end{equation}
where, for each $n$, $d_{n}\left(\mathfrak{z}\right)\in\left\{ 0,1\right\} $
is the $n$th \textbf{$2$-adic digit }of $\mathfrak{z}$. (\ref{eq:Hensel series})
is sometimes called the \textbf{Hensel series }of $\mathfrak{z}$
(see, for example, \cite{Amice}). Hensel series justify the intuitive,
heuristic view of the $2$-adic integers as an extension of the set
of non-negative integers (all of which have finitely many binary digits)
to include binary numbers with infinitely many digits. In particular,
note that a $2$-adic integer $\mathfrak{z}$ will be an element of
$\mathbb{Z}_{2}^{\prime}$ if and only if $\mathfrak{z}$ has infinitely
many non-zero $2$-adic digits (infinitely many $1$s).

Letting $n\in\mathbb{N}_{0}$, if we reduce (\ref{eq:Hensel series})
modulo $2^{n}$, we get:
\begin{equation}
\left[\mathfrak{z}\right]_{2^{n}}=\sum_{m=0}^{n-1}d_{m}\left(\mathfrak{z}\right)2^{m},\textrm{ }\forall n\in\mathbb{N}_{0}\label{eq:Hensel series of z mod 2 to the n}
\end{equation}
where, recall, the sum on the right is defined to be $0$ whenever
$n=0$. Note that:
\begin{equation}
\#_{1}\left(\left[\mathfrak{z}\right]_{2^{n}}\right)=\left|\left\{ m\in\left\{ 0,\ldots,n-1\right\} :d_{m}\left(\mathfrak{z}\right)=1\right\} \right|
\end{equation}
Letting $n\rightarrow\infty$, observe that there are exactly two
possibilities for $\#_{1}\left(\left[\mathfrak{z}\right]_{2^{n}}\right)$:

\vphantom{}

I. There is an integer $k\in\mathbb{N}_{0}$ so that $\#_{1}\left(\left[\mathfrak{z}\right]_{2^{n}}\right)=k$
for all sufficiently large $n\in\mathbb{N}_{0}$.

\vphantom{}

II. $\lim_{n\rightarrow\infty}\#_{1}\left(\left[\mathfrak{z}\right]_{2^{n}}\right)\overset{\mathbb{R}}{=}\infty$,
where the $\mathbb{R}$ over the equality is there to remind us that
this limit is being taken in the topology of $\mathbb{R}$.

\vphantom{}

Examining the Hensel series of $\mathfrak{z}$, we see that (I) is
equivalent (if and only if) to $\mathfrak{z}\in\mathbb{N}_{0}$, and
that (II) is equivalent to $\mathfrak{z}\in\mathbb{Z}_{2}^{\prime}$.
We refer to these equivalences as the \textbf{digit count principle}.

Fixing $\mathfrak{z}\in\mathbb{Z}_{2}^{\prime}$, observe that the
$n$th term of $S_{2,3}\left(\mathfrak{z}\right)$ has a $3$-adic
absolute value of:
\begin{equation}
\left|\frac{3^{\#_{1}\left(\left[\mathfrak{z}\right]_{2^{n}}\right)}}{2^{n}}\right|_{3}=3^{-\#_{1}\left(\left[\mathfrak{z}\right]_{2^{n}}\right)}
\end{equation}
By the digit-count principle, since our $\mathfrak{z}$ has infinitely
many $1$s digits, we obtain the $3$-adic decay estimate:
\begin{equation}
\lim_{n\rightarrow\infty}\left|\frac{3^{\#_{1}\left(\left[\mathfrak{z}\right]_{2^{n}}\right)}}{2^{n}}\right|_{3}=\lim_{n\rightarrow\infty}3^{-\#_{1}\left(\left[\mathfrak{z}\right]_{2^{n}}\right)}\overset{\mathbb{R}}{=}0,\textrm{ }\forall\mathfrak{z}\in\mathbb{Z}_{2}^{\prime}
\end{equation}
One of the perks of working in an ultrametric space like $\mathbb{Z}_{3}$
(the ring of $3$-adic integers) is that a series $a_{1}+a_{2}+a_{3}+\cdots$
converges in $\mathbb{Z}_{3}$ \emph{if and only if its} $n$th term
tends to $0$ in $3$-adic absolute value. As we have shown the $3$-adic
decay of $S_{2,3}\left(\mathfrak{z}\right)$ to $0$ for $\mathfrak{z}\in\mathbb{Z}_{2}^{\prime}$,
we conclude that $S_{2,3}\left(\mathfrak{z}\right)$ \textbf{converges}\footnote{It is not difficult to show that this convergence is strictly point-wise.
We can make the $3$-adic convergence of $S_{2,3}\left(\mathfrak{z}\right)$
occur arbitrarily slowly by choosing a $\mathfrak{z}\in\mathbb{Z}_{2}^{\prime}$
whose $1$s digits are separated by sufficiently large number of $0$s.}\textbf{ in $\mathbb{Z}_{3}$ for all $\mathfrak{z}\in\mathbb{Z}_{2}^{\prime}$}.

On the other hand, suppose $\mathfrak{z}\in\mathbb{N}_{0}$. By the
digit count principle, we see that $\left[\mathfrak{z}\right]_{2^{n}}$
will be equal to $\mathfrak{z}$ for all  large $n$; in particular,
$\left[\mathfrak{z}\right]_{2^{n}}=\mathfrak{z}$ if and only if $n\geq\lambda_{2}\left(\mathfrak{z}\right)$.
Consequently:
\begin{equation}
S_{2,3}\left(\mathfrak{z}\right)=\sum_{n=0}^{\infty}\frac{3^{\#_{1}\left(\left[\mathfrak{z}\right]_{2^{n}}\right)}}{2^{n}}=\sum_{n=0}^{\lambda_{2}\left(\mathfrak{z}\right)-1}\frac{3^{\#_{1}\left(\left[\mathfrak{z}\right]_{2^{n}}\right)}}{2^{n}}+\sum_{n=\lambda_{2}\left(\mathfrak{z}\right)}^{\infty}\frac{3^{\#_{1}\left(\mathfrak{z}\right)}}{2^{n}}
\end{equation}
Since $3^{\#_{1}\left(\mathfrak{z}\right)}$ is a constant in the
right-most sum, the geometric series formula gives:
\begin{equation}
S_{2,3}\left(\mathfrak{z}\right)\overset{\mathbb{R}}{=}\sum_{n=0}^{\lambda_{2}\left(\mathfrak{z}\right)-1}\frac{3^{\#_{1}\left(\left[\mathfrak{z}\right]_{2^{n}}\right)}}{2^{n}}+3^{\#_{1}\left(\mathfrak{z}\right)}\frac{\left(\frac{1}{2}\right)^{\lambda_{2}\left(\mathfrak{z}\right)}}{1-\frac{1}{2}}=\frac{2\cdot3^{\#_{1}\left(\mathfrak{z}\right)}}{2^{\lambda_{2}\left(\mathfrak{z}\right)}}+\sum_{n=0}^{\lambda_{2}\left(\mathfrak{z}\right)-1}\frac{3^{\#_{1}\left(\left[\mathfrak{z}\right]_{2^{n}}\right)}}{2^{n}}\label{eq:Sum of S_2,3 over N_0}
\end{equation}
where, note, the series converges in the topology of $\mathbb{R}$
(and/or $\mathbb{C}$). In fact, this convergence occurs \emph{solely
}in $\mathbb{R}$; the geometric series $\sum_{n=\lambda_{2}\left(\mathfrak{z}\right)}^{\infty}2^{-n}$
does not converge $3$-adically, for its terms do not decay to $0$
in $3$-adic absolute value. By this analysis, we see that \textbf{$S_{2,3}\left(\mathfrak{z}\right)$
converges in the topology of $\mathbb{R}$ when $\mathfrak{z}\in\mathbb{N}_{0}$}.
Moreover, as the far right side of (\ref{eq:Sum of S_2,3 over N_0})
shows, \textbf{$S_{2,3}\left(\mathfrak{z}\right)$} is a rational
number, and not just any rational number, but a \emph{$3$-adic integer}.

As such, we can make sense of the formal sum (\ref{eq:S_2,3}) as
a series representation of a $X_{2,3}:\mathbb{Z}_{2}\rightarrow\mathbb{Z}_{3}$,
defined by the rule:
\begin{equation}
X_{2,3}\left(\mathfrak{z}\right)\overset{\textrm{def}}{=}\begin{cases}
\textrm{The sum of }\sum_{n=0}^{\infty}\frac{3^{\#_{1}\left(\left[\mathfrak{z}\right]_{2^{n}}\right)}}{2^{n}}\textrm{ in }\mathbb{Z}_{3} & \textrm{if }\mathfrak{z}\in\mathbb{Z}_{2}^{\prime}\\
\textrm{The sum of }\sum_{n=0}^{\infty}\frac{3^{\#_{1}\left(\left[\mathfrak{z}\right]_{2^{n}}\right)}}{2^{n}}\textrm{ in }\mathbb{R} & \textrm{if }\mathfrak{z}\in\mathbb{N}_{0}
\end{cases}\label{eq:Definition of X_2,3}
\end{equation}
Regardless of the topology being used to sum $X_{2,3}\left(\mathfrak{z}\right)$,
its convergence is strictly point-wise with respect to $\mathfrak{z}$.
We call $X_{2,3}$ a \textbf{$\left(2,3\right)$-adic function}, on
account of the fact that it goes from the $2$-adics to the $3$-adics.
More generally, a \textbf{$\left(p,q\right)$-adic function} is one
from the $p$-adics to the $q$-adics, where $p$ and $q$ are distinct
primes\footnote{Number theorists may be somewhat uncomfortable with this notation.
Traditionally, given a prime $p$, the symbol $q$ is used to represent
a number of the form $p^{\nu}$ where $\nu\in\mathbb{N}_{1}$. When
working with a given prime $p$, the symbol $\ell$ is preferred to
denote a prime different from $p$. Unfortunately, we will need to
use all three symbols. In general, $p$ and $q$ will be reserved
for the primes associated to the domain and target space of our functions;
these two will be fixed. Meanwhile, $\ell$ will be allowed to vary,
as it is the prime we will use to investigate series convergence.
Heuristically, one could say that our functions will go from $p$-adics
to $q$-adics by way of the $\ell$-adic topology, for appropriate
choices of $\ell$.}.
\end{example}
The original motivation for defining the concept of \textbf{frames
}which will be at the core of this paper is that it quickly become
incredibly tedious to have to keep on specifying the dependence of
the topology of convergence on $\mathfrak{z}$. So, to expedite things,
we adopt the following terminology.
\begin{defn}
Let $\left\{ f_{n}\right\} _{n\geq0}$ be a sequence of functions
$f_{n}:\mathbb{Z}_{2}\rightarrow\mathbb{Q}$. Given a function $F:\mathbb{Z}_{2}\rightarrow\mathbb{Q}_{3}$
(where $\mathbb{Q}_{3}$ is the field of $3$-adic numbers; the fraction
field of $\mathbb{Z}_{3}$), we say \textbf{$f_{n}$ converges to
$F$ with respect to the standard $\left(2,3\right)$-adic frame}
(denoted $\lim_{n\rightarrow\infty}f_{n}\left(\mathfrak{z}\right)\overset{\mathcal{F}_{2,3}}{=}F\left(\mathfrak{z}\right)$)
when the following conditions hold:

I. For each $\mathfrak{z}\in\mathbb{N}_{0}$, $\lim_{n\rightarrow\infty}f_{n}\left(\mathfrak{z}\right)\overset{\mathbb{R}}{=}F\left(\mathfrak{z}\right)$.

II. For each $\mathfrak{z}\in\mathbb{Z}_{2}^{\prime}$, $\lim_{n\rightarrow\infty}f_{n}\left(\mathfrak{z}\right)\overset{\mathbb{Q}_{3}}{=}F\left(\mathfrak{z}\right)$,
where the $\mathbb{Q}_{3}$ over the $=$ means that the convergence
occurs in the topology of $\mathbb{Q}_{3}$.
\end{defn}
With this terminology, we can write:
\begin{equation}
X_{2,3}\left(\mathfrak{z}\right)\overset{\mathcal{F}_{2,3}}{=}\sum_{n=0}^{\infty}\frac{3^{\#_{1}\left(\left[\mathfrak{z}\right]_{2^{n}}\right)}}{2^{n}}\label{eq:F-series for X_2,3}
\end{equation}
The $\mathcal{F}_{2,3}$ over the $=$ reminds us of the summation
convention we adopted in (\ref{eq:Definition of X_2,3}) when defining
(\ref{eq:Definition of X_2,3}).

As defined, $X_{2,3}$ admits a series representation ($S_{2,3}\left(\mathfrak{z}\right)$)
which converges to $X_{2,3}\left(\mathfrak{z}\right)$ at every $\mathfrak{z}\in\mathbb{Z}_{2}$,
albeit with the caveat that the topology used to \emph{sum }$S_{2,3}\left(\mathfrak{z}\right)$
depends on $\mathfrak{z}$. This is the novel convergence phenomenon
of which we spoke.

\vphantom{}

Classically, analysis likes working with functions whose outputs can
be understood in terms of ``numbers'', by which we mean elements
of $\mathbb{R}$, $\mathbb{C}$, as well as other fields and the algebraic
structures built on top of them, especially Banach spaces. The idea
that a sequence of functions $\left\{ f_{n}\right\} _{n\geq0}$ might
converge in topology $\mathcal{T}$ but not in topology $\mathcal{T}^{\prime}$
is not new. Any introductory course in modern analysis will be rife
with examples of sequences of functions that converge in $L^{p}$
norm for some $p\in\left[1,\infty\right]$ and diverge for others.
There is also the well-established notion of \textbf{modes of convergence}:
a sequence might converge point-wise almost everywhere on its domain,
but fail to converge uniformly on its domain; a sequence might converge
uniformly on every compact subset of its domain, but fail to converge
uniformly on the \emph{entirety }of its domain. Comparing different
modes of convergence\textemdash $L^{p}$ convergence, uniform, point-wise,
in measure, almost everywhere, weakly, weak-{*}ly, etc.\textemdash is
a mainstay of mathematical analysis.

Yet, even though the topologies defined by these different modes of
convergence are often extraordinarily different, they are united by
the fact that the \emph{field }over which the $f_{n}$s attain their
values remains fixed. However different the topologies of $L^{p}$
space and $L^{q}$ space might be for distinct $p$ and $q$, it is
implicit in the literature\textemdash and, indeed, even in this very
sentence!\textemdash that the functions in the $L^{p}$ space and
the functions in the $L^{q}$ space take values in the same field
(usually $\mathbb{R}$, or $\mathbb{C}$). As our example with $S_{2,3}\left(\mathfrak{z}\right)$
shows, in $\left(p,q\right)$-adic analysis, we can no longer take
this for granted. In addition to classical distinctions like point-wise
convergence versus uniform convergence, studying $\mathcal{F}$-series
like $S_{2,3}\left(\mathfrak{z}\right)$ requires us to consider topologies
of convergence generated by different absolute values.

The modern understanding of a function, generally attributed to Dirichlet,
is that of a rule $f$ which, to every element $x$ of a set $X$
(the \textbf{domain }of $f$) associates an element $f\left(x\right)$
of a set $Y$ (the \textbf{target space }/ \textbf{co-domain }of $f$).
In this view, a formula for \emph{computing }the value of a function
at certain inputs is merely a \emph{representation }of the function.
For example, $\zeta\left(s\right)$ denotes the Riemann Zeta Function,
which may be computed at any complex number $s$ with $\textrm{Re}\left(s\right)>1$
by summing the series:
\begin{equation}
\sum_{n=1}^{\infty}\frac{1}{n^{s}}\label{eq:Zeta}
\end{equation}
Despite this ``local'' representation, $\zeta\left(s\right)$ remains
perfectly well-defined for all $s\in\mathbb{C}\backslash\left\{ 1\right\} $,
even though the series formula diverges for $\textrm{Re}\left(s\right)<1$.
The lesson here is that we distinguish $\zeta\left(s\right)$ (the
function as a \emph{rule}) from the formula (\ref{eq:Zeta}) used
to compute its values for $\textrm{Re}\left(s\right)>1$.

The same reasoning applies to $X_{2,3}$. It is perfectly well-defined
as an abstract rule which associates a $3$-adic integer $X_{2,3}\left(\mathfrak{z}\right)$
to every $2$-adic integer input $\mathfrak{z}$. The series $S_{2,3}$
is merely a representation of $X_{2,3}$, one which requires the standard
$\left(2,3\right)$-adic frame $\mathcal{F}_{2,3}$ in order to be
properly summed at any given $\mathfrak{z}$. Thus, we distinguish
between the function $X_{2,3}$ and $S_{2,3}$, the road we take to
arrive at the $3$-adic integer $X_{2,3}\left(\mathfrak{z}\right)$
for any given $\mathfrak{z}\in\mathbb{Z}_{2}$.

That we can do this without the vastly different natures of $\mathbb{R}$
and $\mathbb{Z}_{3}$ crashing into one another is thanks to the properties
of the \emph{image }of $X_{2,3}$. Technically, we can realize the
restriction $X_{2,3}\mid_{\mathbb{N}_{0}}$ of $X_{2,3}$ to $\mathbb{N}_{0}$
as a function $\mathbb{N}_{0}\rightarrow\mathbb{R}$. However, the
image of this restriction is contained in $\mathbb{Q}$, and not just
any part of $\mathbb{Q}$, but the part of $\mathbb{Q}$ containing
rational numbers whose denominators are co-prime to $3$. All such
rational numbers are elements of $\mathbb{Z}_{3}$, and so, even though
we use the topology of $\mathbb{R}$ to go from $\mathfrak{z}$ to
$X_{2,3}\left(\mathfrak{z}\right)$ for $\mathfrak{z}\in\mathbb{N}_{0}$,
we end in the same place as we get when we use the topology of $\mathbb{Z}_{3}$
to go from $\mathfrak{z}$ to $X_{2,3}\left(\mathfrak{z}\right)$
for $\mathfrak{z}\in\mathbb{Z}_{2}^{\prime}$.

From a bird's-eye perspective, the reason for this bizarre state of
affairs lies in the way metric completions arise from the field $\mathbb{Q}$.
We can\textemdash and will\textemdash view $\mathbb{Q}$ as being
equipped with the discrete topology, and the associated discrete metric
introduced by the trivial absolute value\footnote{If this terminology is unfamiliar to the reader, a brief review of
absolute values on fields is given at the start of \textbf{Section
\ref{sec:1}}.}. In this view, as long as we add only finitely many terms, everything
is alright. However, as soon as we consider infinite series of rational
numbers, we must pick a metric completion of $\mathbb{Q}$ in order
to make sense of the series' sum. This could be considered a (algebraic)
number theoretic analogue to considering a series $F_{N}\left(x\right)\overset{\textrm{def}}{=}\sum_{n=0}^{N-1}f_{n}\left(x\right)$,
where each $f_{n}:\mathbb{R}\rightarrow\mathbb{R}$ is continuous.
While $F_{N}:\mathbb{R}\rightarrow\mathbb{R}$ is continuous for all
$N$, if we let $N\rightarrow\infty$, this need no longer be true.
In general, we will need to pick a ``completion'' of the space of
continuous functions $\mathbb{R}\rightarrow\mathbb{R}$, such as $L^{p}$
for $p\geq1$, in order to make sense of $\lim_{N\rightarrow\infty}F_{N}\left(x\right)$.
Classically, we choose an appropriate function space to study the
limit of a given sequence of functions. In this paper, we choose an
appropriate metric completion of $\mathbb{Q}$ to study the limit
of a sequence of rational numbers.

It is worth mentioning that the author discovered these phenomena
in the course of the number-theoretic investigations conducted in
his PhD dissertation\footnote{These investigations are also detailed in \cite{first blog post paper},
the first of several planned expository papers based on \cite{my dissertation}.} (\cite{my dissertation}). Briefly, letting $p$ be a prime, $\mathcal{F}$-series
like $S_{p,q}\left(\mathfrak{z}\right)$ emerged from the study of
functions $\chi:\mathbb{N}_{0}\rightarrow\mathbb{Q}$ satisfying systems
of functional equations of the form:
\begin{align}
X\left(pn+j\right) & =a_{j}X\left(n\right)+b_{j},\textrm{ }j\in\mathbb{Z}/p\mathbb{Z}\label{eq:p-hydra map numen functional equation}
\end{align}
for all $n\in\mathbb{N}_{0}$, where $a_{j},b_{j}\in\mathbb{Q}$,
with $a_{j}\neq0$ for any $j$, where $\mathbb{Z}/p\mathbb{Z}$ is
our shorthand for the set $\left\{ 0,\ldots,p-1\right\} $. This connection
is spelled out in \textbf{Proposition \ref{prop:pointwise convergence of generalized S}}
(page \pageref{prop:pointwise convergence of generalized S}), as
well as in \textbf{Examples \ref{exa:5}}, \textbf{\ref{exa:6}},\textbf{
\ref{exa:8}} (pages \pageref{exa:5}, \pageref{exa:6}, and \pageref{exa:8},
respectively). Beyond this connection, however, as discussed in \cite{my dissertation},
the existence of this previously unknown form of convergence appears
to open an brand new frontier in non-archimedean analysis. \cite{my dissertation}
and its methods can be used to show that $\mathcal{F}$-series define
$\left(p,q\right)$-adic measures\textemdash continuous, $\mathbb{C}_{q}$-valued
linear functionals on the space $C\left(\mathbb{Z}_{p},\mathbb{C}_{q}\right)$
of continuous functions $f:\mathbb{Z}_{p}\rightarrow\mathbb{C}_{q}$;
here, $\mathbb{C}_{q}$ is the metric completion of the algebraic
closure of the field of $q$-adic numbers ($\mathbb{Q}_{q}$). By
realizing $\mathcal{F}$-series as measures, it becomes possible to
integrate them, which significantly expands the family of ``integrable''
functions compared to the rather degenerate situation in ``classical''
$\left(p,q\right)$-adic analysis (see \cite{my dissertation} for
a comprehensive explanation; if the reader wants to figure it out
on their own, they can consult the exercise at the end of Appendix
A.5 of \cite{Ultrametric Calculus}). As an illustration, \textbf{Example
\ref{exa:7}} (page \pageref{exa:7}) shows how frames can be used
to realize the function $\mathfrak{z}\in\mathbb{Z}_{p}\mapsto\left|\mathfrak{z}\right|_{q}^{r}\in\mathbb{Q}_{q}$
as a $\left(p,q\right)$-adic measure (here, $r\in\mathbb{Q}\backslash\left\{ 0,-1\right\} $).
This then opens the door to bringing methods of differentiation (in
the sense of distributions) into $\left(p,q\right)$-adic analysis,
a completely novel enterprise.

In terms of this paper's overall structure, the reader is encouraged
to view its sections and examples as the write up of the solutions
for an imaginary set of homework problems given in some university
mathematics course. Some of the problems in the set might be of use
in understanding others in the set, but there is no particular goal
in mind. There is no Parnassus to be climbed here, merely a collection
of hills to be toured. At best, it can be said that, as a work of
exposition, this paper is meant to define the formalism of \textbf{frames
}(developed by the author developed in \cite{my dissertation}), and
to illustrate and motivate this definition through a variety of detailed
examples.

It must be stressed that the presentation given here is incomplete
in two respects. Firstly, as this material is completely novel, there
is a large amount of work left to be done before we will have anything
close to a satisfactory theory to explain the phenomena we will encounter.
Secondly, the author is aware of significantly more detail than what
is presented here, and has been forced to omit multiple important
considerations in order to keep the paper from becoming unpublishably
long. To wit, \cite{my dissertation} presents, refines, and expands
parts of the theory of non-Archimedean Fourier analysis first presented
in Schikhof's doctoral dissertation (\cite{Schikhof's Thesis}), and
a more honest presentation of this paper's contents would utilize
\cite{my dissertation}'s Fourier analytic innovations; alas, this
must be postponed until a later publication.

\vphantom{}

Declarations of interest: none

\vphantom{}

This research did not receive any specific grant from funding agencies
in the public, commercial, or not-for-profit sectors.

\section{\label{sec:1}Point-wise Convergence of $S_{p,q}$}

We begin by reviewing some basic definitions.
\begin{defn}[Absolute values]
\label{def:absolute values}Let $K$ be a field. An absolute value\textbf{
}on $K$ is a function $a:K\rightarrow\left[0,\infty\right)$ defined
by the rules:

I. For $x\in K$, $a\left(x\right)=0$ if and only if $x=0_{K}$.

II. $a\left(xy\right)=a\left(x\right)\cdot a\left(y\right)$ for all
$x,y\in K$.

III. $a\left(x+y\right)\leq a\left(x\right)+a\left(y\right)$ for
all $x,y\in K$.

If, in addition, $a$ satisfies:

IV. $a\left(x+y\right)\leq\max\left\{ a\left(x\right),a\left(y\right)\right\} $
for all $x,y\in K$, with equality whenever $a\left(x\right)\neq a\left(y\right)$,
we say $a$ is a \textbf{non-archimedean }absolute value on $K$.
If $a$ does not satisfy (IV), we say $a$ is an \textbf{archimedean
}absolute value.

The \textbf{trivial absolute value }on $K$ is defined by $\alpha\left(x\right)=1$
if $x\neq0_{K}$ and $a\left(x\right)=0$ if and only if $x=0_{K}$.
Every field can be equipped with the trivial absolute value.

An absolute value $a$ on $K$ gives $K$ the structure of a metric
space by way of the formula $d\left(x,y\right)=a\left(x-y\right)$
for all $x,y\in K$.
\end{defn}
\begin{rem}
We allow for the case where $\left|\cdot\right|_{K}$ is the trivial
absolute value, in which case $K$ has the discrete metric / topology.
\end{rem}
\begin{rem}
Note that we \emph{do not }require $K$ to be of characteristic $0$.
Yet-to-be-published investigations by the author (see the introduction
of \cite{first blog post paper} for a sketch of a sample case) have
revealed that the methods of \cite{my dissertation} apply just as
well to the study of the analogue of functional equations of the shape
(\ref{eq:p-hydra map numen functional equation}) where the solutions
are functions from $\mathbb{N}_{0}$ to $\mathbb{F}\left(x\right)$,
the field of rational functions over $\mathbb{F}$, where $\mathbb{F}$
is a finite field. For simplicity's sake, we will confine our investigations
here to fields of characteristic $0$ that emerge as metric completions
of $\mathbb{Q}$ or an extension thereof.
\end{rem}
\begin{defn}[Complete valued fields]
By a \textbf{complete valued field}, we mean a field $K$ which is:

I. Equipped with an absolute value, denoted $\left|\cdot\right|_{K}$.

II. A complete metric space with respect to the metric on $K$ induced
by $\left|\cdot\right|_{K}$.

If, in addition, $K$ is of characteristic $0$ and contains $\mathbb{Q}$
as a subfield, we say $K$ is a \textbf{complete valued number field}.

Meanwhile, a \textbf{complete valued ring }is the ring of integers
of a complete valued field $K$. A \textbf{complete valued number
ring} is then the ring of integers of a complete valued number field.
(Thus, for example, $\mathbb{Z}_{p}$ is the complete valued (number)
ring associated to $\mathbb{Q}_{p}$.)
\end{defn}
\begin{rem}
We adopt the convention of treating $\mathbb{Q}$, $\overline{\mathbb{Q}}$,
and any intermediary extension as complete valued number fields by
equipping them with the \textbf{trivial absolute value} (and hence,
with the \textbf{discrete topology}\footnote{This is a sensible convention, because the trivial absolute value
is the only absolute value with respect to which these discrete fields
are metrically complete.}). Note then that a sequence $x_{n}$ converges in this topology if
and only if $x_{n}=x_{n+1}$ for all sufficiently large $n$. We write
$\left|\cdot\right|_{0}$ to denote the trivial absolute value (thus,
$\left|x\right|_{0}=1$ if $x\neq0$ and$\left|0\right|_{0}=0$).
\end{rem}
\begin{rem}
Just so that there is no confusion, we shall \emph{always} view the
fields $\mathbb{R}$, $\mathbb{C}$, and $\mathbb{Q}_{p}$ (or $\mathbb{Q}_{q}$
or $\mathbb{Q}_{\ell}$ or extensions thereof) as valued fields equipped
with their standard absolute values.
\end{rem}
\begin{defn}[Places of a field, completions, etc.]
Letting $a$ and $\alpha$ be absolute values on a field $K$, we
say $a$ and $\alpha$ are \textbf{equivalent} if there is a positive
real constant $c$ so that $\alpha\left(x\right)=\left(a\left(x\right)\right)^{c}$
for all $x\in K$. This equivalence defines an equivalence relation
on the set of all absolute values of $K$, and an equivalence class
of absolute values on $K$ is called a \textbf{place }of $K$. For
fields like $K=\mathbb{Q}$, there is a one-to-one correspondence
between places of $K$ and prime (ideals) of $K$. We write $V_{K}$
to denote \textbf{the set of all places} of $K$. Given a place $v\in K$,
we write $\left|\cdot\right|_{v}$ to denote the absolute value on
$K$ represented by $v$, and we write $K_{v}$ to denote \textbf{the
field obtained by completing $K$ with respect to} \textbf{$\left|\cdot\right|_{v}$}.
By abuse of notation, we will speak of $\left|\cdot\right|_{v}$ as
an element of $V_{K}$.

Finally, we write $K_{V}$ to denote the set of all fields obtained
by completing $K$ with respect to some $\left|\cdot\right|_{v}\in V_{K}$.
We call $K_{V}$ the \textbf{set of completions of $K$}. Note that
completing $K$ with respect to the trivial absolute value yields
$K$ itself.
\end{defn}
\begin{defn}[Divisibility, primes, etc.]
Given a field $K$, we use the standard notation $\mathcal{O}_{K}$
to denote the ring of integer elements of $K$; thus, for $\overline{\mathbb{Q}}$,
the field of algebraic numbers, we write $\mathcal{O}_{\overline{\mathbb{Q}}}$
to denote the ring of algebraic integers. Recall that given $z\in\mathcal{O}_{K}$,
we say $x\in\mathcal{O}_{K}$ \textbf{divides }$z$ (written $x\mid z$)
if there is a $y\in\mathcal{O}_{K}$ so that $z=xy$, where juxtaposition
denotes multiplication. Letting $1_{K}$ (a.k.a., $1$) denote the
multiplicative identity element of $K$, an element $u\in\mathcal{O}_{K}$
is said to be a \textbf{unit }of $\mathcal{O}_{K}$ if $u\mid1$,
in which case $1/u$ is an element of $\mathcal{O}_{K}$. We say $\ell\in\mathcal{O}_{K}$
is \textbf{prime }if it cannot be written as a product $\ell=xy$,
where $x,y\in\mathcal{O}_{K}$ are not units of $\mathcal{O}_{K}$.
Given $x,y\in\mathcal{O}_{K}$, we say $x$ and $y$ are \textbf{co-prime
}if there is no prime $\ell\in\mathcal{O}_{K}$ so that $\ell\mid x$
and $\ell\mid y$.

Given a prime $\ell\in\mathcal{O}_{K}$, we write $\left|\cdot\right|_{\ell}$
to denote the standard normalization of the \textbf{$\ell$-adic absolute
value} on $K$. (In case the reader does not know how to compute $\left|\cdot\right|_{\ell}$
when $\ell$ is not an ordinary prime number ($2,3,5,7,\ldots$),
they are advised to consult a source on algebraic number theory, such
as the first chapter of \cite{Childress}.) More generally, every
place of $K$ corresponds to a unique prime (ideal) of $K$ and vice-versa.

Letting $K$ be a number field (a finite-degree extension of $\mathbb{Q}$),
given a prime $\ell\in\mathcal{O}_{K}$, we write $K_{\ell}$ to denote
the non-archimedean field obtained by completing $K$ with respect
to the $\ell$-adic absolute value. $\mathcal{O}_{K_{\ell}}$ then
denotes the ring of integer elements of $K_{\ell}$. A prime number
$p\in\left\{ 2,3,5,\ldots\right\} $ is said to \textbf{lie under
}$\ell$ if $p\mid\ell$ in $\mathcal{O}_{K}$; equivalently, $\ell$
is said to \textbf{lie over }$p$.

It is a basic fact algebraic number theory that, for any number field
$K$ and any prime $\ell\in\mathcal{O}_{K}$, there is a unique prime
$p\in\left\{ 2,3,5,\ldots\right\} $ such that $\ell$ lies over $p$
(see \cite{Marcus}, for example). As such, in an abuse of notation,
given a prime $\ell\in\mathcal{O}_{K}$ for a number field $K$, we
write $\mathbb{C}_{\ell}$ to denote the metric completion of the
algebraic closure of $\mathbb{Q}_{p}$, where $p$ is the unique prime
that lies under $\ell$. $\mathbb{C}_{\ell}$ then contains the metric
completion of the algebraic closure of $K_{\ell}$.

Finally, we speak of the standard absolute value on $\mathbb{R}$
and $\mathbb{C}$ as the \textbf{$\infty$-adic absolute value} (viz.
the ``infinite prime'' spoken of in number theory) and write it
as $\left|\cdot\right|_{\infty}$. We then write $\mathbb{R}$ as
$\mathbb{Q}_{\infty}$, to denote the fact that $\mathbb{R}$ is the
metric completion of $\mathbb{Q}$ with respect to the absolute value
$\left|\cdot\right|_{\infty}$. This makes it easy to talk about series
convergence in a single breath. Also, because places correspond to
primes, we speak of \textbf{finite places }when a place corresponds
to a finite prime and speak of \textbf{infinite places }when a place
corresponds to an infinite prime.
\end{defn}
With those definitions out of the way, we now characterize the topologies
in which $S_{p,q}$ converges.
\begin{prop}
\label{prop:Characterization of S_p,q's convergence}Let $p,q\in\mathcal{O}_{\overline{\mathbb{Q}}}\backslash\left\{ 0\right\} $,
so that $p$ and $q$ are integral elements $K\overset{\textrm{def}}{=}\mathbb{Q}\left(p,q\right)$.
Suppose further that $p$ and $q$ are co-prime as elements of $\mathcal{O}_{K}$.
Then:

I. Let $\mathfrak{z}\in\mathbb{Z}_{2}^{\prime}$, and let $\ell$
be any prime\footnote{Because we will later speak of infinite primes, for pedantry's sake,
it should be mentioned that, here, $\ell$ must be a \emph{finite
}prime.} of $\mathcal{O}_{K}$. Then, $S_{p,q}\left(\mathfrak{z}\right)$
converges in $\mathbb{C}_{\ell}$ if and only if $\left|p\right|_{\ell}=1$
and $\left|q\right|_{\ell}<1$. Moreover, when convergent, $S_{p,q}\left(\mathfrak{z}\right)$
is necessarily an element of $\mathcal{O}_{K_{\ell}}$.

II. Let $\mathfrak{z}\in\mathbb{N}_{0}$. Then, $S_{p,q}\left(\mathfrak{z}\right)$
converges in $\mathbb{C}$ if and only if $\left|p\right|_{\infty}>1$.
Moreover, if convergent, $S_{p,q}\left(\mathfrak{z}\right)$ is an
element of $K$.
\end{prop}
Proof:

I. Let $p,q$ be co-prime non-zero algebraic integers, and let $\ell$
be a prime of $\mathcal{O}_{K}$.

i. Pick $\mathfrak{z}\in\mathbb{Z}_{2}^{\prime}$ and suppose $S_{p,q}\left(\mathfrak{z}\right)$
converges in $\mathbb{C}_{\ell}$. Then, choosing $\mathfrak{z}=-1$,
the series:
\begin{equation}
S_{p,q}\left(-1\right)=\sum_{n=0}^{\infty}\frac{q^{\#_{1}\left(\left[-1\right]_{2^{n}}\right)}}{p^{n}}=\sum_{n=0}^{\infty}\left(\frac{q}{p}\right)^{n}
\end{equation}
converges in $\mathbb{C}_{\ell}$ if and only if $\left|q/p\right|_{\ell}<1$,
and hence, if and only if $\left|q\right|_{\ell}<\left|p\right|_{\ell}$.
Thus, all we need to do is show $\left|p\right|_{\ell}=1$, then the
bound on $\left|q\right|_{\ell}$ will follow automatically.

Now, since $p$ and $q$ are $K$-integers, they are also $K_{\ell}$-integers.
Thus, $\left|p\right|_{\ell},\left|q\right|_{\ell}\leq1$. So, by
way of contradiction, suppose $\left|p\right|_{\ell}<1$. By definition
of the $\ell$-adic absolutue value, this forces $\ell\mid p$. Since
$\left|q\right|_{\ell}<\left|p\right|_{\ell}<1$, this forces $\ell\mid q$
as well, but we assumed $p$ and $q$ were co-prime; thus, no prime
can divide both of them. This is the desired contradiction, and so
$\left|p\right|_{\ell}$ must be $1$, and we are done.

ii. Suppose $\left|p\right|_{\ell}=1$ and $\left|q\right|_{\ell}<1$.
Letting $\mathfrak{z}\in\mathbb{Z}_{2}^{\prime}$ be arbitrary, observe
that:
\begin{equation}
\left|\frac{q^{\#_{1}\left(\left[\mathfrak{z}\right]_{2^{n}}\right)}}{p^{n}}\right|_{\ell}\overset{\mathbb{R}}{=}\left|q\right|_{\ell}^{\#_{1}\left(\left[\mathfrak{z}\right]_{2^{n}}\right)}
\end{equation}
Since $\mathfrak{z}\in\mathbb{Z}_{2}^{\prime}$, $\mathfrak{z}$ has
infinitely many $1$s digits, and so $\#_{1}\left(\left[\mathfrak{z}\right]_{2^{n}}\right)\rightarrow\infty$
in $\mathbb{R}$ as $n\rightarrow\infty$. Since $\left|q\right|_{\ell}<1$,
we then have:
\begin{equation}
\lim_{n\rightarrow\infty}\left|\frac{q^{\#_{1}\left(\left[\mathfrak{z}\right]_{2^{n}}\right)}}{p^{n}}\right|_{\ell}\overset{\mathbb{R}}{=}\lim_{n\rightarrow\infty}\left|q\right|_{\ell}^{\#_{1}\left(\left[\mathfrak{z}\right]_{2^{n}}\right)}\overset{\mathbb{R}}{=}0,\textrm{ }\forall\mathfrak{z}\in\mathbb{Z}_{2}^{\prime}
\end{equation}
As in our introductory example, $K_{\ell}$ is an ultrametric space,
and so the vanishing of this limit is then equivalent to the convergence
of $S_{p,q}\left(\mathfrak{z}\right)$ in $K_{\ell}$. Since $K$
is a field extension of $\mathbb{Q}$, so too is $K_{\ell}$ a field
extension of $\mathbb{Q}_{\ell}$, and so, $K_{\ell}$ is necessarily
contained in $\mathbb{C}_{\ell}$, the algebraic closure of $\mathbb{Q}_{\ell}$.
This shows that the convergence of $S_{p,q}\left(\mathfrak{z}\right)$
in $K_{\ell}$ implies the convergence of $S_{p,q}\left(\mathfrak{z}\right)$
in $\mathbb{C}_{\ell}$.

Finally, keeping $\mathfrak{z}\in\mathbb{Z}_{2}^{\prime}$, note that
$p$ and $q$ are both elements of $\mathcal{O}_{K_{\ell}}$. Moreover,
since $\left|p\right|_{\ell}=1$, $p$ is a unit of $\mathcal{O}_{K_{\ell}}$,
and so $q^{\#_{1}\left(\left[\mathfrak{z}\right]_{2^{n}}\right)}/p^{n}$
is in $\mathcal{O}_{K_{\ell}}$ or all $n\geq0$. So, $S_{p,q}\left(\mathfrak{z}\right)$
is a sum of elements of $\mathcal{O}_{K_{\ell}}$ which converges
with respect to the $\ell$-adic absolute value. Since $\mathcal{O}_{K_{\ell}}$
is complete as a metric space with respect to the $\ell$-adic absolute
value, this guarantees that $S_{p,q}\left(\mathfrak{z}\right)\in\mathcal{O}_{K_{\ell}}$.

\vphantom{}II.

i. Pick $\mathfrak{z}\in\mathbb{N}_{0}$ and suppose $S_{p,q}\left(\mathfrak{z}\right)$
converges in $\mathbb{C}$. Just like in \textbf{Example \ref{exa:1}},
observe that the integer $\left[\mathfrak{z}\right]_{2^{n}}$ is equal
to $\mathfrak{z}$ for all $n\geq\lambda_{2}\left(\mathfrak{z}\right)$.
Thus, pick $N\geq\lambda_{2}\left(\mathfrak{z}\right)$, we can write:
\begin{align*}
\sum_{n=0}^{N}\frac{q^{\#_{1}\left(\left[\mathfrak{z}\right]_{2^{n}}\right)}}{p^{n}} & =\sum_{n=0}^{\lambda_{2}\left(\mathfrak{z}\right)-1}\frac{q^{\#_{1}\left(\left[\mathfrak{z}\right]_{2^{n}}\right)}}{p^{n}}+\sum_{n=\lambda_{2}\left(\mathfrak{z}\right)}^{N}\frac{q^{\#_{1}\left(\left[\mathfrak{z}\right]_{2^{n}}\right)}}{p^{n}}\\
 & =\sum_{n=0}^{\lambda_{2}\left(\mathfrak{z}\right)-1}\frac{q^{\#_{1}\left(\left[\mathfrak{z}\right]_{2^{n}}\right)}}{p^{n}}+q^{\#_{1}\left(\mathfrak{z}\right)}\sum_{n=\lambda_{2}\left(\mathfrak{z}\right)}^{N}\frac{1}{p^{n}}
\end{align*}
The convergence of $S_{p,q}\left(\mathfrak{z}\right)$ in $\mathbb{C}$
then tells us that the geometric series:
\begin{equation}
\sum_{n=\lambda_{2}\left(\mathfrak{z}\right)}^{N}\frac{1}{p^{n}}
\end{equation}
converges in $\mathbb{C}$ to: 
\begin{equation}
\frac{p^{-\lambda_{2}\left(\mathfrak{z}\right)}}{1-\frac{1}{p}}
\end{equation}
as $N\rightarrow\infty$. This forces $\left|p\right|_{\infty}>1$,
and shows that $S_{p,q}\left(\mathfrak{z}\right)$:
\begin{equation}
S_{p,q}\left(\mathfrak{z}\right)\overset{\mathbb{C}}{=}\frac{p}{p-1}\frac{q^{\#_{1}\left(\mathfrak{z}\right)}}{p^{\lambda_{2}\left(\mathfrak{z}\right)}}+\sum_{n=0}^{\lambda_{2}\left(\mathfrak{z}\right)-1}\frac{q^{\#_{1}\left(\left[\mathfrak{z}\right]_{2^{n}}\right)}}{p^{n}},\textrm{ }\forall\mathfrak{z}\in\mathbb{N}_{0}\label{eq:S_p,q at non-negative integer z}
\end{equation}
which is an element of $K$.

ii. Conversely, for $\mathfrak{z}\in\mathbb{N}_{0}$, if $\left|p\right|_{\infty}$,
the above shows $S_{p,q}\left(\mathfrak{z}\right)$ converges in $\mathbb{C}$
to an element of $K$ given by (\ref{eq:S_p,q at non-negative integer z}).

Q.E.D.

\vphantom{}This discussion motivates the following definition:
\begin{defn}
Let $p,q\in\mathcal{O}_{\overline{\mathbb{Q}}}\backslash\left\{ 0\right\} $
be co-prime in $\mathcal{O}_{K}$, where $K=\mathbb{Q}\left(p,q\right)$.
Let $\ell$ be a prime of $\mathcal{O}_{K}$, and suppose that $\left|p\right|_{\infty}>1$,
$\left|p\right|_{\ell}=1$, and $\left|q\right|_{\ell}<1$. Then,
we define the function $X_{p,q}:\mathbb{Z}_{2}\rightarrow K_{\ell}$
by:
\begin{equation}
X_{p,q}\left(\mathfrak{z}\right)\overset{\textrm{def}}{=}\begin{cases}
\textrm{The sum of }\sum_{n=0}^{\infty}\frac{p^{\#_{1}\left(\left[\mathfrak{z}\right]_{2^{n}}\right)}}{q^{n}}\textrm{ in }K_{\ell} & \textrm{if }\mathfrak{z}\in\mathbb{Z}_{2^{\prime}}\\
\textrm{The sum of }\sum_{n=0}^{\infty}\frac{q^{\#_{1}\left(\left[\mathfrak{z}\right]_{2^{n}}\right)}}{p^{n}}\textrm{ in }\mathbb{C} & \textrm{if }\mathfrak{z}\in\mathbb{N}_{0}
\end{cases}\label{eq:Definition of X_p,q}
\end{equation}
Here, \textbf{Proposition \ref{prop:Characterization of S_p,q's convergence}
}guarantees that the sum of the series defining $X_{p,q}\left(\mathfrak{z}\right)$
at any given $\mathfrak{z}\in\mathbb{N}_{0}$ will be an element of
$K$, and hence, an element of the completion $K_{\ell}$.
\end{defn}
Having fully elucidated $S_{p,q}$'s convergence properties, we can
now introduce the notion of the \textbf{standard $\left(p,q\right)$-adic
frame}:
\begin{defn}
Let $p,q$ be distinct primes in $\mathbb{Z}$, and consider a sequence
of functions $f_{n}:\mathbb{Z}_{p}\rightarrow\overline{\mathbb{Q}}$.
We say the $f_{n}$s \textbf{converge to $f:\mathbb{Z}_{p}\rightarrow\mathbb{C}_{q}$
with respect to the} \textbf{standard $\left(p,q\right)$-adic frame
}(denoted $\lim_{n\rightarrow\infty}f_{n}\left(\mathfrak{z}\right)\overset{\mathcal{F}_{p,q}}{=}f\left(\mathfrak{z}\right)$)
whenever the following holds:

I. $f\left(n\right)\in\overline{\mathbb{Q}}$ for all $n\in\mathbb{N}_{0}$.

II. $\lim_{n\rightarrow\infty}f_{n}\left(\mathfrak{z}\right)\overset{\mathbb{C}_{q}}{=}f\left(\mathfrak{z}\right)$
for all $\mathfrak{z}\in\mathbb{Z}_{p}^{\prime}$.

III. $\lim_{n\rightarrow\infty}f_{n}\left(\mathfrak{z}\right)\overset{\mathbb{C}}{=}f\left(\mathfrak{z}\right)$
for all $\mathfrak{z}\in\mathbb{N}_{0}$.

We use the symbol $\mathcal{F}_{p,q}$ to denote this frame. Using
this symbol, we can avoid circumlocutions when discussing the convergence
of functions like $X_{p,q}$. For example, if $p$ and $q$ are distinct
primes of $\mathbb{Z}$, we have:
\begin{equation}
X_{p,q}\left(\mathfrak{z}\right)\overset{\mathcal{F}_{2,q}}{=}\sum_{n=0}^{\infty}\frac{q^{\#_{1}\left(\left[\mathfrak{z}\right]_{2^{n}}\right)}}{p^{n}},\textrm{ }\forall\mathfrak{z}\in\mathbb{Z}_{2}
\end{equation}
which is clearly much more succint than having to repeatedly write
``we sum the series in the topology of $\mathbb{Z}_{q}$ when $\mathfrak{z}\in\mathbb{Z}_{2}^{\prime}$
and sum it in the topology of $\mathbb{C}$ when $\mathfrak{z}\in\mathbb{N}_{0}$''.
But frames turn out to be more than just a notational convenience.
\end{defn}
\begin{rem}
The author refers to $\mathcal{F}_{p,q}$ as the ``standard'' $\left(p,q\right)$-adic
frame simply because it was the one that occurred the most often in
\cite{my dissertation}.
\end{rem}
Classically, the convergence of a sequence or series of functions
at a point is a mostly binary problem: either convergence occurs,
or it does not\textemdash in which case the limit might be some kind
of infinity, or entirely non-existent. However, $\mathcal{F}$-series
defy this dichotomy. Our previous examples show how the topology in
which an $\mathcal{F}$-series converges might change if we change
the point at which we are evaluating the series. But things can get
even more intricate than that: we can have expressions which converge
in \emph{multiple }topologies at a \emph{single }point, as our next
example shows.
\begin{example}
\label{exa:2}Consider the $\mathcal{F}$-series: 
\begin{equation}
S_{7,6}\left(\mathfrak{z}\right)=\sum_{n=0}^{\infty}\frac{6^{\#_{1}\left(\left[\mathfrak{z}\right]_{2^{n}}\right)}}{7^{n}}\label{eq:S_7,6}
\end{equation}
Fixing $n\geq0$, observe that the value of $\#_{1}\left(\left[\mathfrak{z}\right]_{2^{n}}\right)$
for $\mathfrak{z}\in\mathbb{Z}_{2}$ is maximized when $\left[\mathfrak{z}\right]_{2^{n}}=2^{n}-1$,
in which case we get $\#_{1}\left(2^{n}-1\right)=n$ (all of $2^{n}-1$'s
$n$ binary digits are $1$s). Thus, we have the estimate:
\begin{equation}
\frac{6^{\#_{1}\left(\left[\mathfrak{z}\right]_{2^{n}}\right)}}{7^{n}}\leq\frac{6^{\#_{1}\left(\left[-1\right]_{2^{n}}\right)}}{7^{n}}=\left(\frac{6}{7}\right)^{n},\textrm{ }\forall\mathfrak{z}\in\mathbb{Z}_{2}
\end{equation}
Consequently, we have the (real!) bound:
\begin{equation}
\sup_{\mathfrak{z}\in\mathbb{Z}_{2}}S_{7,6}\left(\mathfrak{z}\right)\leq\sum_{n=0}^{\infty}\frac{6^{\#_{1}\left(\left[-1\right]_{2^{n}}\right)}}{7^{n}}\leq\sum_{n=0}^{\infty}\left(\frac{6}{7}\right)^{n}=\frac{1}{1-\frac{6}{7}}=7
\end{equation}
So, by the Weierstrass M-test, the series $S_{7,6}\left(\mathfrak{z}\right)$
converges in $\mathbb{R}$ \emph{uniformly} for $\mathfrak{z}\in\mathbb{Z}_{2}$,
and therefore the (positive!) real-valued function:
\begin{equation}
\mathfrak{z}\in\mathbb{Z}_{2}\mapsto\textrm{sum of }\sum_{n=0}^{\infty}\frac{6^{\#_{1}\left(\left[\mathfrak{z}\right]_{2^{n}}\right)}}{7^{n}}\in\mathbb{R}
\end{equation}
is \emph{uniformly continuous}\footnote{Since this function is uniformly continuous, it is then $L^{1}$-integrable
with respect to the real-valued Haar probability measure on $\mathbb{Z}_{2}$.
So, we can use the theory of Fourier analysis on locally compact abelian
groups (Pontryagin duality, etc.) to study it. This also includes
distributional notions of differentiation (see \cite{Vladimirov - the big paper about complex-valued distributions over the p-adics},
for example).}. We will denote this function by $X_{7,6}^{\left(\infty\right)}:\mathbb{Z}_{2}\rightarrow\mathbb{R}$.

However, this is not the \emph{only} way we can get a function out
of this $\mathcal{F}$-series. As defined, note that, by \textbf{Proposition
\ref{prop:Characterization of S_p,q's convergence}} (setting $K=\mathbb{Q}$),
$S_{7,6}\left(\mathfrak{z}\right)$ converges in both the standard
$\left(2,3\right)$-adic frame \emph{and} the standard $\left(2,2\right)$-adic
frame. This gives us realizations $X_{7,6}^{\left(2\right)}:\mathbb{Z}_{2}\rightarrow\mathbb{Q}_{2}$
and $X_{7,6}^{\left(3\right)}:\mathbb{Z}_{2}\rightarrow\mathbb{Q}_{3}$
defined by:
\begin{equation}
X_{7,6}^{\left(2\right)}\left(\mathfrak{z}\right)\overset{\mathcal{F}_{2,2}}{=}\sum_{n=0}^{\infty}\frac{6^{\#_{1}\left(\left[\mathfrak{z}\right]_{2^{n}}\right)}}{7^{n}},\textrm{ }\forall\mathfrak{z}\in\mathbb{Z}_{2}\label{eq:X_7,6(2)}
\end{equation}
and:
\begin{equation}
X_{7,6}^{\left(3\right)}\left(\mathfrak{z}\right)\overset{\mathcal{F}_{2,3}}{=}\sum_{n=0}^{\infty}\frac{6^{\#_{1}\left(\left[\mathfrak{z}\right]_{2^{n}}\right)}}{7^{n}},\textrm{ }\forall\mathfrak{z}\in\mathbb{Z}_{2}\label{eq:X_7,6(3)}
\end{equation}
respectively. In this way, $X_{7,6}^{\left(\infty\right)}$, $X_{7,6}^{\left(2\right)}$,
and $X_{7,6}^{\left(3\right)}$ represent \emph{three different functions!}

Despite these differences, $X_{7,6}^{\left(\infty\right)}$, $X_{7,6}^{\left(2\right)}$,
and $X_{7,6}^{\left(3\right)}$ nevertheless share a certain amount
of common ground. In particular\textemdash as we will soon show\textemdash for
any $\mathfrak{z}\in\mathbb{Z}_{2}\cap\mathbb{Q}$, there is a rational
number $r\left(\mathfrak{z}\right)$ so that $X_{7,6}^{\left(\infty\right)}\left(\mathfrak{z}\right)$,
$X_{7,6}^{\left(2\right)}\left(\mathfrak{z}\right)$, and $X_{7,6}^{\left(3\right)}\left(\mathfrak{z}\right)$
are all equal to $r\left(\mathfrak{z}\right)$, with their respective
series representations converging to $r\left(\mathfrak{z}\right)$
in their respective topologies. On the other hand, if $\mathfrak{z}$
is an \emph{irrational} $2$-adic integer \footnote{$\mathfrak{z}\in\mathbb{Z}_{2}\backslash\mathbb{Q}$; i.e., the $2$-adic
digits of $\mathfrak{z}$ are not eventually periodic.}, the values to which $X_{7,6}^{\left(\infty\right)}\left(\mathfrak{z}\right)$,
$X_{7,6}^{\left(2\right)}\left(\mathfrak{z}\right)$, and $X_{7,6}^{\left(3\right)}\left(\mathfrak{z}\right)$
converge in their respective topologies will be \emph{incomparable}
to one another. As metric completions of $\mathbb{Q}$, $\mathbb{R}$,
$\mathbb{Q}_{2}$, and $\mathbb{Q}_{3}$ all contain $\mathbb{Q}$
as subfields. Where these completions differ is at their irrational
elements.

That being said, because we will need it later, let us adopt the convention
of defining $X_{7,6}$ (without any superscripts) as a function $X_{7,6}:\mathbb{Z}_{2}\rightarrow\mathbb{Q}_{3}$
with respect to the standard $\left(2,3\right)$-adic frame, with:
\begin{equation}
X_{7,6}\left(\mathfrak{z}\right)\overset{\textrm{def}}{=}\begin{cases}
\textrm{sum of }\sum_{n=0}^{\infty}\frac{6^{\#_{1}\left(\left[\mathfrak{z}\right]_{2^{n}}\right)}}{7^{n}}\textrm{ in }\mathbb{R} & \textrm{if }\mathfrak{z}\in\mathbb{N}_{0}\\
\textrm{sum of }\sum_{n=0}^{\infty}\frac{6^{\#_{1}\left(\left[\mathfrak{z}\right]_{2^{n}}\right)}}{7^{n}}\textrm{ in }\mathbb{Q}_{3} & \textrm{if }\mathfrak{z}\in\mathbb{Z}_{2}^{\prime}
\end{cases}\label{eq:Definition of X_7,6}
\end{equation}
\end{example}
\begin{prop}
\label{prop:2}If $\mathfrak{z}\in\mathbb{Z}_{2}\cap\mathbb{Q}$,
then $X_{7,6}^{\left(\infty\right)}\left(\mathfrak{z}\right)$, $X_{7,6}^{\left(2\right)}\left(\mathfrak{z}\right)$,
$X_{7,6}^{\left(3\right)}\left(\mathfrak{z}\right)$ all represent
the same rational number.
\end{prop}
The reason for this is\textemdash to use the terminology of \cite{Conrad on p-adic series}\textemdash the
\textbf{universality of the geometric series}.
\begin{thm}[Geometric series universality]
\label{thm:geometric series universality}Let $r\in\overline{\mathbb{Q}}$,
$\mathbb{F},\mathbb{K}$ be any two fields obtained by completing
the field $\mathbb{Q}\left(r\right)$ with respect to the absolute
values $\left|\cdot\right|_{\mathbb{F}}$ and $\left|\cdot\right|_{\mathbb{K}}$,
respectively. Suppose that the geometric series:
\begin{equation}
\sum_{n=0}^{\infty}r^{n}
\end{equation}
converges to $s_{\mathbb{F}}\in\mathbb{F}$ with respect to $\left|\cdot\right|_{\mathbb{F}}$
and converges to $s_{\mathbb{K}}$ with respect to $\left|\cdot\right|_{\mathbb{K}}$.
Then, both $s_{\mathbb{F}}$ and $s_{\mathbb{K}}$ are elements of
$\overline{\mathbb{Q}}$, and are in fact equal as elements of $\overline{\mathbb{Q}}$,
with: 
\begin{equation}
s_{\mathbb{F}}=s_{\mathbb{K}}=\frac{1}{1-r}
\end{equation}
\end{thm}
Proof: Let everything be as given, and suppose that the geometric
series converges to $s_{\mathbb{F}}$ in $\mathbb{F}$:
\begin{equation}
s_{\mathbb{F}}=\sum_{n=0}^{\infty}r^{n}\overset{\mathbb{F}}{=}\frac{1}{1-r}
\end{equation}
As such:
\begin{equation}
s_{\mathbb{F}}-\frac{1}{1-r}=0
\end{equation}
and so, $s_{\mathbb{F}}$ is the root of the polynomial $z-\frac{1}{1-r}$.
Since $r$ is algebraic, so is $1/\left(1-r\right)$, and, in fact,
we have that $s_{\mathbb{F}}=1/\left(1-r\right)$. Since $\mathbb{F}$
was arbitrary, we conclude that $s_{\mathbb{F}}=s_{\mathbb{K}}$ for
any $\mathbb{F}$ and $\mathbb{K}$ in which the geometric series
converges.

Q.E.D.

\vphantom{}

This property of geometric series is actually the exception rather
than the norm\footnote{In fact, the assumption that such a universality property holds for
general power series is arguably the oldest mistake in $p$-adic analysis,
dating back to Kurt Hensel himself and his erroneous $p$-adic proof
that the constant $e$ is a transcendental number \cite{my dissertation}.}. In general, just because a series converges in the topologies of
two different completions of $\mathbb{Q}$ (or a finite-degree extension
thereof), there is no reason to conclude that the sum of our series
in one of the two topologies represents the same object obtained by
summing the series in the other topology. \cite{Conrad on p-adic series}
gives examples of this, and discusses the issue in greater detail.

\vphantom{}

Now, to prove \textbf{Proposition \ref{prop:2}}, we need a simple
observation, as well as a useful function worth giving its own symbol\footnote{In this respect, we follow \cite{my dissertation,first blog post paper}}.
\begin{defn}
Define $B_{2}:\mathbb{N}_{0}\rightarrow\mathbb{Q}$ by:
\begin{equation}
B_{2}\left(n\right)\overset{\textrm{def}}{=}\begin{cases}
0 & \textrm{if }n=0\\
\frac{n}{1-2^{\lambda_{2}\left(n\right)}} & \textrm{if }n\geq1
\end{cases}\label{eq:Definition B_2}
\end{equation}
This function has the effect of sending a given $n\in\mathbb{N}_{0}$
to the rational number (and, in fact, $2$-adic integer) obtained
by concatenating infinitely many copies of the sequence of $n$'s
$2$-adic digits. Note that We can extend $B_{2}$ to a function $\mathbb{Z}_{2}\rightarrow\mathbb{Z}_{2}$
by defining: $B_{2}\left(\mathfrak{z}\right)$ to be $\mathfrak{z}$
for all $\mathfrak{z}\in\mathbb{Z}_{2}^{\prime}$.
\end{defn}
\begin{rem}
It may be of interest to note that $B_{2}$ is idempotent over $\mathbb{Z}_{2}$:
$\left(B_{2}\circ B_{2}\right)\left(\mathfrak{z}\right)=B_{2}\left(\mathfrak{z}\right)$
for all $\mathfrak{z}\in\mathbb{Z}_{2}$.
\end{rem}
For us, the key property of $B_{2}$ is that we can explicitly relate
the number of $1$s in the binary representation of $B_{2}\left(n\right)$
mod $2^{k}$ to the number of $1$s in the binary representation of
$n$.
\begin{prop}
\label{prop:3}Let $n\in\mathbb{N}_{1}$. Then, for all $k\in\mathbb{N}_{0}$:
\begin{equation}
\#_{1}\left(\left[B_{2}\left(n\right)\right]_{2^{k}}\right)=\left\lfloor \frac{k}{\lambda_{2}\left(n\right)}\right\rfloor \#_{1}\left(n\right)+\#_{1}\left(\left[n\right]_{2^{\left[k\right]_{\lambda_{2}\left(n\right)}}}\right)\label{eq:Number of 1s in the projection of B_2 of n mod 2 to the k}
\end{equation}
where $\left[n\right]_{2^{\left[k\right]_{\lambda_{2}\left(n\right)}}}$
is the value of $n$ mod $2^{\left[k\right]_{\lambda_{2}\left(n\right)}}$,
where $\left[k\right]_{\lambda_{2}\left(n\right)}$ is the value of
$k$ mod $\lambda_{2}\left(n\right)$.
\end{prop}
Proof: Write out the Hensel series of $B_{2}\left(n\right)$ and figure
out the pattern.

Q.E.D.

\vphantom{}With \textbf{Proposition \ref{prop:3}}, we can find a
closed-form expression for $S_{p,q}\left(\mathfrak{z}\right)$ for
any $\mathfrak{z}\in\mathbb{Z}_{2}\cap\mathbb{Q}$ at which the series
converges in \emph{some} metric completion of $\overline{\mathbb{Q}}$.
\begin{lem}
\label{lem:1}Let $p,q\in\overline{\mathbb{Q}}\backslash\left\{ 0\right\} $,
and let $K$ be the completion of $\mathbb{Q}\left(p,q\right)$ with
respect to an arbitrary non-trivial absolute value $\left|\cdot\right|_{K}$.
If $S_{p,q}\left(\mathfrak{z}\right)$ converges in $K$ at $\mathfrak{z}\in\mathbb{Z}_{2}\cap\mathbb{Q}$,
then $S_{p,q}\left(\mathfrak{z}\right)$ is an element of $\mathbb{Q}\left(p,q\right)$
of the form:
\begin{equation}
A+\frac{B}{1-r}
\end{equation}
where $A,B,r$ are elements of $\mathbb{Q}\left(p,q\right)$, with
$r\neq1$.
\end{lem}
Proof: Fix $\mathfrak{z}\in\mathbb{Z}_{2}\cap\mathbb{Q}$, and suppose
$S_{p,q}\left(\mathfrak{z}\right)$ converges in $K$. Since $\mathfrak{z}$
is rational, its sequence of $2$-adic digits is eventually periodic.
As such, there are\footnote{Specifically, $m$ is the non-negative integer whose $2$-adic digits
are the initial string of $\mathfrak{z}$'s $2$-adic digits \emph{before}
$\mathfrak{z}$'s digits become periodic. Meanwhile, $n$ is the non-negative
integer whose $2$-adic digits generate the periodic part of $\mathfrak{z}$'s
digit sequence. For example, if the periodic part of $\mathfrak{z}$
is $01010101\ldots$, then the generating sequence is $01$, which
represents $n=2$. On the other hand, for something like $110110110110\ldots$,
note that the generating sequence is $1101$ (which represents $n=11$),
rather than $110$. As a $2$-adic integer, we have $110=1+1\cdot2+0\cdot2^{2}=3=11$.
This is because, when writing a non-zero $2$-adic integers' digits
left-to-right in order of increasing powers of $2$, the right-most
digit has to be $1$.} non-negative integers $m,n$ so that $m+2^{\lambda_{2}\left(m\right)}B_{2}\left(n\right)$.
Since:
\begin{equation}
\#_{1}\left(\left[\mathfrak{z}\right]_{2^{k}}\right)=\#_{1}\left(\left[m+2^{\lambda_{2}\left(m\right)}B_{2}\left(n\right)\right]_{2^{k}}\right)=\#_{1}\left(m\right)+\#_{1}\left(\left[B_{2}\left(n\right)\right]_{2^{k}}\right),\textrm{ }\forall k\geq\lambda_{2}\left(m\right)
\end{equation}
we can then write:
\begin{align*}
\sum_{k=0}^{\infty}\frac{q^{\#_{1}\left(\left[\mathfrak{z}\right]_{2^{k}}\right)}}{p^{k}} & \overset{K}{=}\overbrace{\sum_{k=0}^{\lambda_{2}\left(m\right)-1}\frac{q^{\#_{1}\left(\left[m+2^{\lambda_{2}\left(m\right)}B_{2}\left(n\right)\right]_{2^{k}}\right)}}{p^{k}}}^{0\textrm{ if }m=0}+\sum_{k=\lambda_{2}\left(m\right)}^{\infty}\frac{q^{\#_{1}\left(\left[m+2^{\lambda_{2}\left(m\right)}B_{2}\left(n\right)\right]_{2^{k}}\right)}}{p^{k}}\\
 & \overset{K}{=}\sum_{k=0}^{\lambda_{2}\left(m\right)-1}\frac{q^{\#_{1}\left(\left[m\right]_{2^{k}}\right)}}{p^{k}}+\sum_{k=\lambda_{2}\left(m\right)}^{\infty}\frac{q^{\#_{1}\left(m\right)+\#_{1}\left(\left[B_{2}\left(n\right)\right]_{2^{k}}\right)}}{p^{k}}\\
\left(\textrm{\textbf{Proposition }\textbf{\ref{prop:3}}}\right); & \overset{K}{=}\sum_{k=0}^{\lambda_{2}\left(m\right)-1}\frac{q^{\#_{1}\left(\left[m\right]_{2^{k}}\right)}}{p^{k}}+q^{\#_{1}\left(m\right)}\sum_{k=\lambda_{2}\left(m\right)}^{\infty}\frac{q^{\left\lfloor \frac{k}{\lambda_{2}\left(n\right)}\right\rfloor \#_{1}\left(n\right)+\#_{1}\left(\left[n\right]_{2^{\left[k\right]_{\lambda_{2}\left(n\right)}}}\right)}}{p^{k}}
\end{align*}
Here, note that the summand:
\begin{equation}
\frac{q^{\left\lfloor \frac{k}{\lambda_{2}\left(n\right)}\right\rfloor \#_{1}\left(n\right)+\#_{1}\left(\left[n\right]_{2^{\left[k\right]_{\lambda_{2}\left(n\right)}}}\right)}}{p^{k}}
\end{equation}
is periodic in $k$, with period $\lambda_{2}\left(n\right)$.

After pulling out a finite number of terms from the infinite sum to
get some $a\in\mathbb{Q}\left(p,q\right)$, we get:
\begin{equation}
\sum_{k=\lambda_{2}\left(m\right)}^{\infty}\frac{q^{\left\lfloor \frac{k}{\lambda_{2}\left(n\right)}\right\rfloor \#_{1}\left(n\right)+\#_{1}\left(\left[n\right]_{2^{\left[k\right]_{\lambda_{2}\left(n\right)}}}\right)}}{p^{k}}=a+\sum_{j=0}^{\lambda_{2}\left(n\right)-1}q^{\#_{1}\left(\left[n\right]_{2^{j}}\right)}\sum_{k=N}^{\infty}\left(\frac{q^{\#_{1}\left(n\right)}}{p}\right)^{k}
\end{equation}
for some $N\geq0$. Since the geometric series converges in $K$,
summing it yields:
\begin{align*}
\sum_{k=0}^{\infty}\frac{q^{\#_{1}\left(\left[\mathfrak{z}\right]_{2^{k}}\right)}}{p^{k}} & =\sum_{k=0}^{\lambda_{2}\left(m\right)-1}\frac{q^{\#_{1}\left(\left[m\right]_{2^{k}}\right)}}{p^{k}}+q^{\#_{1}\left(m\right)}\left(a+\sum_{j=0}^{\lambda_{2}\left(n\right)-1}q^{\#_{1}\left(\left[n\right]_{2^{j}}\right)}\sum_{k=N}^{\infty}\frac{q^{k\#_{1}\left(n\right)}}{p^{k}}\right)\\
 & =\underbrace{q^{\#_{1}\left(m\right)}a+\sum_{k=0}^{\lambda_{2}\left(m\right)-1}\frac{q^{\#_{1}\left(\left[m\right]_{2^{k}}\right)}}{p^{k}}}_{A}+\frac{\overbrace{q^{\#_{1}\left(m\right)}\left(\frac{q^{\#_{1}\left(n\right)}}{p}\right)^{N}\sum_{j=0}^{\lambda_{2}\left(n\right)-1}q^{\#_{1}\left(\left[n\right]_{2^{j}}\right)}}^{B}}{1-\underbrace{\frac{q^{\#_{1}\left(n\right)}}{p}}_{r}}
\end{align*}
as desired.

Q.E.D.
\begin{rem}
This result is very characteristic of frame theory and its applications.
Rather than treating our underlying fields as being fixed, we fix
a base field $K$ and then let ourselves freely move from one completion
of $K$ to another, choosing whichever completion(s) best serves our
needs.
\end{rem}
Finally, we prove \textbf{Proposition \ref{prop:2}}:

\vphantom{}

\textbf{Proof of Proposition \ref{prop:2}}: If $\mathfrak{z}\in\mathbb{Z}_{2}\cap\mathbb{Q}$,
then, by \textbf{Lemma} \textbf{\ref{lem:1}}, $S_{7,6}\left(\mathfrak{z}\right)$
is a geometric series over $\mathbb{Q}$. Since $X_{7,6}^{\left(\infty\right)}\left(\mathfrak{z}\right)$,
$X_{7,6}^{\left(2\right)}\left(\mathfrak{z}\right)$, and $X_{7,6}^{\left(3\right)}\left(\mathfrak{z}\right)$
are sums of this geometric series in topologies obtained by completing
$\mathbb{Q}$ with respect to three different absolute values, \textbf{Geometric
Series Universality }guarantees that the three functions take the
same value at $\mathfrak{z}$, which is some rational number given
by $S_{7,6}\left(\mathfrak{z}\right)$.

Q.E.D.

\section{\label{sec:2}An Introduction to Frames}

In \cite{my dissertation}'s terminology, all of the frames in this
section are examples of \textbf{one-dimensional frames}. The notion
of the \textbf{degree }of a frame is not present in \cite{my dissertation};
the author only devised it after the fact, during the writing of this
paper. In our inclusion of the degree of the frame, in allowing our
valued fields to have positive characteristic, and in explicitly framing
the discussion (pun intended) in terms of the set of all completions
of a given field $K$, our presentation here is more general than
\cite{my dissertation}'s. At the same time, our terminology is not
as comprehensive as \cite{my dissertation}'s, in that we omit Fourier-theoretic
issues such as classifications of measures (scalar or vector-valued)
and Fourier multipliers.

As a rule, where \cite{my dissertation} conflicts with this paper,
this paper should be given precedence, at least until such time as
the author gets around to re-writing \cite{my dissertation} and,
ideally, getting it published it as a textbook. \textbf{Section 3
}of this paper revists the defintions given in \textbf{Section 2}
in slightly greater generality. To that end, in the terminology of
\textbf{Section 3}, the frames considered in \textbf{Section 2 }are
all \textbf{frames of} \textbf{dimension }$1$; the dimension of a
frame is defined in \textbf{Section 3}.
\begin{defn}[Frames of dimension $1$]
Given a topological space $X$ and a valued field $K$, a\textbf{
$K$-frame $\mathcal{F}$ on $X$ of dimension $1$ }is a map $\mathcal{F}:X\rightarrow K_{V}$
which, to each $x\in X$, outputs a completion of $K$, denoted $\mathcal{F}\left(x\right)$.
We will sometimes omit reference to the underlying field $K$ and
simply speak of a \textbf{frame }when either $K$ is known by context,
or is not relevant to the discussion. In practice, $K$ will be a
discrete field with the trivial absolute value (examples include $\mathbb{Q}$,
$\overline{\mathbb{Q}}$, a finite field $\mathbb{F}$, and its algebraic
closure $\overline{\mathbb{F}}$).
\end{defn}
\begin{rem}
From this point on, unless specifically stated otherwise, we fix a
topological space $X$, a valued field $K$, and a $K$-frame $\mathcal{F}$
on $X$.
\end{rem}
In order to use frames to speak of functions on $X$, we need a set
to serve as our functions' co-domain (a.k.a. target space). This is
the \textbf{image }of our frame.
\begin{defn}
The \textbf{image }of $\mathcal{F}$, denoted $I\left(\mathcal{F}\right)$,
is defined as the union of all the $\mathcal{F}\left(x\right)$s:
\begin{equation}
I\left(\mathcal{F}\right)\overset{\textrm{def}}{=}\bigcup_{x\in X}\mathcal{F}\left(x\right)\label{eq:Image of a frame}
\end{equation}
\end{defn}
\begin{rem}
Here, $I\left(\mathcal{F}\right)$ is only a set. \emph{It does not
come with a topology. }Nevertheless, since $\mathcal{F}\left(x\right)$
is a completion of $K$, observe that $\bigcap_{x\in X}\mathcal{F}\left(x\right)$
\emph{is} a field\textemdash in fact, it is a \emph{topological field},
and it contains $K$ as a subfield. Consequently, $I\left(\mathcal{F}\right)$
contains $K$, as well.
\end{rem}
\begin{example}
Let $p$ and $q$ be distinct primes, and consider the standard $\left(p,q\right)$-adic
frame. $\mathcal{F}_{p,q}$ satisfies:
\begin{equation}
\mathcal{F}_{p,q}\left(\mathfrak{z}\right)=\begin{cases}
\mathbb{C} & \textrm{if }\mathfrak{z}\in\mathbb{N}_{0}\\
\mathbb{C}_{q} & \textrm{if }\mathfrak{z}\in\mathbb{Z}_{2}^{\prime}
\end{cases},\textrm{ }\forall\mathfrak{z}\in\mathbb{Z}_{2}
\end{equation}
Moreover, $\mathbb{C}$ and $\mathbb{C}_{q}$ are both completions
of $\overline{\mathbb{Q}}$. Thus, the standard $\left(p,q\right)$-adic
frame is a $\overline{\mathbb{Q}}$-frame. Its image is $\mathbb{C}\cup\mathbb{C}_{q}$.
\end{example}
Frames have at least two main uses. The first is that they provide
us with a convenient shorthand for referring to sequences of sequences
of functions on $X$ whose topologies of convergence depend on the
point in $X$ at which the functions in the sequence are being evaluated.
But they also serve a second, deeper purpose. Much of the power of
classical notions of convergence such as uniform convergence on a
set or $L^{p}$ convergence lies in the fact that they come from norms,
with which we can define Banach spaces of well-behaved functions.
By classical standards, the peculiar convergence properties of $\mathcal{F}$-series
makes them strongly pathological. As such, both for practicality's
sake and in the hope of attaining a greater understanding, it makes
sense to look for function spaces where we can work with $\mathcal{F}$-series
in greater generality. Frames are the first step toward that goal.
Using frames, we can define a notion of a \textbf{compatible function},
which we can use to obtain function spaces with desirable properties.
\begin{defn}
Given a frame $\mathcal{F}$ on $X$, we say a function $f:X\rightarrow I\left(\mathcal{F}\right)$
is \textbf{compatible with $\mathcal{F}$} if $f\left(x\right)\in\mathcal{F}\left(x\right)$
for all $x\in X$. We then write $C\left(\mathcal{F}\right)$ to denote
\textbf{the set of all $\mathcal{F}$-compatible functions}. Additionally,
given a function $f:X\rightarrow Y$, where $Y$ is a set, we say
a frame $\mathcal{F}$ on $X$ is \textbf{compatible }with $f$ if
$Y\subseteq I\left(\mathcal{F}\right)$ and $f\in C\left(\mathcal{F}\right)$.
\end{defn}
With the notion of compatible functions, we can define convergence
with respect to a frame.
\begin{defn}
Consider a frame $\mathcal{F}$ and a function $f\in C\left(\mathcal{F}\right)$.
We say a sequence $\left\{ f_{n}\right\} _{n\geq0}\subseteq C\left(\mathcal{F}\right)$
\textbf{converges to $f$ with respect to $\mathcal{F}$ }if, for
every $x\in X$, $f_{n}\left(x\right)$ converges to $f\left(x\right)$
in the topology of $\mathcal{F}\left(x\right)$. We denote this convergence
by writing:
\begin{equation}
f\left(x\right)\overset{\mathcal{F}}{=}\lim_{n\rightarrow\infty}f_{n}\left(x\right)
\end{equation}
\end{defn}
\begin{rem}
In practice, the compatible functions we work with arise in the following
way: we have a valued field $K$ (typically a number field, equipped
with the trivial absolute value), and an infinite series $f\left(x\right)=\sum_{n=0}^{\infty}f_{n}\left(x\right)$,
where $f_{n}:X\rightarrow K$. Note that $f_{n}$ is then $\mathcal{F}$-compatible
for all $K$-frames $\mathcal{F}$ for all $n$; the same is true
of the partial sums $\sum_{n=0}^{N}f_{n}\left(x\right)$. The idea
here is that we pick a $K$-frame $\mathcal{F}$ on $X$ so that,
for every $x$, $\mathcal{F}\left(x\right)$ is the field completion
of $K$ in which $f\left(x\right)$ happens to converge.
\end{rem}
In many ``nice'' cases, $I\left(\mathcal{F}\right)$ will turn out
to be overkill. In these cases, the outputs of our compatible function
$f$ will be confined to some completion $K^{\prime}\in K_{V}$; i.e.,
 $f\left(x\right)\in K^{\prime}$ for all $x\in X$. In this way,
$f$ can be realized as a function $X\rightarrow K^{\prime}$. Meanwhile,
the frame $\mathcal{F}$ keeps track of the topologies we need to
sum the $\mathcal{F}$-series representation $F\left(x\right)$ at
any given $x\in X$. This idea is important enough that it deserves
its own name: \textbf{extension}.
\begin{defn}
Let $K$ be a field, let $\mathcal{F}$ be a $K$-frame on $X$, and
let $K^{\prime}$ be a completion of $K$ (or the ring of integers
of such a completion). Given $f\in C\left(\mathcal{F}\right)$, we
say $f$ \textbf{extends to $K^{\prime}$} if $f\left(x\right)\in K^{\prime}$
for all $x\in X$. As such, we can view $f$ as a function $X\rightarrow K^{\prime}$.
We then write $C_{K^{\prime}}\left(\mathcal{F}\right)$ to denote
the set of all $f\in C\left(\mathcal{F}\right)$ which extend to $K^{\prime}$.
\end{defn}
\begin{example}
Let $q$ be an odd prime, and consider the standard $\left(2,q\right)$-adic
frame $\mathcal{F}_{2,q}$. By \textbf{Proposition \ref{prop:Characterization of S_p,q's convergence}},
$X_{2,q}\left(\mathfrak{z}\right)$ is compatible with $\mathcal{F}_{2,q}$.
Note that $X_{2,q}\left(\mathfrak{z}\right)\in\mathbb{Z}_{q}$ for
all $\mathfrak{z}\in\mathbb{Z}_{2}^{\prime}$, and that $X_{2,q}\left(\mathfrak{z}\right)\in\mathbb{Q}\cap\mathbb{Z}_{q}$
for all $\mathfrak{z}\in\mathbb{N}_{0}$; i.e., the real number to
which $X_{2,q}\left(\mathfrak{z}\right)$ converges at $\mathfrak{z}\in\mathbb{N}_{0}$
is a rational number which is \emph{also} a $q$-adic integer. Since
$\mathbb{Z}_{q}$ and $\mathbb{Q}\cap\mathbb{Z}_{q}$ are both contained
in $\mathbb{Z}_{q}$, we then have that $X_{2,q}\left(\mathfrak{z}\right)$
extends to $\mathbb{Z}_{q}$. $X_{2,q}\left(\mathfrak{z}\right)$
also extends to $\mathbb{Q}_{q}$, and to any other complete valued
ring (or field) containing $\mathbb{Z}_{q}$.
\end{example}
All of these definitions are not without purpose. Many of the most
important developments in analysis involve function spaces, particularly
\textbf{Banach spaces} and related structures ($L^{p}$ spaces, Hardy
spaces, $C^{*}$-algebras, tempered distributions, etc.). Ideally,
one would like to be able to create similar theories for $\mathcal{F}$-series
and related functions. However, any theory-building of that sort is
immediately derailed by $\mathcal{F}$-series' unprecedented convergence
properties. Banach spaces and their relatives are all\textbf{ vector
spaces}, and\textemdash as we all know\textemdash vector spaces work
with a single field of coefficients. At first glance, this would seem
to directly conflict with the frame theoretic idea of studying convergence
from multiple different metric completions of our base field $K$.

Frames are just what we need to bridge these two worlds. We just need
to keep in mind the following mantra: \emph{think point-wise}. Case
in point: let $\mathcal{F}$ be a $K$-frame on $X$, and fix $x\in X$.
Then, given any $\mathcal{F}$-compatible functions $f,g\in C\left(\mathcal{F}\right)$,
note that $f\left(x\right)$ and $g\left(x\right)$ are both elements
of the field $\mathcal{F}\left(x\right)$. As such, we can make sense
of expressions like $f\left(x\right)+g\left(x\right)$ and $f\left(x\right)g\left(x\right)$
by using the addition and multiplication operations of the field $\mathcal{F}\left(x\right)$,
and the result will be an element of $\mathcal{F}\left(x\right)$.
Keeping our mantra in mind, note that this recipe works perfectly
well for any $x\in X$. As such, we can make sense of expressions
$f+g$ and $fg$ by having them output the elements $f\left(x\right)+g\left(x\right)$
and $f\left(x\right)g\left(x\right)$ in $\mathcal{F}\left(x\right)$
for each $x\in X$. In doing so, we have shown that $f+g$ and $fg$
are functions $X\rightarrow I\left(\mathcal{F}\right)$. This proves:
\begin{prop}
\label{prop:compatible functions form a ring}$C\left(\mathcal{F}\right)$
has the structure of a commutative, unital ring under the operations
of point-wise addition and point-wise multiplication.
\end{prop}
Proof: The additive identity element of $C\left(\mathcal{F}\right)$
is the constant function which, to each $x$, outputs the $0$ of
$\mathcal{F}\left(x\right)$. Since $\mathcal{F}\left(x\right)$ is
the completion of $K$, this constant function is simply the function
which sends every $x\in X$ to the $0$ of $K$. In addition to the
above, $C\left(\mathcal{F}\right)$ has a multiplicative identity
element; this is the function which, to each $x$, outputs the $1$
of $\mathcal{F}\left(x\right)$. Just like with $0$, we can make
sense of this function as the constant function from $X$ to $K$
whose output is the $1$ of $K$.

Q.E.D.
\begin{rem}
Likewise, note that the additive identity element of $C\left(\mathcal{F}\right)$
is the constant function from $X$ to $K$ which outputs the $0$
of $K$, since that $0$ is also the $0$ of every completion of $K$.
\end{rem}
The proof of the previous proposition brings up a key point: because
$\mathcal{F}\left(x\right)$ is a completion of $K$ for all $x\in X$,
$K$ ends up being the ``common denominator'' relating the values
taken by any $f\in C\left(\mathcal{F}\right)$. As a result, we can
upgrade $C\left(\mathcal{F}\right)$ from a ring to a $K$-algebra.
Moreover, this same argument works for any subfield of $K$: anything
which can act on $K$ can act on a completion of $K$.
\begin{prop}
\label{prop:5}Let $K$ be a field, let $\mathcal{F}$ be a $K$-frame,
and let $\mathbb{F}$ be $K$ or any subfield thereof. Then, $C\left(\mathcal{F}\right)$
is a unital algebra over $\mathbb{F}$ with point-wise addition of
functions, scalar multiplication of functions by elements of $\mathbb{F}$,
and point-wise multiplication of functions as the addition, scalar
multiplication, and algebra multiplication operations of $C\left(\mathcal{F}\right)$,
respectively.
\end{prop}
Proof: Use \textbf{Proposition \ref{prop:compatible functions form a ring}},
along with the observation that for any $\alpha\in\mathbb{F}$, and
any $f\in C\left(\mathcal{F}\right)$, note that $\alpha f\left(x\right)$
is an element of $\mathcal{F}\left(x\right)$ for all $x\in X$. This
makes $\alpha f\in C\left(\mathcal{F}\right)$.

Q.E.D.

\vphantom{}By the same reasoning, we can get an algebra structure
for $C_{K^{\prime}}\left(\mathcal{F}\right)$\textemdash and a normed
one, at that.
\begin{prop}
Let $\mathcal{F}$ be a $K$-frame, and let $K^{\prime}$ be a completion
of $K$. Then, $C_{K^{\prime}}\left(\mathcal{F}\right)$ is a unital
algebra over $K^{\prime}$, with point-wise addition of functions,
scalar multiplication of functions by elements of $K^{\prime}$, and
point-wise multiplication of functions as the addition, scalar multiplication,
and algebra multiplication operations of $C_{K^{\prime}}\left(\mathcal{F}\right)$,
respectively. If in addition $X$ is compact, $C_{K^{\prime}}\left(\mathcal{F}\right)$
is a normed vector space under the $K^{\prime}$-sup norm:
\begin{equation}
\left\Vert f\right\Vert _{K^{\prime}}\overset{\textrm{def}}{=}\sup_{x\in X}\left|f\left(x\right)\right|_{K^{\prime}}
\end{equation}
where $\left|\cdot\right|_{K^{\prime}}$ is the absolute value of
$K^{\prime}$.
\end{prop}
Proof: The argument that $C_{K^{\prime}}\left(\mathcal{F}\right)$
is an algebra over $K^{\prime}$ is all but identical to the proof
of that $C\left(\mathcal{F}\right)$ is an algebra over $\mathbb{F}$,
where $\mathbb{F}$ is $K$ or a subfield thereof. Since $K^{\prime}$
is a complete valued field, the set of all $K^{\prime}$-valued functions
on $X$ is a Banach space under the norm $\left\Vert \cdot\right\Vert _{K^{\prime}}$.
$C_{K^{\prime}}\left(\mathcal{F}\right)$ is a vector space contained
in this Banach space, and therefore inherits the norm\textemdash though
not necessarily the completeness.

Q.E.D.
\begin{question}
What is the completion of $C_{K^{\prime}}\left(\mathcal{F}\right)$?
\end{question}
As mentioned previously, the terminology of the ``standard'' $\left(p,q\right)$-adic
frame arose simply because that frame was the one the author encountered
most frequently in his dissertation. One significant difference between
\cite{my dissertation} and the present paper is that we \emph{will
}consider examples of non-standard frames. First, however, we need
to extend our terminology for $\mathbb{Z}_{2}$ to $\mathbb{Z}_{p}$,
for an arbitrary prime $p$.
\begin{defn}
Fix a prime $p$. We write $\mathbb{Z}_{p}^{\prime}\overset{\textrm{def}}{=}\mathbb{Z}_{p}\backslash\mathbb{N}_{0}$.
For any $\mathfrak{z}\in\mathbb{Z}_{p}$, we write $\left[\mathfrak{z}\right]_{p^{n}}$
to denote the projection of $\mathfrak{z}$ mod $p^{n}$, viewed as
an integer in the set $\left\{ 0,\ldots,p^{n}-1\right\} $. (Note
we have $\left[\mathfrak{z}\right]_{p^{0}}=0$ for all $\mathfrak{z}\in\mathbb{Z}_{p}$.)
For any $j\in\left(\mathbb{Z}/p\mathbb{Z}\right)\backslash\left\{ 0\right\} $
and any $m\in\mathbb{N}_{0}$, we define $\#_{p:j}\left(m\right)$
to be the number of $j$s in the $p$-adic digits of $m$. We write
$\lambda_{p}\left(m\right)$ to denote the number of $p$-adic digits
of $m$; note that $\lambda_{p}\left(m\right)=\left\lceil \log_{p}\left(m+1\right)\right\rceil $.
We then define $\#_{p:0}\left(m\right)$ by:
\begin{equation}
\#_{p:0}\left(m\right)\overset{\textrm{def}}{=}\lambda_{p}\left(m\right)-\sum_{j=1}^{p-1}\#_{p:j}\left(m\right)\label{eq:Definition of number of zeroes}
\end{equation}
so that:
\begin{equation}
\sum_{j=0}^{p-1}\#_{p:j}\left(m\right)=\lambda_{p}\left(m\right),\textrm{ }\forall m\in\mathbb{N}_{0}
\end{equation}
For all $j\in\left\{ 1,\ldots,p-1\right\} $, we have the functional
equations:
\begin{equation}
\#_{p:j}\left(pm+k\right)=\#_{p:j}\left(m\right)+\left[k\overset{p}{\equiv}j\right],\textrm{ }\forall m\in\mathbb{N}_{0},\textrm{ }\forall k\in\mathbb{Z}/p\mathbb{Z}\label{eq:pound p colon j functional equations}
\end{equation}
where $\overset{p}{\equiv}$ means ``congruent modulo $p$'', and
where $\left[k\overset{p}{\equiv}j\right]$ is an example of \textbf{Iverson
bracket notation}, being equal to $1$ whenever the enclosed statement
(in this case, $k\overset{p}{\equiv}j$) is true and being equal to
$0$ otherwise. We also have the functional equations:
\begin{equation}
\lambda_{p}\left(pm+j\right)=\lambda_{p}\left(m\right)+1,\textrm{ }\forall m\in\mathbb{N}_{0},\textrm{ }\forall j\in\mathbb{Z}/p\mathbb{Z}
\end{equation}
Note that these functional equations uniquely characterize $\#_{p:j}$
and $\lambda_{p}$ as functions $\mathbb{N}_{0}\rightarrow\mathbb{N}_{0}$.

The functional equation:
\begin{equation}
\#_{p:0}\left(pm+k\right)=\#_{p:0}\left(m\right)+\left[k\overset{p}{\equiv}0\right]
\end{equation}
is true for all integers $m\in\mathbb{N}_{0}$ and $k\in\mathbb{Z}/p\mathbb{Z}$
\emph{except }the case $m=k=0$; there, the left-hand side is $0$,
but the right-hand side is $1$.
\end{defn}
Now, to some examples of non-standard frames. Unsurprisingly, these
examples emerge from $\mathcal{F}$-series, albeit of a form more
general that $S_{p,q}$. For simplicity, we will work solely with
positive integer parameters, though our results (such as \textbf{Proposition
\ref{prop:7}}) can be extended for much more general parameters in
the manner of \textbf{Proposition \ref{prop:Characterization of S_p,q's convergence}}.
\begin{defn}
Let $p$ be a prime number, and for each $j\in\mathbb{Z}/p\mathbb{Z}$,
let $q_{j}$ be a positive integer. Also, choose a positive rational
number $d$. Then, we define the formal series $S_{d;q_{0},\ldots,q_{p-1}}\left(\mathfrak{z}\right)$
by:
\begin{equation}
S_{d;q_{0},\ldots,q_{p-1}}\left(\mathfrak{z}\right)\overset{\textrm{def}}{=}\sum_{n=0}^{\infty}\frac{1}{d^{n}}\prod_{j=0}^{p-1}q_{j}^{\#_{p:j}\left(\left[\mathfrak{z}\right]_{p^{n}}\right)}\label{eq:Def of generalized S}
\end{equation}
where $\mathfrak{z}$ is a formal $p$-adic integer variable. Like
$S_{p,q}$, we call (\ref{eq:Def of generalized S}) an \textbf{$\mathcal{F}$-series}.
\end{defn}
Our next result is the analogue of \textbf{Proposition \ref{prop:Characterization of S_p,q's convergence}}
for $S_{d;q_{0},\ldots,q_{p-1}}\left(\mathfrak{z}\right)$.
\begin{prop}
\label{prop:7}Let $p$ be a prime, and for each $j\in\mathbb{Z}/p\mathbb{Z}$,
choose a positive rational number $q_{j}$. Also, choose a positive
rational number $d$. With these hypotheses:

I. Let $\ell$ be any prime, including the infinite prime. If $\mathfrak{z}\in\mathbb{N}_{0}$,
then $S_{d;q_{0},\ldots,q_{p-1}}\left(\mathfrak{z}\right)$ converges
in $\mathbb{Q}_{\ell}$ if and only if $\left|d\right|_{\ell}>1$.

II. Let $\ell$ be any finite prime. If $\mathfrak{z}\in\mathbb{Z}_{p}^{\prime}$,
then $S_{d;q_{0},\ldots,q_{p-1}}\left(\mathfrak{z}\right)$ converges
in $\mathbb{Q}_{\ell}$ if and only if:
\begin{equation}
\lim_{n\rightarrow\infty}\frac{1}{\left|d\right|_{\ell}^{n}}\prod_{j=0}^{p-1}\left|q_{j}\right|_{\ell}^{\#_{p:j}\left(\left[\mathfrak{z}\right]_{p^{n}}\right)}\overset{\mathbb{R}}{=}0
\end{equation}
\end{prop}
Proof: If $\mathfrak{z}\in\mathbb{N}_{0}$, then $\mathfrak{z}$ has
finitely many $p$-adic digits, and $S_{d;q_{0},\ldots,q_{p-1}}\left(\mathfrak{z}\right)$
is a geometric series in $1/d$, and hence we get the condition from
(I). If $\mathfrak{z}\in\mathbb{Z}_{p}^{\prime}$, observe that the
limit condition given in (II) is then equivalent to the convergence
of $S_{d;q_{0},\ldots,q_{p-1}}\left(\mathfrak{z}\right)$ in $\mathbb{Q}_{\ell}$
for any \emph{finite }prime $\ell$, thanks to the basic principles
of ultrametric analysis.

Q.E.D.

\vphantom{}As mentioned in the Introduction, the author discovered
$\mathcal{F}$-series and frames in the context of solving certain
systems of functional equations. Before proving a formal result about
the connections between these topics, let us consider an example\footnote{Excluding the explicit discussion of frames,\textbf{ Example \ref{exa:5}}
is essentially a sketch of the content of \cite{first blog post paper}.} to get the main ideas across.
\begin{example}
\label{exa:5}Let $p$ be an odd prime, and fix non-negative rational
constants $a_{j},b_{j}\in\mathbb{Q}$ for $j\in\mathbb{Z}/p\mathbb{Z}$,
where $0<a_{0}<1$, $b_{0}=0$, and $a_{j}\neq0$ for any $j\in\mathbb{Z}/p\mathbb{Z}$.
Then, consider:
\begin{align}
f\left(pn+j\right) & =a_{j}f\left(n\right)+b_{j},\textrm{ }\forall n\in\mathbb{N}_{0},\textrm{ }\forall j\in\mathbb{Z}/p\mathbb{Z}\label{eq:System}
\end{align}
Setting $n=j=0$, we obtain: 
\begin{equation}
f\left(0\right)=a_{0}f\left(0\right)+\underbrace{b_{0}}_{0}=a_{0}f\left(0\right)
\end{equation}
and hence, $\left(1-a_{0}\right)f\left(0\right)=0$. Since $a_{0}\neq0$,
this forces any solution of (\ref{eq:System}) to satisfy $f\left(0\right)=0$.
By induction, one can easily show that there is a unique function
$\chi:\mathbb{N}_{0}\rightarrow\mathbb{Q}$ satisfying (\ref{eq:System}).
This is done by using (\ref{eq:System}) and the initial condition
$\chi\left(0\right)=0$ to recursively compute $\chi\left(m\right)$
for all $m\in\mathbb{N}_{0}$. For example, for $j\in\mathbb{Z}/p\mathbb{Z}$,
we have:
\begin{equation}
\chi\left(j\right)=a_{j}\chi\left(0\right)+b_{j}=b_{j}
\end{equation}
And so:
\begin{equation}
\chi\left(p\cdot1+j\right)=a_{j}\chi\left(1\right)+b_{j},\forall j\in\mathbb{Z}/p\mathbb{Z}
\end{equation}
and so on and so forth.

Now, \textbf{suppose for each $j\in\left\{ 1,\ldots,p-1\right\} $
there is a prime number $q_{j}$ so that the $q_{j}$-adic absolute
value of $a_{j}$ is less than $1$}. By induction, by recursively
applying (\ref{eq:System}), we can formally express $\lim_{N\rightarrow\infty}\chi\left(\left[\mathfrak{z}\right]_{p^{N}}\right)$
as an infinite series:
\begin{align*}
\chi\left(j_{0}+j_{1}p+j_{2}p^{2}+\cdots\right) & =b_{j_{0}}+a_{j_{0}}\chi\left(j_{1}+j_{2}p+j_{3}p^{2}+\cdots\right)\\
 & =b_{j_{0}}+a_{j_{0}}\left(b_{j_{1}}+a_{j_{1}}\chi\left(j_{2}+j_{3}p+j_{4}p^{2}+\cdots\right)\right)\\
 & =b_{j_{0}}+a_{j_{0}}\left(b_{j_{1}}+a_{j_{1}}\left(b_{j_{2}}+a_{j_{2}}\chi\left(j_{3}+j_{4}p+j_{5}p^{2}+\cdots\right)\right)\right)\\
 & \vdots\\
 & =b_{j_{0}}+a_{j_{0}}b_{j_{1}}+a_{j_{0}}a_{j_{1}}b_{j_{2}}+a_{j_{0}}a_{j_{1}}a_{j_{2}}b_{j_{3}}+\cdots
\end{align*}
In order to determine this series' convergence properties, let us
first note that the $n$th term of this series is 
\begin{equation}
b_{j_{n}}\times\prod_{k=0}^{n-1}a_{j_{k}},\textrm{ }\forall n\in\mathbb{N}_{0}
\end{equation}
where, as per our conventions, the product is $1$ whenever $n=0$.
We can express this as a function of $\mathfrak{z}$ by noting that,
whne $\mathfrak{z}\in\mathbb{Z}_{p}^{\prime}$ the $j_{k}$s are precisely
the $p$-adic digits of $\mathfrak{z}$, viz. $\mathfrak{z}$'s Hensel
series representation:
\begin{equation}
\mathfrak{z}=\sum_{k=0}^{\infty}j_{k}p^{k}
\end{equation}
This then gives us:
\begin{equation}
b_{j_{n}}\times\prod_{j=0}^{p-1}a_{j}^{\#_{p:j}\left(\left[\mathfrak{z}\right]_{p^{n}}\right)},\textrm{ }\forall n\in\mathbb{N}_{0}\label{eq:nth term}
\end{equation}
as the $n$th term of our series, with:
\begin{equation}
\lim_{N\rightarrow\infty}\chi\left(\left[\mathfrak{z}\right]_{p^{N}}\right)=\lim_{N\rightarrow\infty}\sum_{n=0}^{N-1}b_{j_{n}}\prod_{j=0}^{p-1}a_{j}^{\#_{p:j}\left(\left[\mathfrak{z}\right]_{p^{n}}\right)}\label{eq:rising limit as limit of series}
\end{equation}
That is, for any given $N\in\mathbb{N}_{1}$, $j\in\mathbb{Z}/p\mathbb{Z}$,
and $n\in\left\{ 0,\ldots,N-1\right\} $the exponent of $a_{j}$ in
the $n$th term of the infinite series representing $\lim_{N\rightarrow\infty}\chi\left(\left[\mathfrak{z}\right]_{p^{N}}\right)$
is the number of times $j$ occurs in the $p$-adic digits of the
non-negative integer $\left[\mathfrak{z}\right]_{p^{N}}$.

First, let us deal with (\ref{eq:rising limit as limit of series})
in the case where $\mathfrak{z}\in\mathbb{N}_{0}$\textemdash i.e.,
the case where $\mathfrak{z}$ has only finitely many non-zero $p$-adic
digits. For such a $\mathfrak{z}$, the recursive application of (\ref{eq:System})
will enable us to express $\lim_{N\rightarrow\infty}\chi\left(\left[\mathfrak{z}\right]_{p^{N}}\right)$
as a sum of finitely many terms. Indeed, for any $\mathfrak{z}\in\mathbb{N}_{0}$,
the limit $\lim_{N\rightarrow\infty}\chi\left(\left[\mathfrak{z}\right]_{p^{N}}\right)$
converges in the discrete topology on $\mathbb{Q}$ and any field
extensions thereof, since $\chi\left(\left[\mathfrak{z}\right]_{p^{N}}\right)=\chi\left(\left[\mathfrak{z}\right]_{p^{N+1}}\right)$
for all $N\geq\lambda_{p}\left(\mathfrak{z}\right)$ whenever $\mathfrak{z}\in\mathbb{N}_{0}$.
As such, a frame for interpreting $\lim_{N\rightarrow\infty}\chi\left(\left[\mathfrak{z}\right]_{p^{N}}\right)$
can have \emph{any }topology associated to the $\mathfrak{z}\in\mathbb{N}_{0}$s.

It is worth noting that this freedom of choice of the topology of
convergence of $\lim_{N\rightarrow\infty}\chi\left(\left[\mathfrak{z}\right]_{p^{N}}\right)$
when $\mathfrak{z}\in\mathbb{N}_{0}$ is actually \emph{built into}
the series formula (\ref{eq:rising limit as limit of series}); it
follows entirely by (completely justified) formal manipulations of
(\ref{eq:rising limit as limit of series}). Fixing $\mathfrak{z}\in\mathbb{N}_{0}$,
let $N\geq\lambda_{p}\left(\mathfrak{z}\right)$. Then:
\begin{align*}
\chi\left(\left[\mathfrak{z}\right]_{p^{N}}\right) & =\sum_{n=0}^{N-1}b_{j_{n}}\prod_{j=0}^{p-1}a_{j}^{\#_{p:j}\left(\left[\mathfrak{z}\right]_{p^{n}}\right)}\\
\left(\textrm{pull out terms with }n\geq\lambda_{p}\left(\mathfrak{z}\right)\right); & =\sum_{n=0}^{\lambda_{p}\left(\mathfrak{z}\right)-1}b_{j_{n}}\prod_{j=0}^{p-1}a_{j}^{\#_{p:j}\left(\left[\mathfrak{z}\right]_{p^{n}}\right)}+\sum_{n=\lambda_{p}\left(\mathfrak{z}\right)}^{N-1}b_{j_{n}}\prod_{j=0}^{p-1}a_{j}^{\#_{p:j}\left(\left[\mathfrak{z}\right]_{p^{n}}\right)}
\end{align*}
Here, note that $\left[\mathfrak{z}\right]_{p^{n}}=\mathfrak{z}$
and $j_{n}=0$ for all $n\geq\lambda_{p}\left(\mathfrak{z}\right)$.
Since we assumed $b_{0}=0$, we have:
\begin{equation}
b_{j_{n}}\prod_{j=0}^{p-1}a_{j}^{\#_{p:j}\left(\left[\mathfrak{z}\right]_{p^{n}}\right)}=\underbrace{b_{0}}_{0}\times\prod_{j=0}^{p-1}a_{j}^{\#_{p:j}\left(\mathfrak{z}\right)}=0,\textrm{ }\forall\mathfrak{z}\in\mathbb{N}_{0},\textrm{ }\forall n\geq\lambda_{p}\left(\mathfrak{z}\right)
\end{equation}
and so the tail of the series vanishes, giving:
\begin{equation}
\chi\left(\mathfrak{z}\right)=\chi\left(\left[\mathfrak{z}\right]_{p^{N}}\right)=\sum_{n=0}^{\lambda_{p}\left(\mathfrak{z}\right)-1}b_{j_{n}}\prod_{j=0}^{p-1}a_{j}^{\#_{p:j}\left(\left[\mathfrak{z}\right]_{p^{n}}\right)},\textrm{ }\forall\mathfrak{z}\in\mathbb{N}_{0},\textrm{ }\forall N\geq\lambda_{p}\left(\mathfrak{z}\right)
\end{equation}
This is a ``formulaic'' verification of the fact that $\lim_{N\rightarrow\infty}\chi\left(\left[\mathfrak{z}\right]_{p^{N}}\right)$
converges in the discrete topology for $\mathfrak{z}\in\mathbb{N}_{0}$.

So, in order to find a frame capable of making sense of (\ref{eq:rising limit as limit of series}),
we need only deal with those $p$-adic integers in $\mathbb{Z}_{p}^{\prime}$.
For this, we use a $p$-adic (as opposed to $2$-adic) version of
the \textbf{digit-counting principle} from \textbf{Example \ref{exa:1}}.

Let $q$ be any prime number and let $\mathfrak{z}\in\mathbb{Z}_{p}$.
Then, by the ultrametric structure of $\mathbb{Q}_{q}$, since (\ref{eq:rising limit as limit of series})
expresses $\lim_{N\rightarrow\infty}\chi\left(\left[\mathfrak{z}\right]_{p^{N}}\right)$
as the limit of a sequence of partial sums, it follows that $\lim_{N\rightarrow\infty}\chi\left(\left[\mathfrak{z}\right]_{p^{N}}\right)$
will converge in $\mathbb{Q}_{q}$ if and only if its $n$th term's
$q$-adic absolute value decays to $0$. Examining our formula (\ref{eq:nth term})
for the $n$th term of this series, we see that in order for $q$-adic
convergence to occur, it suffices that:
\begin{equation}
\lim_{n\rightarrow\infty}\prod_{j=0}^{p-1}\left|a_{j}\right|_{q}^{\#_{p:j}\left(\left[\mathfrak{z}\right]_{p^{n}}\right)}\overset{\mathbb{R}}{=}0\label{eq:decay condition}
\end{equation}
where, as indicated, the convergence is in the topology of $\mathbb{R}$.
Our \textbf{$p$-adic} \textbf{digit counting principle }is then the
observation that, given $\mathfrak{z}\in\mathbb{Z}_{p}$, (\ref{eq:decay condition})
occurs whenever there is a single $j\in\left\{ 1,\ldots,p-1\right\} $
such that both of the following occur:

I. $j$ occurs infinitely many times in $\mathfrak{z}$'s sequence
of $p$-adic digits.

II. $a_{j}$ has $q$-adic absolute value $<1$.

Recall that, for all $j\in\left\{ 1,\ldots,p-1\right\} $, the system
of functional equations we started with came with primes $q_{j}$
and integer constants $a_{j}$ so that $\left|a_{j}\right|_{q_{j}}<1$.
Thus, given $\mathfrak{z}\in\mathbb{Z}_{p}$, the digit-counting principle
tells us that $\lim_{N\rightarrow\infty}\chi\left(\left[\mathfrak{z}\right]_{p^{N}}\right)$
will converge in $\mathbb{Q}_{q_{j}}$ whenever $\mathfrak{z}$'s
digit sequence contains infinitely many $j$s. Since every $\mathfrak{z}\in\mathbb{Z}_{p}^{\prime}$
has \emph{at least} one number in $\left\{ 1,\ldots,p-1\right\} $
which occurs in $\mathfrak{z}$'s digit sequence infinitely many times,
every $\mathfrak{z}\in\mathbb{Z}_{p}^{\prime}$ has at least one $j\in\left\{ 1,\ldots,p-1\right\} $
so that $\lim_{N\rightarrow\infty}\chi\left(\left[\mathfrak{z}\right]_{p^{N}}\right)$
converges $q_{j}$-adically. (This $j$ need not be unique, but at
the moment, that is not a problem for us.)

For simplicity's sake, we will define a $\mathbb{Q}$-frame $\mathcal{F}$
for $\lim_{N\rightarrow\infty}\chi\left(\left[\mathfrak{z}\right]_{p^{N}}\right)$
by constructing sets $U_{0},\ldots,U_{p-1}$ which partition $\mathbb{Z}_{p}$,
and then have $\mathcal{F}\left(\mathfrak{z}\right)$ be determined
by which of the $U_{j}$s $\mathfrak{z}$ lives in. Letting $q_{0}\overset{\textrm{def}}{=}1$
and $\mathbb{Q}_{q_{0}}\overset{\textrm{def}}{=}\mathbb{C}$, let
us \emph{define }$U_{0}$ to be $\mathbb{N}_{0}$. We then define
$\mathcal{F}\left(\mathfrak{z}\right)$ to be $\mathbb{Q}_{q_{0}}$
for all $\mathfrak{z}\in U_{0}$. Next, for $j\in\left\{ 0,\ldots,p-2\right\} $,
let $U_{j+1}$ be the set of all $\mathfrak{z}\in\mathbb{Z}_{p}\backslash U_{j}$
with infinitely many $\left(j+1\right)$s in their $p$-adic digits\footnote{Thus, $U_{1}$ is the set of all $p$-adic integers with infinitely
many $1$s digits; $U_{2}$ is the set of all $p$-adic integers with
infinitely many $2$s digits, but only finitely many $1$s digits;
$U_{3}$ is the set of all $p$-adic integers with infinitely many
$3$s digits, but only finitely many $1$s or $2$s digits; etc.}. We then define our $\mathbb{Q}$-frame $\mathcal{F}$ on $\mathbb{Z}_{p}$
by the rule:
\begin{equation}
\mathcal{F}\left(\mathfrak{z}\right)\overset{\textrm{def}}{=}\begin{cases}
\mathbb{Q}_{q_{0}} & \textrm{if }\mathfrak{z}\in U_{0}\\
\mathbb{Q}_{q_{1}} & \textrm{if }\mathfrak{z}\in U_{1}\\
\vdots & \vdots\\
\mathbb{Q}_{q_{p-1}} & \textrm{if }\mathfrak{z}\in U_{p-1}
\end{cases}
\end{equation}
Then, we define a function $\chi:\mathbb{Z}_{p}\rightarrow I\left(\mathcal{F}\right)$
(here, $I\left(\mathcal{F}\right)=\bigcup_{j=0}^{p-1}\mathbb{Q}_{q_{j}}$)
by the rule:
\begin{equation}
\chi\left(\mathfrak{z}\right)\overset{\textrm{def}}{=}\begin{cases}
\lim_{N\rightarrow\infty}\chi\left(\left[\mathfrak{z}\right]_{p^{N}}\right)\textrm{ in }\mathbb{Q}_{q_{j}} & \textrm{if }\mathfrak{z}\in U_{j}\textrm{ for }j\in\mathbb{Z}/p\mathbb{Z}\end{cases}
\end{equation}
As defined, $\chi$ is compatible with $\mathcal{F}$. Note that this
is \emph{not }the only frame we could have chosen to make sense of
$\lim_{N\rightarrow\infty}\chi\left(\left[\mathfrak{z}\right]_{p^{N}}\right)$;
we will address this issue in \textbf{Section \ref{sec:3}}.
\end{example}
With this example in our hands, we can now directly\emph{ }connect
this work to $S_{d;q_{0},\ldots,q_{p-1}}\left(\mathfrak{z}\right)$.
This will be done by demonstrating $S_{d;q_{0},\ldots,q_{p-1}}\left(\mathfrak{z}\right)$
satisfies a $\chi$-type system of functional equations. To accomplish
this, we need two simple results from \cite{my dissertation}, which
we quote in full\textemdash proofs and all. Note that in the proof,
we use \textbf{Iverson bracket }notation; thus the function $n\in\mathbb{Z}\mapsto\left[n\overset{p}{\equiv}m\right]\in\left\{ 0,1\right\} $
outputs $1$ if $n$ is congruent to $m$ mod $p$ and outputs $0$
otherwise.
\begin{lem}
\label{lem:truncation functional equation lemma}Let $p$ be a rational
integer $\geq2$, let $K$ be an abelian group\footnote{In practice, $K$ will be a field.},
written additively, and consider a function $\chi:\mathbb{N}_{0}\rightarrow K$.
Suppose that for $j\in\mathbb{Z}/p\mathbb{Z}$ there are functions
$\Phi_{j}:\mathbb{N}_{0}\times K\rightarrow K$ so that: 
\begin{equation}
\chi\left(pn+j\right)=\Phi_{j}\left(n,\chi\left(n\right)\right),\textrm{ }\forall n\in\mathbb{N}_{0},\textrm{ }\forall j\in\mathbb{Z}/p\mathbb{Z}\label{eq:Relation between truncations and functional equations - Hypothesis}
\end{equation}
Then:
\begin{equation}
\chi\left(\left[pn+j\right]_{p^{N}}\right)=\Phi_{j}\left(\left[n\right]_{p^{N-1}},\chi\left(\left[n\right]_{p^{N-1}}\right)\right),\textrm{ }\forall n\in\mathbb{N}_{0},\textrm{ }\forall j\in\mathbb{Z}/p\mathbb{Z}\label{eq:Relation between truncations and functional equations, version 2}
\end{equation}
\end{lem}
Proof: Fix $N\geq0$, $n\in\mathbb{N}_{0}$, and $j\in\mathbb{Z}/p\mathbb{Z}$.
Then: 
\begin{align*}
\chi\left(\left[pn+j\right]_{p^{N}}\right) & =\sum_{m=0}^{p^{N}-1}\chi\left(m\right)\left[pn+j\overset{p^{N}}{\equiv}m\right]\\
\left(\textrm{split }m\textrm{ mod }p\right); & =\sum_{\ell=0}^{p^{N-1}-1}\sum_{k=0}^{p-1}\chi\left(p\ell+k\right)\left[pn+j\overset{p^{N}}{\equiv}p\ell+k\right]\\
 & =\sum_{\ell=0}^{p^{N-1}-1}\sum_{k=0}^{p-1}\Phi_{k}\left(\ell,\chi\left(\ell\right)\right)\underbrace{\left[n\overset{p^{N-1}}{\equiv}\ell+\frac{k-j}{p}\right]}_{0\textrm{ }\forall n\textrm{ if }k\neq j}\\
 & =\sum_{\ell=0}^{p^{N-1}-1}\Phi_{j}\left(\ell,\chi\left(\ell\right)\right)\left[n\overset{p^{N-1}}{\equiv}\ell\right]\\
 & =\Phi_{j}\left(\left[n\right]_{p^{N-1}},\chi\left(\left[n\right]_{p^{N-1}}\right)\right)
\end{align*}

Q.E.D.

\vphantom{}Using \textbf{Lemma \ref{lem:truncation functional equation lemma}},\textbf{
}we can see how $\#_{p:j}$'s functional equations are altered when
we pre-compose $\#_{p:j}$ with a projection mod $p^{n}$.
\begin{prop}
\label{prop:pound truncation functional equations}Let $p$ be a prime.
Then, for all $N\in\mathbb{N}_{1}$ and all $j\in\mathbb{Z}/p\mathbb{Z}$:
\begin{equation}
\#_{p:j}\left(\left[pn+k\right]_{p^{N}}\right)=\#_{p:j}\left(\left[n\right]_{p^{N-1}}\right)+\left[k\overset{p}{\equiv}j\right],\textrm{ }\forall n\in\mathbb{N}_{0},\textrm{ }\forall k\in\mathbb{Z}/p\mathbb{Z}
\end{equation}
Also: 
\begin{equation}
\lambda_{p}\left(\left[pn+k\right]_{p^{N}}\right)=\lambda_{p}\left(\left[n\right]_{p^{N-1}}\right)+1,\textrm{ }\forall n\in\mathbb{N}_{0},\textrm{ }\forall k\in\mathbb{Z}/p\mathbb{Z}
\end{equation}
\end{prop}
Proof: We will handle $\#_{p:0}$ separately. First, let $j\in\left\{ 1,\ldots,p-1\right\} $,
let $\chi\left(n\right)=\#_{p:j}\left(n\right)$ (resp. $\lambda_{p}\left(n\right)$),
and let: $\Phi_{k}\left(n,\chi\left(n\right)\right)=\chi\left(n\right)+\left[k\overset{p}{\equiv}j\right]$
(resp. $\chi\left(n\right)+1$). Applying \textbf{Lemma \ref{lem:truncation functional equation lemma},}
we conclude:
\begin{equation}
\chi\left(\left[pn+k\right]_{p^{N}}\right)=\Phi_{k}\left(\left[n\right]_{p^{N-1}},\chi\left(\left[n\right]_{p^{N-1}}\right)\right)=\chi\left(\left[n\right]_{p^{N-1}}\right)+\left[k\overset{p}{\equiv}j\right]
\end{equation}
resp.:
\begin{equation}
\chi\left(\left[pn+k\right]_{p^{N}}\right)=\Phi_{k}\left(\left[n\right]_{p^{N-1}},\chi\left(\left[n\right]_{p^{N-1}}\right)\right)=\chi\left(\left[n\right]_{p^{N-1}}\right)+1
\end{equation}
For $\#_{p:0}$, we use (\ref{eq:Definition of number of zeroes})
to write:
\begin{equation}
\#_{p:0}\left(\left[pn+k\right]_{p^{N}}\right)=\lambda_{p}\left(\left[pn+k\right]_{p^{N}}\right)-\sum_{j=1}^{p-1}\#_{p:j}\left(\left[pn+k\right]_{p^{N}}\right)
\end{equation}
and then evaluate the right-hand side using the formulae for $\lambda_{p}$
and $\#_{p:j}$s for non-zero $j$:
\begin{align*}
\#_{p:0}\left(\left[pn+k\right]_{p^{N}}\right) & =\lambda_{p}\left(\left[n\right]_{p^{N-1}}\right)+1-\sum_{j=1}^{p-1}\left(\#_{p:j}\left(\left[n\right]_{p^{N-1}}\right)+\left[k\overset{p}{\equiv}j\right]\right)\\
 & =\underbrace{\lambda_{p}\left(\left[n\right]_{p^{N-1}}\right)-\sum_{j=1}^{p-1}\#_{p:j}\left(\left[n\right]_{p^{N-1}}\right)}_{\#_{p:0}\left(\left[n\right]_{p^{N-1}}\right)}+1-\sum_{j=1}^{p-1}\left[k\overset{p}{\equiv}j\right]
\end{align*}
Noting that:
\begin{equation}
\sum_{j=0}^{p-1}\left[k\overset{p}{\equiv}j\right]=1,\textrm{ }\forall k\in\mathbb{Z}
\end{equation}
(because every integer $k$ is congruent mod $p$ to exactly one $j\in\mathbb{Z}/p\mathbb{Z}$),
we have:
\begin{equation}
1-\sum_{j=1}^{p-1}\left[k\overset{p}{\equiv}j\right]=\left[k\overset{p}{\equiv}0\right]
\end{equation}
and so:
\begin{align*}
\#_{p:0}\left(\left[pn+k\right]_{p^{N}}\right) & =\#_{p:0}\left(\left[n\right]_{p^{N-1}}\right)+1-\sum_{j=1}^{p-1}\left[k\overset{p}{\equiv}j\right]\\
 & =\#_{p:0}\left(\left[n\right]_{p^{N-1}}\right)+\left[k\overset{p}{\equiv}0\right]
\end{align*}

Q.E.D.

\vphantom{}Now the connection between $S_{d;q_{0},\ldots,q_{p-1}}\left(\mathfrak{z}\right)$
and $\chi$.
\begin{prop}
\label{prop:pointwise convergence of generalized S}$S_{d;q_{0},\ldots,q_{p-1}}\left(\mathfrak{z}\right)$
is a formal solution to the system of functional equations:
\begin{equation}
S_{d;q_{0},\ldots,q_{p-1}}\left(p\mathfrak{z}+k\right)=\frac{q_{k}}{d}S_{d;q_{0},\ldots,q_{p-1}}\left(\mathfrak{z}\right)+1,\textrm{ }\forall\mathfrak{z}\in\mathbb{Z}_{p},\textrm{ }\forall k\in\mathbb{Z}/p\mathbb{Z}\label{eq:Generalized S functional equations}
\end{equation}
Moreover, for any $\mathfrak{z}\in\mathbb{Z}_{p}$ and any prime $\ell$
of $\mathbb{Z}$ (including the infinite prime), the following hold:

I. $S_{d;q_{0},\ldots,q_{p-1}}\left(\mathfrak{z}\right)$ converges
in $\mathbb{Q}_{\ell}$ if and only if $S_{d;q_{0},\ldots,q_{p-1}}\left(p\mathfrak{z}+k\right)$
converges in $\mathbb{Q}_{\ell}$ for \emph{at least one} $k\in\mathbb{Z}/p\mathbb{Z}$.

II. If $S_{d;q_{0},\ldots,q_{p-1}}\left(p\mathfrak{z}+k\right)$ converges
in $\mathbb{Q}_{\ell}$ for \emph{at least one} $k\in\mathbb{Z}/p\mathbb{Z}$,
then $S_{d;q_{0},\ldots,q_{p-1}}\left(p\mathfrak{z}+k\right)$ converges
in $\mathbb{Q}_{\ell}$ \emph{for all $k\in\mathbb{Z}/p\mathbb{Z}$.}
\end{prop}
Proof: Proceeding formally, let $\mathfrak{z}\in\mathbb{Z}_{p}$,
$k\in\mathbb{Z}/p\mathbb{Z}$, and apply \textbf{Proposition \ref{prop:pound truncation functional equations}}:
\begin{align*}
S_{d;q_{0},\ldots,q_{p-1}}\left(p\mathfrak{z}+k\right) & =\sum_{n=0}^{\infty}\frac{1}{d^{n}}\prod_{j=0}^{p-1}q_{j}^{\#_{p:j}\left(\left[p\mathfrak{z}+k\right]_{p^{n}}\right)}\\
\left(\left[\mathfrak{y}\right]_{p^{0}}=0\textrm{ }\forall\mathfrak{y}\in\mathbb{Z}_{p}\right); & =\prod_{j=0}^{p-1}q_{j}^{\overbrace{\#_{p:j}\left(0\right)}^{0}}+\sum_{n=0}^{\infty}\frac{1}{d^{n}}\prod_{j=0}^{p-1}q_{j}^{\#_{p:j}\left(\left[p\mathfrak{z}+k\right]_{p^{n}}\right)}\\
\left(\textrm{\textbf{Prop. }\textbf{\ref{exa:5}}}\right); & =1+\sum_{n=1}^{\infty}\frac{1}{d^{n}}\prod_{j=0}^{p-1}q_{j}^{\left(\#_{p:j}\left(\left[\mathfrak{z}\right]_{p^{n-1}}\right)+\left[j\overset{p}{\equiv}k\right]\right)}\\
\left(\textrm{re-index}\right); & =1+\frac{q_{k}}{d}\sum_{n=0}^{\infty}\frac{1}{d^{n}}\prod_{j=0}^{p-1}q_{j}^{\#_{p:j}\left(\left[\mathfrak{z}\right]_{p^{n}}\right)}\\
 & =1+\frac{q_{k}}{d}S_{d;q_{0},\ldots,q_{p-1}}\left(\mathfrak{z}\right)
\end{align*}
This proves (\ref{eq:Generalized S functional equations}).

For (I) and (II), fix $\mathfrak{z}\in\mathbb{Z}_{p}$ and suppose
$\ell$ is a prime (possibly infinite) so that $S_{d;q_{0},\ldots,q_{p-1}}\left(\mathfrak{z}\right)$
converges in $\mathbb{Q}_{\ell}$. Then, (\ref{eq:Generalized S functional equations})
shows that for any $k\in\mathbb{Z}/p\mathbb{Z}$, the value of $S_{d;q_{0},\ldots,q_{p-1}}\left(p\mathfrak{z}+k\right)$
is given by applying the affine linear map $\phi_{k}:\mathbb{Q}_{\ell}\rightarrow\mathbb{Q}_{\ell}$
defined by:
\begin{equation}
\phi_{k}\left(\mathfrak{x}\right)\overset{\textrm{def}}{=}\frac{q_{k}\mathfrak{x}}{d}+1,\textrm{ }\forall\mathfrak{x}\in\mathbb{Q}_{\ell}
\end{equation}
to $S_{d;q_{0},\ldots,q_{p-1}}\left(\mathfrak{z}\right)$. Since $\phi_{k}$
is a homeomorphism of $\mathbb{Q}_{\ell}$, it sends the $\mathbb{Q}_{\ell}$-convergent
series $S_{d;q_{0},\ldots,q_{p-1}}\left(\mathfrak{z}\right)$ to the
$\mathbb{Q}_{\ell}$-convergent series $S_{d;q_{0},\ldots,q_{p-1}}\left(p\mathfrak{z}+k\right)$,
and its inverse $\phi_{k}^{-1}$ does the opposite. Moreover, for
any $j,k\in\mathbb{Z}/p\mathbb{Z}$, $\phi_{j}\circ\phi_{k}^{-1}$
is a homeomorphism of $\mathbb{Q}_{\ell}$ which sends the $\mathbb{Q}_{\ell}$-convergent
$S_{d;q_{0},\ldots,q_{p-1}}\left(p\mathfrak{z}+k\right)$ series to
the the $\mathbb{Q}_{\ell}$-convergent series $S_{d;q_{0},\ldots,q_{p-1}}\left(p\mathfrak{z}+j\right)$.

Q.E.D.

\vphantom{}In this viewpoint, we can regard $\mathcal{F}$-series
as natural series representations for solutions to these systems of
affine linear functional equations. In \cite{my dissertation,first blog post paper},
these functions arose as a result of considering solutions $\mathbb{N}_{0}\rightarrow\mathbb{Q}$
to the functional equations and then extending them to functions of
a $p$-adic variable. In all cases, the fact that affine linear maps
like $\phi_{k}$ are homeomorphisms of $\mathbb{Q}_{\ell}$ plays
a pivotal role in \cite{my dissertation,first blog post paper}'s
arguments. The reason this all works is as follows. Given a $K$-frame
$\mathcal{F}$, let $\alpha,\beta\in K$ be constants so that the
affine linear map:
\begin{align*}
\phi & :K\rightarrow K\\
\phi\left(x\right) & \overset{\textrm{def}}{=}\alpha x+\beta,\textrm{ }\forall x\in K
\end{align*}
is invertible on $K$. Crucially, not only is $\phi$ a homeomorphism
of $K$, it also extends to be a homeomorphism of any \emph{completion}
of $K$. By a judicious selection of such an $\phi$, we can then
conjugate $S_{d;q_{0},\ldots,q_{p-1}}\left(\mathfrak{z}\right)$'s
functional equations however we please, as the next example shows.
\begin{example}
\label{exa:6}Choose a $\mathbb{Q}$-frame $\mathcal{F}$ on $\mathbb{Z}_{p}$
compatible with $S_{d;q_{0},\ldots,q_{p-1}}$ (i.e., $S_{d;q_{0},\ldots,q_{p-1}}\left(\mathfrak{z}\right)$
converges in $\mathcal{F}\left(\mathfrak{z}\right)$ for every $\mathfrak{z}\in\mathbb{Z}_{p}$).
Then, letting $\alpha,\beta\in\mathbb{Q}$ (or, in $\mathbb{F}$,
if $\mathbb{F}$ is a completion of $\mathbb{Q}$ to which $S_{d;q_{0},\ldots,q_{p-1}}\left(\mathfrak{z}\right)$
extends) we get a compatible function $f:\mathbb{Z}_{p}\rightarrow I\left(\mathcal{F}\right)$
by defining:
\begin{equation}
f\left(\mathfrak{z}\right)\overset{\textrm{def}}{=}\alpha S_{d;q_{0},\ldots,q_{p-1}}\left(\mathfrak{z}\right)+\beta
\end{equation}
Note that:
\begin{align*}
f\left(p\mathfrak{z}+j\right) & =\alpha S_{d;q_{0},\ldots,q_{p-1}}\left(p\mathfrak{z}+j\right)+\beta\\
 & =\alpha\left(1+\frac{q_{j}}{d}S_{d;q_{0},\ldots,q_{p-1}}\left(\mathfrak{z}\right)\right)+\beta\\
 & =\frac{\alpha q_{j}}{d}S_{d;q_{0},\ldots,q_{p-1}}\left(\mathfrak{z}\right)+\alpha+\beta\\
\left(\pm\frac{\beta q_{j}}{d}\right); & =\frac{q_{j}}{d}\left(\alpha S_{d;q_{0},\ldots,q_{p-1}}\left(\mathfrak{z}\right)+\beta\right)+\alpha+\beta-\frac{\beta q_{j}}{d}\\
 & =\frac{q_{j}}{d}f\left(\mathfrak{z}\right)+\alpha+\beta\left(1-\frac{q_{j}}{d}\right)
\end{align*}
As such, for any $\alpha,\beta\in\mathbb{Q}$ with $\alpha\neq0$,
$f$ satisfies the system of functional equations:
\begin{equation}
f\left(p\mathfrak{z}+j\right)=\frac{q_{j}}{d}f\left(\mathfrak{z}\right)+\gamma_{j}\left(\alpha,\beta\right)
\end{equation}
where:
\begin{equation}
\gamma_{j}\left(\alpha,\beta\right)\overset{\textrm{def}}{=}\alpha+\beta\left(1-\frac{q_{j}}{d}\right)
\end{equation}
As a particular case, pick $p=d=2$, $q_{0}=1$, and $q_{1}=3$. Then
$S_{d;q_{0},\ldots,q_{p-1}}\left(\mathfrak{z}\right)$ becomes:
\begin{equation}
S_{2;1,3}\left(\mathfrak{z}\right)=\sum_{n=0}^{\infty}\frac{3^{\#_{1}\left(\left[\mathfrak{z}\right]_{2^{n}}\right)}}{2^{n}}
\end{equation}
which is our old friend $S_{2,3}$. Using the standard $\left(2,3\right)$-adic
frame to realize this series as $X_{2,3}$, \textbf{Proposition \ref{prop:pointwise convergence of generalized S}}'s
functional equations yield the system:
\begin{align}
X_{2,3}\left(2\mathfrak{z}\right) & =\frac{1}{2}X_{2,3}\left(\mathfrak{z}\right)+1\\
X_{2,3}\left(2\mathfrak{z}+1\right) & =\frac{3}{2}X_{2,3}\left(\mathfrak{z}\right)+1
\end{align}
Since $X_{2,3}$ extends to $\mathbb{Q}_{3}$ with $\mathcal{F}_{2,3}$,
so the above functional equations hold over $\mathbb{Q}_{3}$ for
all $\mathfrak{z}\in\mathbb{Z}_{2}$. Setting:
\begin{equation}
\chi_{3}\left(\mathfrak{z}\right)\overset{\textrm{def}}{=}\frac{X_{2,3}\left(\mathfrak{z}\right)-2}{4}\label{eq:Chi_3}
\end{equation}
we have that $\chi_{3}:\mathbb{Z}_{2}\rightarrow\mathbb{Z}_{3}$ admits
the $\mathcal{F}$-series representation:
\begin{equation}
\chi_{3}\left(\mathfrak{z}\right)\overset{\mathcal{F}_{2,3}}{=}-\frac{1}{2}+\frac{1}{4}\sum_{n=0}^{\infty}\frac{3^{\#_{1}\left(\left[\mathfrak{z}\right]_{2^{n}}\right)}}{2^{n}}\label{eq:F-series for Chi_3}
\end{equation}
and satisfies the functional equations:
\begin{align}
\chi_{3}\left(2\mathfrak{z}\right) & =\frac{1}{2}\chi_{3}\left(\mathfrak{z}\right)\label{eq:Even branch functional equation for Chi_3}\\
\chi_{3}\left(2\mathfrak{z}+1\right) & =\frac{3\chi_{3}\left(\mathfrak{z}\right)+1}{2}
\end{align}
for all $\mathfrak{z}\in\mathbb{Z}_{2}$. Not coincidentally, the
affine linear maps on the right-hand side\textemdash $x\mapsto x/2$
and $x\mapsto\left(3x+1\right)/2$\textemdash are the branches of
the \textbf{Shortened Collatz Map}, $T_{3}:\mathbb{Z}\rightarrow\mathbb{Z}$
defined by:
\begin{equation}
T_{3}\left(n\right)\overset{\textrm{def}}{=}\begin{cases}
\frac{n}{2} & \textrm{if }n\overset{2}{\equiv}0\\
\frac{3n+1}{2} & \textrm{if }n\overset{2}{\equiv}1
\end{cases}
\end{equation}
The infamous \textbf{Collatz Conjecture }is equivalent to the assertion
that every positive integer is eventually mapped to $1$ under repeated
iteration by $T_{3}$. The author originally derived the $\mathcal{F}$-series
representation (\ref{eq:Chi_3}) in \cite{my dissertation} by proving
that $\chi_{3}$ possesses a $\left(2,3\right)$-adic Fourier-transform.
This is a complete novelty. Prior to \cite{my dissertation}, the
only $\left(2,3\right)$-adic functions (or, more generally, $\left(p,q\right)$-adic
functions $\mathbb{Z}_{p}\rightarrow\mathbb{C}_{q}$) known to have
Fourier transforms were the continuous functions, but $\chi_{3}$
is not $\left(2,3\right)$-adically continuous. Indeed, (\ref{eq:Even branch functional equation for Chi_3})
shows that $\chi_{3}\left(2^{n-1}\right)=1/2^{n}$, and so, if $\chi_{3}$
were $\left(2,3\right)$-adically continuous, the convergence of $2^{n-1}$
to $0$ in $\mathbb{Z}_{2}$ as $n\rightarrow\infty$ would force
$1/2^{n}$ to converge to $\chi_{3}\left(0\right)$ in $\mathbb{Z}_{3}$.
However, $\chi_{3}\left(0\right)=0$, and $\left\{ 1/2^{n}\right\} _{n\geq1}$
does not converge $3$-adically!
\end{example}
Frames naturally arise in Fourier analytic investigations of $\chi_{3}$,
$X_{p,q}$, and other $\mathcal{F}$-series. By realizing certain
$\left(p,q\right)$-adic functions as measures, the capabilities of
$\left(p,q\right)$-adic analysis can be expanded into areas previously
believed to lie outside of its purview. At the risk of making a long
paper slightly longer, it is worth briefly introducing just enough
of the pertinent terminology so that we can consider an example of
potentially immense import to the future development of $\left(p,q\right)$-adic
analysis as a subspecialty of non-archimedean/ultrametric analysis.
\begin{example}
\label{exa:7}Following Jordan Bell's exposition (\cite{Bell's notes})
of \cite{Automorphic Representations}'s exposition of methods of
Fourier analysis for functions $\mathbb{Z}_{p}\rightarrow\mathbb{C}$,
we write $\hat{\mathbb{Z}}_{p}$ to denote $\mathbb{Z}\left[1/p\right]/\mathbb{Z}$,
the Pontryagin Dual of $\mathbb{Z}_{p}$, identified here with the
set of rational numbers in $\left[0,1\right)$ with powers of $p$
as their denominators, made into a group by way of addition modulo
$1$. We write $\left\{ \cdot\right\} _{p}:\mathbb{Q}_{p}\rightarrow\hat{\mathbb{Z}}_{p}$
to denote the $p$-adic fractional part. Every unitary character $\mathbb{Z}_{p}\rightarrow\mathbb{T}$
is of the form $\mathfrak{z}\in\mathbb{Z}_{p}\mapsto e^{2\pi i\left\{ t\mathfrak{z}\right\} _{p}}\in\mathbb{C}$
for some unique $t\in\hat{\mathbb{Z}}_{p}$. Multiplication of $t$
by $\mathfrak{z}$ is justified by viewing $\hat{\mathbb{Z}}_{p}$
as the subset of $\mathbb{Q}_{p}$ containing all $p$-adic numbers
whose integer part (their $\mathbb{Z}_{p}$ component) is $0$. We
write $v_{p}$ to denote the $p$-adic valuation (so that $\left|\cdot\right|_{p}=p^{-v_{p}\left(\cdot\right)}$
is the $p$-adic absolute value); recall $v_{p}\left(0\right)\overset{\textrm{def}}{=}+\infty$.

Now, let $p$ be a prime, and let $c$ be an algebraic constant\textemdash say,
$c\in\mathbb{\overline{Q}}\backslash\left\{ 0,1,-1,p\right\} $\textemdash let
$\mathfrak{z}\in\mathbb{Z}_{p}$, and let $N$ be an integer $\geq0$.
Then, a simple computation shows that:

\begin{equation}
\sum_{0<\left|t\right|_{p}\leq p^{N}}c^{v_{p}\left(t\right)}e^{2\pi i\left\{ t\mathfrak{z}\right\} _{p}}=\begin{cases}
\frac{p-1}{p-c}\left(\left(\frac{p}{c}\right)^{N}-1\right) & \textrm{if }\mathfrak{z}=0\\
\left(\frac{p}{c}\right)^{N}\left[N\leq v_{p}\left(\mathfrak{z}\right)\right]+\frac{c-1}{p-c}\left(\frac{p}{c}\right)^{1+\min\left\{ N-1,v_{p}\left(\mathfrak{z}\right)\right\} }-\frac{p-1}{p-c} & \textrm{else}
\end{cases}\label{eq:Pre-Second Main Identity}
\end{equation}
Here, $\left[N\leq v_{p}\left(\mathfrak{z}\right)\right]$ is an Iverson
bracket, which evaluates to $1$ for all $N$ and $\mathfrak{z}$
for which $N\leq v_{p}\left(\mathfrak{z}\right)$ is true, and which
evaluates to $0$ otherwise.

For fixed $\mathfrak{z}\in\mathbb{Z}_{p}\backslash\left\{ 0\right\} $,
observe that $v_{p}\left(\mathfrak{z}\right)$ is a non-negative integer.
As such, for said $\mathfrak{z}$, we have:
\begin{align}
\lim_{N\rightarrow\infty}\left[N\leq v_{p}\left(\mathfrak{z}\right)\right] & \overset{\mathbb{Z}}{=}0\\
\lim_{N\rightarrow\infty}\min\left\{ N-1,v_{p}\left(\mathfrak{z}\right)\right\}  & \overset{\mathbb{Z}}{=}v_{p}\left(\mathfrak{z}\right)
\end{align}
where $\overset{\mathbb{Z}}{=}$ means the limits converge in the
topology of $\mathbb{Z}$, which is to say, the discrete topology.
Consequently, taking limits, we have that, if $c\neq p$ and $\mathfrak{z}\in\mathbb{Z}_{p}\backslash\left\{ 0\right\} $:
\begin{align*}
\sum_{t\in\mathbb{Z}_{p}\backslash\left\{ 0\right\} }c^{v_{p}\left(t\right)}e^{2\pi i\left\{ t\mathfrak{z}\right\} _{p}} & \overset{\textrm{def}}{=}\lim_{N\rightarrow\infty}\sum_{0<\left|t\right|_{p}\leq p^{N}}c^{v_{p}\left(t\right)}e^{2\pi i\left\{ t\mathfrak{z}\right\} _{p}}\\
 & \overset{\mathbb{Q}}{=}\frac{c-1}{p-c}\left(\frac{p}{c}\right)^{1+v_{p}\left(\mathfrak{z}\right)}-\frac{p-1}{p-c}
\end{align*}
where the topology of convergence is that of the discrete topology
on $\mathbb{Q}$, and is point-wise with respect to $\mathfrak{z}\in\mathbb{Z}_{p}\backslash\left\{ 0\right\} $.
If we let $c=p^{\alpha+1}$, for some constant $\alpha$\textemdash say,
$\alpha\in\mathbb{Q}\backslash\left\{ -1\right\} $\textemdash observe
that:

\begin{equation}
\sum_{t\in\mathbb{Z}_{p}\backslash\left\{ 0\right\} }p^{\left(\alpha+1\right)v_{p}\left(t\right)}e^{2\pi i\left\{ t\mathfrak{z}\right\} _{p}}\overset{\mathbb{Q}}{=}\frac{p^{-\alpha}\left(p^{\alpha+1}-1\right)}{p-p^{\alpha+1}}\left|\mathfrak{z}\right|_{p}^{\alpha}-\frac{p-1}{p-p^{\alpha+1}},\textrm{ }\forall\mathfrak{z}\in\mathbb{Z}_{p}\backslash\left\{ 0\right\} 
\end{equation}
and hence (since $p^{\left(\alpha+1\right)v_{p}\left(t\right)}=\left|t\right|_{p}^{-\alpha-1}$):
\begin{equation}
\left|\mathfrak{z}\right|_{p}^{\alpha}\overset{\mathbb{Q}}{=}\frac{p^{\alpha}\left(p-1\right)}{p^{\alpha+1}-1}+\frac{p^{\alpha}\left(p-p^{\alpha+1}\right)}{p^{\alpha+1}-1}\sum_{t\in\mathbb{Z}_{p}\backslash\left\{ 0\right\} }\left|t\right|_{p}^{-\alpha-1}e^{2\pi i\left\{ t\mathfrak{z}\right\} _{p}},\textrm{ }\forall\mathfrak{z}\in\mathbb{Z}_{p}\backslash\left\{ 0\right\} 
\end{equation}
Consequently, the map $\hat{\mu}_{\alpha}:\hat{\mathbb{Z}}_{p}\rightarrow\overline{\mathbb{Q}}$
defined by:
\begin{equation}
\hat{\mu}_{\alpha}\left(t\right)\overset{\textrm{def}}{=}\begin{cases}
\frac{p^{\alpha}\left(p-1\right)}{p^{\alpha+1}-1} & \textrm{if }t=0\\
\frac{p^{\alpha}\left(p-p^{\alpha+1}\right)}{p^{\alpha+1}-1}\left|t\right|_{p}^{-\alpha-1} & \textrm{else}
\end{cases}
\end{equation}
for any $\alpha\in\mathbb{Q}\backslash\left\{ 0\right\} $ satisfies
the property that the function $\mu_{\alpha,N}:\mathbb{Z}_{p}\rightarrow\mathbb{Q}\left(p^{\alpha}\right)$
defined by:
\begin{equation}
\mu_{\alpha,N}\left(\mathfrak{z}\right)\overset{\textrm{def}}{=}\sum_{0<\left|t\right|_{p}\leq p^{N}}\hat{\mu}_{\alpha}\left(t\right)e^{2\pi i\left\{ t\mathfrak{z}\right\} _{p}}
\end{equation}
is locally constant and $\mathbb{Q}\left(p^{\alpha}\right)$-valued
on $\mathbb{Z}_{p}$ for all $N\in\mathbb{N}_{1}$. Moreover, for
every $\mathfrak{z}\in\mathbb{Z}_{p}\backslash\left\{ 0\right\} $,
$\mu_{\alpha,N}\left(\mathfrak{z}\right)$ converges to $\left|\mathfrak{z}\right|_{p}^{\alpha}$
in the discrete topology\footnote{Note that we could just as easily say that the sequence $\left\{ \mu_{\alpha,N}\left(\mathfrak{z}\right)\right\} _{N\geq1}$
converges to $\left|\mathfrak{z}\right|_{p}^{\alpha}$ on $\mathbb{Z}_{p}\backslash\left\{ 0\right\} $
with respect to the frame that assigns the field $\overline{\mathbb{Q}}$
(equipped with the discrete topology) to each $\mathfrak{z}\in\mathbb{Z}_{p}\backslash\left\{ 0\right\} $.
Note also that this notion of convergence is compatible with that
of any other frame $\mathcal{F}$ on $\mathbb{Z}_{p}\backslash\left\{ 0\right\} $,
so long as $\mathcal{F}\left(\mathfrak{z}\right)$ is an algebraically
closed field extension of $\mathbb{Q}$ for all $\mathfrak{z}\in\mathbb{Z}_{p}\backslash\left\{ 0\right\} $;
any such field extension necessarily contains an isomorphic copy of
$\overline{\mathbb{Q}}$. Admittedly, certain Galois-theoretic technicalities
arise, but these can be resolved by adopting a simple standard convention
of how to write the $p$-adic digits of a $q$-power root of unity,
where $q$ is a prime divisor of $p-1$. For brevity's sake, a thorough
discussion of this issue, as well as $\left(p,q\right)$-adic Fourier
analysis in general, will be deferred for a later publication, though
the reader can always consult \cite{my dissertation} for all the
pertinent details.}.

So, letting $q$ be a prime distinct from $p$ (including the infinite
prime), after embedding $\mathbb{Q}\left(p^{\alpha}\right)$ in $\mathbb{C}_{q}$,
we realize $\mu_{\alpha,N}:\mathbb{Z}_{p}\rightarrow\mathbb{C}_{q}$
as a locally constant\textemdash and hence, continuous\textemdash function.
Consequently, it has a Fourier transform $\hat{\mu}_{\alpha,N}:\hat{\mathbb{Z}}_{p}\rightarrow\mathbb{C}_{q}$.
Noting that:
\begin{equation}
\hat{\mu}_{\alpha,N}\left(t\right)=\begin{cases}
\hat{\mu}_{\alpha}\left(t\right) & \textrm{if }\left|t\right|_{p}\leq p^{N}\\
0 & \textrm{else}
\end{cases}
\end{equation}
we have that the Fourier transforms of the $\mu_{\alpha,N}$s converge
point-wise in the discrete topology on $\mathbb{C}_{q}$ (and hence,
in \emph{any }topology on $\mathbb{C}_{q}$) to $\hat{\mu}_{\alpha}$.
As long as we choose $q$ to be a prime so that $\hat{\mu}_{\alpha}\left(t\right)$
is bounded in $q$-adic absolute value, we can whip up a $\left(p,q\right)$-adic
measure by miming the Parseval-Plancherel formula and writing:
\begin{equation}
f\mapsto\sum_{t\in\hat{\mathbb{Z}}_{p}}\hat{f}\left(-t\right)\hat{\mu}_{\alpha}\left(t\right)
\end{equation}
where $f:\mathbb{Z}_{p}\rightarrow\mathbb{C}_{q}$ is any continuous
$\left(p,q\right)$-adic function. We then formally denote the right-hand
side as:
\begin{equation}
\int_{\mathbb{Z}_{p}}f\left(\mathfrak{z}\right)\left|\mathfrak{z}\right|_{p}^{\alpha}d\mathfrak{z}
\end{equation}
where, here, $\left|\mathfrak{z}\right|_{p}^{\alpha}d\mathfrak{z}$
is defined to be the unique $\left(p,q\right)$-adic measure possessing
$\hat{\mu}_{\alpha}\left(t\right)$ as its Fourier-Stieltjes transform.
Note that, for $\alpha\in\mathbb{Q}\backslash\left\{ 0,-1\right\} $,
$\mathfrak{z}\mapsto\left|\mathfrak{z}\right|_{p}^{\alpha}$ is \emph{not}
$\left(p,q\right)$-adically continuous when $p\neq q$, and thus
is not $\left(p,q\right)$-adically integrable in the classical sense
(e.g., as a limit of Riemann sums). Nevertheless, by \emph{defining}
$\hat{\mu}_{\alpha}$ as the Fourier transform of $\left|\mathfrak{z}\right|_{p}^{\alpha}$,
we can realize $\left|\mathfrak{z}\right|_{p}^{\alpha}$ as the $\left(p,q\right)$-adic
measure\footnote{Or, perhaps, as the ``Radon-Nikodym derivative'' of the measure
$\left|\mathfrak{z}\right|_{p}^{\alpha}d\mathfrak{z}$, for an appropriate
sense of ``Radon-Nikodym derivative''; the issue is not entirely
straightforward in the $\left(p,q\right)$-adic context, particularly
when $p\neq q$; see for example, \cite{Schikhof - Radon-Nikodym}.} $\left|\mathfrak{z}\right|_{p}^{\alpha}d\mathfrak{z}$, and thereby
integrate it against any continuous $\left(p,q\right)$-adic function.

This recipe is of interest, seeing as integration against $\left|\mathfrak{z}\right|_{p}^{\alpha}$
is the basis for the theory of distributional differentiation of complex-valued
(generalized) functions of a $p$-adic variable (see \cite{First 30 years of p-adic mathematical physics,Vladimirov - the big paper about complex-valued distributions over the p-adics}
for details), a construction which is at the heart of $p$-adic approaches
to quantum mechanics, among other applications of note. This strongly
suggests that we should be able to extend the distributional notion
of differentiation given by \cite{Vladimirov - the big paper about complex-valued distributions over the p-adics}
to $\left(p,q\right)$-adic functions, thereby opening the door to
entirely new forms of analysis. Particularly tantalizing is the possibility
of adèlic-valued analysis, where functions are studied from all possible
metric completions of $\mathbb{Q}$ \emph{simultaneously}.
\end{example}

\section{\label{sec:3}Frames of Higher Degree; or, an Adèlic Mindset}

An astute reader will have noticed that a certain arbitrary decision
was made in constructing the frame used to make sense of the function
$f$ from \textbf{Example \ref{exa:5}}. There were other frames we
could have chosen at the end was \emph{not }the only frame we \emph{could}
have chosen.
\begin{example}
\label{exa:8}Consider the $\mathcal{F}$-series:
\begin{equation}
S_{2;3,5}\left(\mathfrak{z}\right)\overset{\textrm{def}}{=}\sum_{n=0}^{\infty}\frac{3^{\#_{3:1}\left(\left[\mathfrak{z}\right]_{3^{n}}\right)}5^{\#_{3:2}\left(\left[\mathfrak{z}\right]_{3^{n}}\right)}}{2^{n}}\label{eq:Definition of S 2 3 5}
\end{equation}
for $\mathfrak{z}\in\mathbb{Z}_{3}$, with the associated system of
functional equations:
\begin{align}
\chi\left(3n\right) & =\frac{\chi\left(n\right)}{2}+1\\
\chi\left(3n+1\right) & =\frac{3\chi\left(n\right)}{2}+1\\
\chi\left(3n+2\right) & =\frac{5\chi\left(n\right)}{2}+1
\end{align}
By \textbf{Proposition \ref{prop:pointwise convergence of generalized S}},
$S_{2;3,5}\left(\mathfrak{z}\right)$ will converge:

\textbullet{} In the topology of $\mathbb{Z}_{3}$ if (and only if)
$\mathfrak{z}$ has infinitely many $1$s in its $3$-adic digits.

\textbullet{} In the topology of $\mathbb{Z}_{5}$ if (and only if)
$\mathfrak{z}$ has infinitely many $2$s in its $3$-adic digits.

\textbullet{} In the \emph{both} topology of $\mathbb{Z}_{3}$ \emph{and
}$\mathbb{Z}_{5}$ if (and only if) $\mathfrak{z}$ has infinitely
many $1$s \emph{and }infinitely many $2$s in its $3$-adic digits.

\textbullet{} In the topology of $\mathbb{R}$ whenever $\mathfrak{z}\in\mathbb{N}_{0}$.

Thus, in using a $\mathbb{Q}$-frame to make sense of $S_{2;3,5}\left(\mathfrak{z}\right)$'s
infinite series representation for inputs in $\mathbb{Z}_{3}^{\prime}$,
we can choose any $\mathbb{Q}$-frame which sends $\mathfrak{z}$s
with infinitely many $1$s (resp. infinitely many $2$s) to $\mathbb{Z}_{3}$
(resp. $\mathbb{Z}_{5}$). We can get different choices of frames
based on whether we send a given $\mathfrak{z}$ with infinitely many
$1$s \emph{and }$2$s to either $\mathbb{Z}_{3}$ or $\mathbb{Z}_{5}$.

We can say even more by using geometric series universality (\textbf{Theorem
\ref{thm:geometric series universality}}). By proving an analogue
of \textbf{Lemma \ref{lem:1}} for $S_{2;3,5}\left(\mathfrak{z}\right)$,
it is easy to show that $S_{2;3,5}\left(\mathfrak{z}\right)$ will
be a geometric series whenever $\mathfrak{z}\in\mathbb{Z}_{3}\cap\mathbb{Q}$.
Let $r\left(\mathfrak{z}\right)$ denote the common ratio of this
series; i.e., so that $S_{2;3,5}\left(\mathfrak{z}\right)$ is a formal
series of the form:
\begin{equation}
S_{2;3,5}\left(\mathfrak{z}\right)=A+B\sum_{n=0}^{\infty}\left(r\left(\mathfrak{z}\right)\right)^{n}
\end{equation}
for some rational numbers $A,B$ and $r\left(\mathfrak{z}\right)$.

By geometric series universality (\textbf{Theorem \ref{thm:geometric series universality}}),
for any $\mathfrak{z}\in\mathbb{Z}_{3}\cap\mathbb{Q}$ and any prime
$p$ (including the infinite prime) such that the $\left|r\left(\mathfrak{z}\right)\right|_{p}<1$
(so that $S_{2;3,5}\left(\mathfrak{z}\right)$ converges in $\mathbb{Q}_{p}$),
if there is a (possibly infinite) prime $\ell\neq p$ so that $\left|r\left(\mathfrak{z}\right)\right|_{\ell}<1$,
then $S_{2;3,5}\left(\mathfrak{z}\right)$ will converge in $\mathbb{Q}_{\ell}$
as well, and the rational number to which it converges in $\mathbb{Q}_{\ell}$
will be the same as the rational number to which it converges in $\mathbb{Q}_{p}$.
So, for example, we could make sense of $S_{2;3,5}\left(\mathfrak{z}\right)$
using a $\mathbb{Q}$-frame $\mathcal{F}$ which assigns the topology
of $\mathbb{R}$ to all $\mathfrak{z}\in\mathbb{Z}_{3}\cap\mathbb{Q}$
for which $S_{2;3,5}\left(\mathfrak{z}\right)$ converges in $\mathbb{R}$.
Fixing such a $\mathfrak{z}$\textemdash let us call it $\mathfrak{z}_{0}$\textemdash observe
that for any other $\mathbb{Q}$-frame $\mathcal{F}^{\prime}$ so
that $S_{2;3,5}\in C\left(\mathcal{F}^{\prime}\right)$, the value
to which $S_{2;3,5}\left(\mathfrak{z}_{0}\right)$ converges in $\mathcal{F}\left(\mathfrak{z}_{0}\right)$
(a.k.a. in $\mathbb{R}$) will be the same rational number as the
one that $S_{2;3,5}\left(\mathfrak{z}_{0}\right)$ converges to in
$\mathcal{F}^{\prime}\left(\mathfrak{z}_{0}\right)$.
\end{example}
The phenomenon observed in \textbf{Example \ref{exa:8}} also occurred
in our analysis of $S_{7,6}\left(\mathfrak{z}\right)$ in \textbf{Example
\ref{exa:2}}. The functions $X_{7,6}^{\left(\infty\right)}\left(\mathfrak{z}\right)$,
$X_{7,6}^{\left(2\right)}\left(\mathfrak{z}\right)$, and $X_{7,6}^{\left(3\right)}\left(\mathfrak{z}\right)$
obtained in \textbf{Example \ref{exa:8}} can be viewed as the result
of interpreting $S_{7,6}\left(\mathfrak{z}\right)$ the standard $\left(2,\infty\right)$-adic,
$\left(2,2\right)$-adic, and $\left(2,3\right)$-adic frames, respectively.
The lesson here is that, unlike classical studies of infinite series,
$\mathcal{F}$-series need to be examined from the vantage of more
than one topology of convergence, sometimes even multiple topologies
at a single point. Of course, this is an incredibly unusual matter
to consider, and, as such, it is natural to wonder if there is anything
to be gained by considering different topologies. The answer to this
question is ``yes''. A basic example of a useful application of
this viewpoint comes from a straightforward generalization of one
of the principal results of \cite{my dissertation}. First, however,
a limited definition:
\begin{defn}[Degree of a frame of dimension $1$]
\label{def:degree of a 1-dim frame}Let $\mathcal{F}$ be a $K$-frame
on $X$ of dimension\footnote{The reader can just assume the frame $\mathcal{F}$ is of the kind
that we have currently been working with. The definition of dimension
is given later on in this section.} $1$. Then, for every $x\in X$, there is a place $v$ of $K$ so
that $\mathcal{F}\left(x\right)=K_{v}$. Now, let: 
\begin{equation}
V_{\mathcal{F}}\overset{\textrm{def}}{=}\left\{ v\in V_{K}:\mathcal{F}\left(x\right)=K_{v}\textrm{ for some }x\in X\right\} \label{eq:places of F}
\end{equation}
That is, $V_{\mathcal{F}}$ is the set of all places of $K$ that
occur in $\mathcal{F}$. We call $V_{\mathcal{F}}$ the set of \textbf{places
of $\mathcal{F}$}.\textbf{ }The \textbf{degree }of $\mathcal{F}$,
denoted $\deg\mathcal{F}$, is defined as the cardinality of $V_{\mathcal{F}}$.

Additionally, regardless of the dimension of a frame, we say $\mathcal{F}$
is:

I. ``\textbf{of lower degree}'' if \textbf{$\deg\mathcal{F}\leq2$}.

II. ``\textbf{of higher degree}'' if $\deg\mathcal{F}\geq3$.
\end{defn}
\begin{example}
$\mathcal{F}$ is a frame of degree $1$ if and only if there is a
single completion $K^{\prime}$ of $K$ so that $\mathcal{F}\left(x\right)=K^{\prime}$
for all $x\in X$. The standard $\left(p,q\right)$-adic frame has
degree $2$; it features two completions of the underlying field $\mathbb{Q}$:
$\mathbb{C}$ and $\mathbb{C}_{q}$. The frames proposed in \textbf{Example
\ref{exa:8}} to make sense of $S_{2;3,5}\left(\mathfrak{z}\right)$
would have degree $3$; the places there are the archimedean place
of $\mathbb{Q}$, and the places of $\mathbb{Q}$ corresponding to
the primes $3$ and $5$.
\end{example}
With this terminology, we can discuss a case where frames of higher
degree naturally arise. In brief, given a Collatz-type map $H:\mathbb{Z}\rightarrow\mathbb{Z}$
satisfying certain mild conditions, one can, following the lines\footnote{The functional equations are those generated by the affine linear
maps contained in $H$.} of \textbf{Example \ref{exa:5}}, construct a function on $\chi_{H}:\mathbb{N}_{0}\rightarrow\mathbb{Q}$
whose domain can be expanded from $\mathbb{N}_{0}$ to $\mathbb{Z}_{p}$
(again, in the manner of \textbf{Example \ref{exa:5}}), where $p$
is a prime number such that $H$ is defined by $p$ distinct ``branches''/''rules''.
In our terminology, \cite{my dissertation} called $H$ \textbf{monogenic
}when $\chi_{H}$ was compatible with a frame of degree $2$, and
said $H$ was \textbf{polygenic }when $\chi_{H}$ was compatible with
a frame of degree $3$ or higher. At the time \cite{my dissertation}
was written, the author had yet to realize that frames of higher degree
could be used to deal with the polygenic case. As such, \cite{my dissertation}
focused the bulk of its attention on monogenic $H$s.

One of the principal discoveries\footnote{The results given below are stated and proven in \cite{first blog post paper}
for the special case of the so-called \textbf{Shortened $qx+1$ map
}$T_{q}:\mathbb{Z}\rightarrow\mathbb{Z}$ defined by:
\begin{equation}
T_{q}\left(n\right)\overset{\textrm{def}}{=}\begin{cases}
\frac{n}{2} & \textrm{if }n\overset{2}{\equiv}0\\
\frac{qn+1}{2} & \textrm{if }n\overset{2}{\equiv}1
\end{cases}
\end{equation}
where $q$ is a fixed odd prime.} of \cite{my dissertation} was a result the author calls the \textbf{Correspondence
Principle} (\textbf{CP}). This consists of the following two statements:
\begin{enumerate}
\item $x\in\mathbb{Z}$ is a \textbf{periodic point}\footnote{That is, there is an integer $n\geq1$ so that $H^{\circ n}\left(x\right)=x$,
where $H^{\circ n}\left(x\right)\overset{\textrm{def}}{=}\underbrace{H\circ\cdots\circ H}_{n\textrm{ times}}$.} of $H$ if and only if there is a $\mathfrak{z}\in\mathbb{Z}_{p}\cap\mathbb{Q}$
so that $\chi_{H}\left(\mathfrak{z}\right)=x$.
\item Let $\mathfrak{z}\in\mathbb{Z}_{p}\backslash\mathbb{Q}$. If $\chi_{H}\left(\mathfrak{z}\right)\in\mathbb{Z}$,
then $\chi_{H}\left(\mathfrak{z}\right)$ is a \textbf{divergent point}\footnote{$x\in\mathbb{Z}$ is a divergent point when the ordinary absolute
value $\left|H^{\circ n}\left(x\right)\right|_{\infty}$ tends to
$\infty$ as $n\rightarrow\infty$. The author conjectures in \cite{my dissertation}
that the converse of (2) is also true: if $x\in\mathbb{Z}$ is a divergent
point of $H$, then there is a $\mathfrak{z}\in\mathbb{Z}_{p}\backslash\mathbb{Q}$
so that $x=\chi_{H}\left(\mathfrak{z}\right)$.}\textbf{ }of $H$.
\end{enumerate}
\cite{my dissertation} only proves the \textbf{CP }for \emph{monogenic}
$H$. However, subsequent work by the author (yet to be published)
shows that the \textbf{CP} holds for \emph{polygenic} $H$ as well,
and in a strikingly frame-theoretic way. Specifically, letting $H$
be polygenic with $p$ branches (where $p$ is a prime number), and
supposing $H$ satisfies certain mild non-degeneracy conditions, and
viewing $\chi_{H}$ as a formal $\mathcal{F}$-series\textemdash one
yet to be equipped with any particular frame\textemdash we have the
following polygenic \textbf{Correspondence Principle}:
\begin{enumerate}
\item An integer $x\in\mathbb{Z}$ is a periodic point of $H$ if and only
if there exists a $\mathfrak{z}\in\mathbb{Z}_{p}\cap\mathbb{Q}$ and
a prime $\ell$ of $\mathbb{Q}$ (possibly the infinite prime) so
that the $\mathcal{F}$-series $\chi_{H}\left(\mathfrak{z}\right)$
converges to $x$ in the $\ell$-adic topology.
\item Let $\mathfrak{z}\in\mathbb{Z}_{p}\backslash\mathbb{Q}$. If there
is a prime $\ell$ of $\mathbb{Q}$ (possibly the infinite prime)
so that the $\mathcal{F}$-series $\chi_{H}\left(\mathfrak{z}\right)$
converges in the $\ell$-adic topology to an element of $\mathbb{Z}$,
then the integer $\chi_{H}\left(\mathfrak{z}\right)$ is a divergent
point of $H$.
\end{enumerate}
Moreover, just like with the monogenic case, we conjecture that the
converse also holds. Specifically, addressing the polygenic and monogenic
cases in one fell swoop, we have:
\begin{conjecture}
Let $H:\mathbb{Z}\rightarrow\mathbb{Z}$ be a $p$-Hydra map (where
$p$ is a finite prime) so that its numen $\chi_{H}:\mathbb{N}_{0}\rightarrow\mathbb{Q}$
can be continued to a function $\chi_{H}\in C\left(\mathcal{F}\right)$,
where $\mathcal{F}$ is a $\mathbb{Q}$-frame $\mathcal{F}$ on $\mathbb{Z}_{p}$
of degree $\geq2$. Then, $x\in\mathbb{Z}$ is a divergent point of
$H$ if and only if there is a $\mathfrak{z}\in\mathbb{Z}_{p}\backslash\mathbb{Q}$
and a prime $\ell$ of $\mathbb{Q}$ (including possibly the infinite
prime) so that $\lim_{n\rightarrow\infty}\chi_{H}\left(\left[\mathfrak{z}\right]_{p^{n}}\right)$
converges to $x$ in the $\ell$-adic topology.
\end{conjecture}
The moral here is that the CP only seems to care if there exists \emph{some
}topology and some $\mathfrak{z}\in\mathbb{Z}_{p}$ at which the $\mathcal{F}$-series
$\chi_{H}\left(\mathfrak{z}\right)$ converges to an element of $\mathbb{Z}$.
As a result, if we only consider the behavior of a function like $\chi_{H}$
with regard to one topology of convergence, we immediately close ourselves
off to valuable information that could be obtained by changing our
topological perspective.

\vphantom{}

At this point, it should (one hopes) feel natural to want to turn
toward \textbf{adèlic }considerations. Indeed, one of the main themes
of this paper\textemdash considering convergence of sequences of rational
numbers with respect to different metric completions of $\mathbb{Q}$\textemdash has
a decidedly adélic smell to it. In case the reader is not familiar
with algebraic number theory and its construction of the \textbf{adèle
ring }of a number field, we briefly and informally review it here\footnote{A more rigorous, comprehensive exposition of the adèles is given in
\cite{adelic analysis}; that resource also details $p$-adic and
adèlic distributions and the Schwartz-Bruhat functions which are used
as test functions. It should be mentioned, however, that, like virtually
all of the extant literature, to the extent analysis is done with
the adèles, it is with functions that accept adèlic inputs and produce
real or complex numbers as outputs. This is the opposite (or, one
might say, \emph{dual}) of what we do here, which is to consider functions
that produce adèles as \emph{outputs}. That being said, there may
be a way to combine these two approaches, and have functions whose
inputs and outputs are both adèlic.}.

Given a number field $K$, a \textbf{$K$-adèle} is a tuple whose
entries are indexed by the ``places'' (a.k.a. prime (ideals)) of
$K$, such that for each place $v$, the $v$th entry of a $K$-adèle
is an element of $K_{v}$ (the metric completion of $K$ with respect
the absolute value induced by $v$). Additionally, because we want
to be able to use the Poisson Summation Formula, we have to stipulate
that the $v$th entry of any given he $K$-adèle is an element of
$\mathcal{O}_{K_{v}}$ (the ring of integers of the completion of
$K$ with respect to $v$) for all but finitely many $v$s. This is
called the \textbf{restricted product condition}.\textbf{ }The adéle
ring of $K$ is then denoted $\mathbb{A}_{K}$. This is a ring under
the operations of component-wise addition and multiplication.

For example, an element of the rational adèles $\mathbb{A}_{\mathbb{Q}}$
is a tuple $\mathbf{x}=\left(x_{\infty},x_{2},x_{3},x_{5},x_{7},\ldots\right)$
such that $x_{\infty}$ is a real number, $x_{p}$ is in $\mathbb{Q}_{p}$
for all primes $p$, with $x_{p}\in\mathbb{Z}_{p}$ for all but finitely
many primes $p$. The restricted product condition can be thought
of as a way of encoding the fundamental theorem of arithmetic\textemdash that
every integer is uniquely factorable as a product of primes.

Given any rational number $x$, $x$ is naturally represented in adèlic
form by the tuple whose every entry is $x$. If $x$ is an integer,
then $x$ is a $p$-adic integer for all primes $p$. However, if
$x$ is a non-integer rational number, then $x$ will be a $p$-adic
integer if and only if $p$ is not a divisor of $x$'s denominator.
Thus, the at most finitely many primes $p$ for which $x$ is not
a $p$-adic integer (i.e., an element of $\mathbb{Q}_{p}\backslash\mathbb{Z}_{p}$)
correspond to the prime divisors of $x$'s denominator. We don't allow
adèles to have non-integer components for infinitely many primes because
that would correspond to a non-zero rational number whose denominator
was divisible by infinitely many primes\textemdash but then the denominator
is infinite, and so our rational number is actually $0$.

The adèle ring of a number field $K$ is eminently compatible with
$K$-frames. It provides us with an accounting tool for keeping track
of all the topological viewpoints we use to study a given $\mathcal{F}$-series,
as the following example shows.
\begin{example}
\label{exa:Let-us-return}Let us return to $S_{7,6}\left(\mathfrak{z}\right)$
from \textbf{Example \ref{exa:2}}. When $\mathfrak{z}\in\mathbb{N}_{0}$,
we can use the formula (\ref{eq:S_p,q at non-negative integer z})
to write:
\begin{equation}
S_{7,6}\left(\mathfrak{z}\right)\overset{\mathbb{C}}{=}\frac{7}{6}\frac{6^{\#_{1}\left(\mathfrak{z}\right)}}{7^{\lambda_{2}\left(\mathfrak{z}\right)}}+\sum_{n=0}^{\lambda_{2}\left(\mathfrak{z}\right)-1}\frac{6^{\#_{1}\left(\left[\mathfrak{z}\right]_{2^{n}}\right)}}{7^{n}},\textrm{ }\forall\mathfrak{z}\in\mathbb{N}_{0}
\end{equation}
This is a rational number, and as such, can be embedded in $\mathbb{Q}_{\ell}$
for $\ell\in\left\{ 2,3,\infty\right\} $. On the other hand, for
each $\mathfrak{z}\in\mathbb{Z}_{2}^{\prime}$, there are two possibilities:

\textbullet{} $\mathfrak{z}$ is rational, in which case $S_{7,6}\left(\mathfrak{z}\right)$
will be a rational number $x$, and $X_{7,6}^{\left(\ell\right)}\left(\mathfrak{z}\right)$
converges to $x$ in $\mathbb{Q}_{\ell}$ for $\ell\in\left\{ 2,3,\infty\right\} $.

\textbullet{} $\mathfrak{z}$ is irrational, in which case, for each
$\ell\in\left\{ 2,3,\infty\right\} $, $S_{7,6}\left(\mathfrak{z}\right)$
will converge to an irrational element of $\mathbb{Q}_{\ell}$. Moreover,
these three irrational numbers will be completely incomparable with
one another; unlike rational numbers, the irrational corresponding
to one $\ell$ does not live in the completion of $\mathbb{Q}$ for
any \emph{other} $\ell$.

Instead of treating the $X_{7,6}^{\left(\ell\right)}$s as completely
different realizations of $S_{7,6}$, let us take a page from the
adèles and combine these realizations into a single object, a function
from $\mathbb{Z}_{2}$ to $\mathbb{R}\times\mathbb{Q}_{2}\times\mathbb{Q}_{3}$.
We do this by way of the rule:
\begin{equation}
X_{7,6;\textrm{tuple}}\left(\mathfrak{z}\right)\overset{\textrm{def}}{=}\left(X_{7,6}^{\left(\infty\right)}\left(\mathfrak{z}\right),X_{7,6}^{\left(2\right)}\left(\mathfrak{z}\right),X_{7,6}^{\left(3\right)}\left(\mathfrak{z}\right)\right),\textrm{ }\forall\mathfrak{z}\in\mathbb{Z}_{2}
\end{equation}
Note that, for any $\mathfrak{z}\in\mathbb{Z}_{2}$, the triple on
the right will converge in the ``product topology'' $\mathbb{R}\times\mathcal{F}_{2,2}\times\mathcal{F}_{2,3}$,
where we sum all three series in the topology of $\mathbb{R}$ if
$\mathfrak{z}\in\mathbb{N}_{0}$ and sum the series in the topologies
of $\mathbb{R}$, $\mathbb{Q}_{2}$, and $\mathbb{Q}_{3}$, respectively,
for $\mathfrak{z}\in\mathbb{Z}_{2}^{\prime}$. But we can go further.
By embedding $\mathbb{R}\times\mathbb{Q}_{2}\times\mathbb{Q}_{3}$
in $\mathbb{A}_{\mathbb{Q}}$, we can realize $S_{7,6}$ as a function
$X_{7,6;\textrm{adèle}}:\mathbb{Z}_{2}\rightarrow\mathbb{A}_{\mathbb{Q}}$.
There are many ways to do this embedding; likely the simplest would
be to define $X_{7,6;\textrm{adèle}}\left(\mathfrak{z}\right)$ as
the $\mathbb{Q}$-adèle:
\begin{equation}
X_{7,6;\textrm{adèle}}\left(\mathfrak{z}\right)\overset{\textrm{def}}{=}\left(X_{7,6;\textrm{adèle}}^{\left(\infty\right)}\left(\mathfrak{z}\right),X_{7,6;\textrm{adèle}}^{\left(2\right)}\left(\mathfrak{z}\right),X_{7,6;\textrm{adèle}}^{\left(3\right)}\left(\mathfrak{z}\right),X_{7,6;\textrm{adèle}}^{\left(5\right)}\left(\mathfrak{z}\right),\ldots\right)
\end{equation}
where, for every $\mathfrak{z}\in\mathbb{Z}_{2}$ and every prime
$\ell$ (including the infinite prime):
\begin{equation}
X_{7,6;\textrm{adèle}}^{\left(\ell\right)}\left(\mathfrak{z}\right)\overset{\textrm{def}}{=}\begin{cases}
X_{7,6}^{\left(\ell\right)}\left(\mathfrak{z}\right) & \textrm{if }\ell\in\left\{ \infty,2,3\right\} \\
X_{7,6}^{\left(\infty\right)}\left(\mathfrak{z}\right) & \textrm{if }\ell\notin\left\{ \infty,2,3\right\} \textrm{ \& }\mathfrak{z}\in\mathbb{Q}\\
0 & \textrm{if }\ell\notin\left\{ \infty,2,3\right\} \textrm{ \& }\mathfrak{z}\notin\mathbb{Q}
\end{cases}
\end{equation}
Alternatively, we might want to make $X_{7,6;\textrm{adèle}}^{\left(\ell\right)}\left(\mathfrak{z}\right)$
equal to the point at infinity for irrational $\mathfrak{z}$ and
$\ell\notin\left\{ \infty,2,3\right\} $. We could do this by adjoining
a point at infinity to the adèles, or, perhaps by defining \textbf{projective
($\mathbb{Q}$-)adèlic space}, denoted $\mathbb{P}^{1}\left(\mathbb{A}_{\mathbb{Q}}\right)$\textemdash the
set of equivalence classes of $\mathbb{A}_{\mathbb{Q}}$ generated
by the relation $\sim$ defined by:
\[
\vec{\mathfrak{x}}\sim\vec{\mathfrak{y}}\Leftrightarrow\exists\lambda\in\mathbb{Q}\backslash\left\{ 0\right\} :\vec{\mathfrak{y}}=\lambda\vec{\mathfrak{x}}
\]
Note that scalar multiplication of $\vec{\mathfrak{x}}\in\mathbb{A}_{\mathbb{Q}}$
by $\lambda\in\mathbb{Q}$ is well-defined because every $\lambda\in\mathbb{Q}$
is uniquely represented adèlically as the $\mathbb{Q}$-adéle whose
$\ell$th component is $\lambda$ for all primes $\ell$. In this
case, the scalar product $\lambda\vec{\mathfrak{x}}$ is the same
as the product of the adèle $\vec{\mathfrak{x}}$ with the adèle representing
$\lambda$.
\end{example}
\begin{rem}
We mention the possibility of projective ($\mathbb{Q}$-)adèlic space
only because of a tantalizing connection between $\chi_{3}$ self-intersections
of curves. As addressed in \cite{first blog post paper} (specifically,
\textbf{Remark 23} of that paper), defining the bijection: $\eta_{2}:\mathbb{Z}_{2}\rightarrow\left[0,1\right]$
by:
\begin{equation}
\eta_{2}\left(\sum_{n=0}^{\infty}j_{n}2^{n}\right)\overset{\textrm{def}}{=}\sum_{n=0}^{\infty}\frac{j_{n}}{2^{n+1}}\label{eq:Definition of Eta_2}
\end{equation}
with inverse $\eta_{2}^{-1}$, the function $C_{3}:\left[0,1\right]\rightarrow\mathbb{Z}_{3}$
defined by $C_{3}\overset{\textrm{def}}{=}\chi_{3}\circ\eta_{2}^{-1}$
then gives us a $3$-adic ``curve''. This is not a continuous function,
but it does satisfy a one-sided continuity condition, with:
\begin{equation}
\lim_{n\rightarrow\infty}C_{3}\left(t_{n}\right)\overset{\mathbb{Q}_{3}}{=}C_{3}\left(t\right)
\end{equation}
for all $t\in\left[0,1\right]$ and all sequences $\left\{ t_{n}\right\} _{n\geq1}$
in $\left[0,1\right]$ converging to $t$ satisfying $0\leq t_{1}<t_{2}<t_{3}<\cdots<t$.
The \textbf{Correspondence Principle }\cite{my dissertation,first blog post paper}
then yields that $x\in\mathbb{Z}$ is a periodic point of $T_{3}$
if and only if $C_{3}\left(t\right)=x$ for some $t\in\left[0,1\right]\cap\mathbb{Q}$,
and that any $t\in\left[0,1\right]\backslash\mathbb{Q}$ which makes
$C_{3}\left(t\right)\in\mathbb{Z}$ then makes $C_{3}\left(t\right)$
a divergent point of $T_{3}$. Conjecturally, the divergent points
of $T_{3}$ (if they exist) are precisely those $x\in\mathbb{Z}$
which are branch points / self-intersection points of $C_{3}\left(t\right)$:
i.e., there is more than one $t\in\left[0,1\right]$ for which $C_{3}\left(t\right)=x$.

If we were to attempt to realize $\chi_{3}$ (and, by pre-composition
with $\eta_{2}^{-1}$, $C_{3}$) as an adèle-valued function in the
manner of \textbf{Example \ref{exa:Let-us-return}}, perhaps the use
of projective adèles would come in handy. Additionally, the definition
of the projective adèles by equivalence relations might end up preserving
arithmetically significant information.
\end{rem}
Things get more complicated when we consider $S_{2;3,5}\left(\mathfrak{z}\right)$.
\begin{example}
To make sense of $S_{2;3,5}\left(\mathfrak{z}\right)$, it helps to
write out the frame explicitly by listing, for each topology $T$,
every $\mathfrak{z}\in\mathbb{Z}_{3}$ at which $S_{2;3,5}\left(\mathfrak{z}\right)$
converges in $T$. That is to say, we will be writing out the pre-images
of $\mathbb{Q}_{\ell}$ under $\mathcal{F}$, for $\ell\in\left\{ 3,5,\infty\right\} $.

\textbullet{} (Convergence in $\mathbb{R}$): $S_{2;3,5}\left(\mathfrak{z}\right)$
clearly converges in $\mathbb{R}$ if $\mathfrak{z}\in\mathbb{N}_{0}$,
where $S_{2;3,5}\left(\mathfrak{z}\right)$ can be written in terms
of a geometric series with common ratio $1/2$. Note also that $S_{2;3,5}\left(\mathfrak{z}\right)$
will converge in $\mathbb{R}$ for any $\mathfrak{z}\in\mathbb{Z}_{3}$
whose non-zero $3$-adic digits are separated by sufficiently many
$0$s as to guarantee that the partial sums of $S_{2;3,5}\left(\mathfrak{z}\right)$
converge in $\mathbb{R}$. We can do this explicitly for $\mathfrak{z}$s
with (eventually) periodic digits sequences. Here, for any finite
prime $\ell$, we define $B_{\ell}:\mathbb{N}_{0}\rightarrow\mathbb{Q}$
by:
\begin{equation}
B_{\ell}\left(n\right)\overset{\textrm{def}}{=}\begin{cases}
0 & \textrm{if }n=0\\
\frac{n}{1-\ell^{\lambda_{\ell}\left(n\right)}} & \textrm{if }n\geq1
\end{cases}\label{eq:Definition of B_p}
\end{equation}
and then derive the formula:
\begin{equation}
\#_{\ell:j}\left(\left[B_{\ell}\left(n\right)\right]_{\ell^{k}}\right)=\left\lfloor \frac{k}{\lambda_{\ell}\left(n\right)}\right\rfloor \#_{\ell:j}\left(n\right)+\#_{\ell:j}\left(\left[n\right]_{\ell^{\left[k\right]_{\lambda_{\ell}\left(n\right)}}}\right)\label{eq:number of js in B_p mod p to the k}
\end{equation}
Then example, for $\mathfrak{z}$ with the digit sequence: 
\begin{equation}
\mathfrak{z}=000120001200012000\ldots=1\cdot3^{3}+2\cdot3^{4}+1\cdot3^{8}+2\cdot3^{9}+1\cdot3^{13}+2\cdot3^{14}+\cdots
\end{equation}
we get:
\begin{align*}
S_{2;3,5}\left(\overline{00012}\ldots\right) & =S_{2;3,5}\left(B_{3}\left(63\right)\right)\\
 & =\sum_{k=0}^{\infty}\frac{3^{\#_{3:1}\left(\left[63\right]_{3^{k}}\right)}5^{\#_{3:2}\left(\left[63\right]_{3^{k}}\right)}}{2^{k}}\\
\left(\lambda_{3}\left(63\right)=5\right); & =\sum_{k=0}^{\infty}\frac{3^{\left\lfloor \frac{k}{5}\right\rfloor \#_{3:1}\left(63\right)+\#_{3:1}\left(\left[63\right]_{3^{\left[k\right]_{5}}}\right)}5^{\left\lfloor \frac{k}{5}\right\rfloor \#_{3:2}\left(63\right)+\#_{3:2}\left(\left[63\right]_{3^{\left[k\right]_{5}}}\right)}}{2^{k}}\\
\left(\textrm{split }k\textrm{ mod }5\right); & =\sum_{n=0}^{\infty}\sum_{k=0}^{4}\frac{3^{n\#_{3:1}\left(63\right)+\#_{3:1}\left(\left[63\right]_{3^{k}}\right)}5^{n\#_{3:2}\left(63\right)+\#_{3:2}\left(\left[63\right]_{3^{k}}\right)}}{2^{5n+k}}\\
 & =\left(\sum_{n=0}^{\infty}\left(\frac{15}{32}\right)^{n}\right)\sum_{k=0}^{4}\frac{3^{\#_{3:1}\left(\left[63\right]_{3^{k}}\right)}5^{\#_{3:2}\left(\left[63\right]_{3^{k}}\right)}}{2^{k}}\\
 & \overset{!}{=}\frac{1}{1-\frac{15}{32}}\left(\frac{3^{0}\cdot5^{0}}{1}+\frac{3^{0}\cdot5^{0}}{2}+\frac{3^{0}\cdot5^{0}}{4}+\frac{3^{1}\cdot5^{0}}{8}+\frac{3^{1}\cdot5^{1}}{16}\right)\\
 & =\frac{1}{1-\frac{15}{32}}\cdot\frac{49}{16}\\
 & =\frac{98}{17}
\end{align*}
In the step marked ($!$), observe that the geometric series:
\begin{equation}
\sum_{n=0}^{\infty}\left(\frac{15}{32}\right)^{n}=\frac{1}{1-\frac{15}{32}}
\end{equation}
converges simultaneously in $\mathbb{R}$, $\mathbb{Z}_{3}$, and
$\mathbb{Z}_{5}$. A bit of thought shows that the structure of the
set of $\mathfrak{z}\in\mathbb{Z}_{3}\cap\mathbb{Q}$ for which $S_{2;3,5}\left(\mathfrak{z}\right)$
converges in $\mathbb{R}$ is related to the value of linear logarithmic
form:
\begin{equation}
a\ln3+b\ln5-c\ln2
\end{equation}
in $\mathbb{R}$, where $a,b,c\in\mathbb{Z}$; in particular, for
those $a,b,c$ for which the linear logarithmic form is negative.
In this way, we see that the frames compatible with a given $\mathcal{F}$-series
of the form $S_{d;q_{0},\ldots,q_{p-1}}$ relate to the behavior of
the linear logarithmic form:
\begin{equation}
-b\ln d+\sum_{k=0}^{p-1}a_{k}\ln q_{k}
\end{equation}
Returning to $S_{2;3,5}$, we will let $U_{\infty}$ denote the set
of all $\mathfrak{z}\in\mathbb{Z}_{3}$ for which $S_{2;3,5}\left(\mathfrak{z}\right)$
converges in $\mathbb{R}$. This Estimates on the absolute values
of linear logarithmic forms is a mainstay of transcendental number
theory, and it seems worth investigating if frames and $\mathcal{F}$-series
may be of use in that notoriously difficult area of number theory.

\textbullet{} (Convergence in $\mathbb{Z}_{3}$ and/or $\mathbb{Z}_{5}$)
Whereas the convergence of $S_{2;3,5}\left(\mathfrak{z}\right)$ in
$\mathbb{R}$ is subtle, its non-archimedean convergence is quite
simple. $U_{3}$ and $U_{5}$\textemdash the sets of all $\mathfrak{z}\in\mathbb{Z}_{3}$
at which $S_{2;3,5}\left(\mathfrak{z}\right)$ converges $3$-adically
and $5$-adically, respectively\textemdash are precisely the sets
of all $\mathfrak{z}\in\mathbb{Z}_{3}$ with infinitely many $1$s
digits and infinitely many $2$s digits, respectively, in their $3$-adic
digits.

With this analysis, we can then construct a frame $\mathcal{F}$ for
$S_{2;3,5}$ which captures the totality of $S_{2;3,5}$'s convergence
behavior.

\begin{equation}
\mathcal{F}\left(\mathfrak{z}\right)\overset{\textrm{def}}{=}\prod_{\ell\in\left\{ 3,5,\infty\right\} :\mathfrak{z}\in U_{\ell}}\mathbb{Q}_{\ell}\label{eq:Frame for S_2;3,5}
\end{equation}
We can realize $S_{2;3,5}$ as a function $X_{2;3,5;\textrm{adèle}}:\mathbb{Z}_{3}\rightarrow\mathbb{A}_{\mathbb{Q}}\cup\left\{ \hat{\infty}\right\} $
($\hat{\infty}$ is the point at infinity) by defining for all $\ell\in\left\{ 3,5,\infty\right\} $:

\begin{equation}
X_{2;3,5}^{\left(\ell\right)}\left(\mathfrak{z}\right)\overset{\textrm{def}}{=}\begin{cases}
\textrm{Formal sum of }S_{2;3,5}\left(\mathfrak{z}\right) & \textrm{if }\mathfrak{z}\in\mathbb{Z}_{3}\cap\mathbb{Q}\\
\textrm{Sum of }S_{2;3,5}\left(\mathfrak{z}\right)\textrm{ in }\mathbb{Q}_{\ell} & \textrm{if }\mathfrak{z}\in U_{\ell}\backslash\mathbb{Q}\\
\hat{\infty} & \textrm{if }\mathfrak{z}\in\mathbb{Z}_{3}\backslash\left(\mathbb{Q}\cup U_{\ell}\right)
\end{cases}\label{eq:Definition of X_2;3,5's pth component}
\end{equation}
where the ``formal sum'' of $S_{2;3,5}\left(\mathfrak{z}\right)$
simply means writing $S_{2;3,5}\left(\mathfrak{z}\right)$ as geometric
series of the form:
\begin{equation}
A+B\sum_{n=0}^{\infty}r^{n}
\end{equation}
for constants $A,B,r$ depending on $\mathfrak{z}$ (note: this can
be done only because $\mathfrak{z}\in\mathbb{Z}_{3}\cap\mathbb{Q}$)
and formally evaluating the series as:
\begin{equation}
A+\frac{B}{1-r}
\end{equation}
The reason why this is justified is because, for any $\mathfrak{z}\in\mathbb{Z}_{3}\cap\mathbb{Q}$,
there exists at least one $\ell\in\left\{ 3,5,\infty\right\} $ so
that $\mathfrak{z}\in U_{\ell}$, and hence, so that $S_{2;3,5}\left(\mathfrak{z}\right)$
is rigorously convergent. Using the ``formal sum''\footnote{It would be interesting to see if there was an \emph{algebraic} way
to define this formal summation process, perhaps by playing around
with quotients of rings of formal power series?} to define (\ref{eq:Definition of X_2;3,5's pth component}) is just
a short-hand. A more precise way of defining that part of $X_{2;3,5}^{\left(\ell\right)}$
would be to say that for any $\mathfrak{z}\in\mathbb{Z}_{3}\cap\mathbb{Q}$,
the number $X_{2;3,5}^{\left(\ell\right)}\left(\mathfrak{z}\right)$
is the rational number produced by summing $S_{2;3,5}\left(\mathfrak{z}\right)$
in any topology in which it converges.

We then define $X_{2;3,5;\textrm{adèle}}:\mathbb{Z}_{3}\rightarrow\mathbb{A}_{\mathbb{Q}}\cup\left\{ \hat{\infty}\right\} $
as the (projective?) $\mathbb{Q}$-adèle whose $\ell$th entry, for
any prime $\ell$, is:
\begin{equation}
\begin{cases}
X_{2;3,5}^{\left(\ell\right)}\left(\mathfrak{z}\right) & \textrm{if }\ell\in\left\{ 3,5,\infty\right\} \\
\textrm{Formal sum of }S_{2;3,5}\left(\mathfrak{z}\right) & \textrm{if }\ell\notin\left\{ 3,5,\infty\right\} \textrm{ and }\mathfrak{z}\in\mathbb{Z}_{3}\cap\mathbb{Q}\\
\hat{\infty} & \textrm{if }\ell\notin\left\{ 3,5,\infty\right\} \textrm{ and }\mathfrak{z}\notin\mathbb{Z}_{3}\backslash\mathbb{Q}
\end{cases}
\end{equation}
Note that because the formal sum of $S_{2;3,5}\left(\mathfrak{z}\right)$
at any $\mathfrak{z}\in\mathbb{Z}_{3}\cap\mathbb{Q}$ is an element
of $\mathbb{Q}$, this definition ensures that whenever $\mathfrak{z}\in\mathbb{Z}_{3}\cap\mathbb{Q}$,
the rational number produced by $S_{2;3,5}\left(\mathfrak{z}\right)$
is then properly represented adèlically; $X_{2;3,5;\textrm{adèle}}\left(\mathfrak{z}\right)$
is then the $\mathbb{Q}$-adèle whose every entry is the rational
number $S_{2;3,5}\left(\mathfrak{z}\right)$.
\end{example}
We now give the definition of a frame of arbitrary degree and dimension.
First, however, a technical definition:
\begin{defn}
Let $K$ be a field. We say a Banach space $\mathcal{B}$ is \textbf{$K$-completed}
if the underlying field of $\mathcal{B}$ is the completion of $K$
with respect to an absolute value. We write $\mathcal{B}_{K}$ to
denote the set of all $K$-completed Banach spaces.
\end{defn}
\begin{rem}
\emph{We allow for our completions to be with respect to the trivial
absolute value!} Consequently, we can have Banach spaces over fields
like $\mathbb{Q}$ or $\overline{\mathbb{Q}}$ and any intermediary
extension thereof. In the case where $\mathcal{B}$'s field is equipped
with the trivial absolute value, a sequence$\left\{ x_{n}\right\} _{n\geq0}$
in $\mathcal{B}$ converges if and only if $x_{n+1}=x_{n}$ for all
sufficiently large $n$.
\end{rem}
\begin{defn}
Given a topological space $X$, a field $K$, and an integer $d\geq1$,
a \textbf{$d$-dimensional} \textbf{$K$}-\textbf{frame} $\mathcal{F}$
on $X$ is a map $\mathcal{F}:X\rightarrow\mathcal{B}_{K}$ so that,
for each $x\in X$, the Banach space $\mathcal{F}\left(x\right)$
is a Cartesian product of the form:
\begin{equation}
\mathcal{F}\left(x\right)=\prod_{n=1}^{\deg_{x}\mathcal{F}}K_{v_{x,n}}^{d}
\end{equation}
where $\deg_{x}\mathcal{F}$ is an integer $\geq1$ depending on $x$
called the \textbf{degree of $\mathcal{F}$ at $x$}, and where $v_{x,1},\ldots,v_{x,\deg_{x}\mathcal{F}}$
are distinct places of $K$, depending on both $x$ and $n$. We write
$V_{\mathcal{F}}\left(x\right)$ to denote the set $\left\{ v_{x,1},\ldots,v_{x,\deg_{x}\mathcal{F}}\right\} $,
which we call \textbf{the set of places of $\mathcal{F}$ at $x$}.\textbf{
}We then define $V_{\mathcal{F}}$\textemdash \textbf{the set of places
of $\mathcal{F}$}\textemdash by:
\begin{equation}
V_{\mathcal{F}}\overset{\textrm{def}}{=}\bigcup_{x\in X}V_{\mathcal{F}}\left(x\right)
\end{equation}
and define \textbf{degree of $\mathcal{F}$}, denoted $\deg\mathcal{F}$,
as the cardinality\footnote{When $\mathcal{F}$ is $1$-dimensional, $V_{\mathcal{F}}\left(x\right)$
is a set with one element, so this definition generalizes the one
given in \textbf{Definition \ref{def:degree of a 1-dim frame}}.} of $V_{\mathcal{F}}$.

We give each $K_{v_{x,n}}^{d}$ the topology of the $\ell^{\infty}$-norm
induced by $K_{v_{x,n}}$'s absolute value ($\left|\cdot\right|_{v_{x,n}}$):
\begin{equation}
\left\Vert \mathbf{a}\right\Vert _{K_{v_{x,n}}}\overset{\textrm{def}}{=}\max\left\{ \left|a_{1}\right|_{v_{x,n}},\ldots,\left|a_{d}\right|_{v_{x,n}}\right\} ,\textrm{ }\forall\mathbf{a}=\left(a_{1},\ldots,a_{d}\right)\in K_{v_{x,n}}^{d}
\end{equation}
That is, an element of $K_{v_{x,n}}^{d}$ is a $d$-tuple $\mathbf{a}$
whose entries belong to the field $K_{v_{x,n}}$, and the norm of
$\mathbf{a}$ is the maximum of the $v_{x,n}$-adic absolute values
of its entries. We give $\mathcal{F}\left(x\right)$ the topology
induced by the norm:
\begin{equation}
\left\Vert \vec{\mathbf{a}}\right\Vert _{\mathcal{F}\left(x\right)}\overset{\textrm{def}}{=}\max\left\{ \left\Vert \mathbf{a}_{1}\right\Vert _{K_{v_{x,1}}},\ldots,\left\Vert \mathbf{a}_{\deg_{x}\mathcal{F}}\right\Vert _{K_{v_{x,\deg_{x}\mathcal{F}}}}\right\} 
\end{equation}
Here, an element of $\mathcal{F}\left(x\right)$ is a $\deg_{x}\mathcal{F}$-tuple
$\vec{\mathbf{a}}=\left(\mathbf{a}_{1},\ldots,\mathbf{a}_{\deg_{x}\mathcal{F}}\right)$
(or, if the reader prefers, a $\left(\deg_{x}\mathcal{F}\right)\times d$
matrix) whose $m$th entry is the tuple $\mathbf{a}_{m}=\left(a_{m,1},\ldots,a_{m,d}\right)\in K_{v_{x,m}}^{d}$.
The norm of an $\vec{\mathbf{a}}\in\mathcal{F}\left(x\right)$ outputs
the maximum of the norms of the tuples $\mathbf{a}_{1},\ldots,\mathbf{a}_{\deg_{x}\mathcal{F}}$.
This then makes $\mathcal{F}\left(x\right)$ into a Banach space.

The \textbf{image }of $\mathcal{F}$, denoted $I\left(\mathcal{F}\right)$,
is defined as the union of the $\mathcal{F}\left(x\right)$, just
as in the $1$-dimensional case:
\begin{equation}
I\left(\mathcal{F}\right)\overset{\textrm{def}}{=}\bigcup_{x\in X}\mathcal{F}\left(x\right)\label{eq:Image of a frame, multi-dimensional}
\end{equation}
Likewise, we say a function $f:X\rightarrow I\left(\mathcal{F}\right)$
is \textbf{compatible with $\mathcal{F}$} if $f\left(x\right)\in\mathcal{F}\left(x\right)$
for all $x\in X$, and then write $C\left(\mathcal{F}\right)$ to
denote \textbf{the set of all $\mathcal{F}$-compatible functions}.
\end{defn}
\begin{prop}
Let $\mathcal{F}$ be a $K$-frame, and let $\mathbb{F}$ be $K$
or any subfield thereof. Then, $C\left(\mathcal{F}\right)$ is a unital
algebra over $\mathbb{F}$, with point-wise addition of functions,
scalar multiplication of functions by elements of $\mathbb{F}$, and
point-wise multiplication of functions as the addition, scalar multiplication,
and algebra multiplication operations of $C\left(\mathcal{F}\right)$,
respectively.
\end{prop}
Proof: Just do the multi-dimensional analogue of the proof of \textbf{Proposition
\ref{prop:5}}.

Q.E.D.

\vphantom{}Finally, we conclude with the definition of convergence
with respect to a frame of arbitrary dimension.
\begin{defn}
Let $K$ be a field, and let $\mathcal{F}$ be a $d$-dimensional
$K$-frame . Given a function $F\in C\left(\mathcal{F}\right)$, we
say a sequence $\left\{ f_{n}\right\} _{n\geq0}$ of $\mathcal{F}$-compatible
functions \textbf{converges to $F$ with respect to $\mathcal{F}$
}if, for every $x\in X$, $f_{n}\left(x\right)$ converges to $F\left(x\right)$
in the topology of $\mathcal{F}\left(x\right)$. We denote this convergence
by writing:
\begin{equation}
F\left(x\right)\overset{\mathcal{F}}{=}\lim_{n\rightarrow\infty}f_{n}\left(x\right)
\end{equation}
Consequently, in a frame of dimension $d$, compatible functions output
$d$-tuples, and convergence with respect to the frame means that,
at each $x\in X$, the sequence of $d$-tuples $f_{n}\left(x\right)$
converge to the $d$-tuple $F\left(x\right)$ in the topology of $\mathcal{F}\left(x\right)$.
\end{defn}

\section*{Acknowledgements}

As mentioned, this paper presents and builds upon the author's PhD
(Mathematics) dissertation \cite{my dissertation} done at the University
of Southern California under the supervision of Professors Sheldon
Kamienny and Nicolai Haydn. Thanks must also be given to the author's
friends and family, as well as Jeffery Lagarias, Steven J. Miller,
Alex Kontorovich, Andrei Khrennikov, K.R. Matthews, Susan Montgomery,
Amy Young and all the helpful staff of the USC Mathematics Department;
all the kindly strangers became and acquainted with along the way,
and probably lots of other people, too.

\vphantom{}

\end{document}